%% file: bshthesis.tex
\documentclass[11pt]{ucthesis}
\usepackage{amssymb,amsmath,amsthm, amscd}
\usepackage{verbatim}
\usepackage{mathrsfs}
\usepackage[small,nohug,heads=littlevee]{diagrams}
\diagramstyle[labelstyle=\scriptstyle]

\def\dsp{\def\baselinestretch{2.0}\large\normalsize}
\dsp
\usepackage{simplemargins}
    \setrightmargin{1.2 in}

\usepackage{color}

\DeclareMathOperator{\Spec}{Spec} \DeclareMathOperator{\Aut}{Aut}
\DeclareMathOperator{\PGL}{PGL} \DeclareMathOperator{\GL}{GL}
\DeclareMathOperator{\PSL}{PSL} 
\newcommand{\Z}{{\mathbb Z}}
\newcommand{\R}{{\mathbb R}}
\newcommand{\C}{{\mathbb C}}
\newcommand{\F}{{\mathbb F}}

\newcommand{\Q}{{\mathbb Q}}

\DeclareMathOperator{\SL}{SL}\DeclareMathOperator{\GU}{GU}\DeclareMathOperator{\SU}{SU}\DeclareMathOperator{\PGU}{PGU}
\DeclareMathOperator{\PSU}{PSU} \DeclareMathOperator{\Gal}{Gal}
\DeclareMathOperator{\Stab}{Stab} \DeclareMathOperator{\Deg}{deg}
\newcommand{\G}{{\mathfrak G}}

\newcommand{\el}[1]{{}^{#1\negthinspace}}

    \newtheorem{thm}{Theorem}
    \numberwithin{thm}{section}
    \newtheorem{cor}[thm]{Corollary} 
    \newtheorem{lem}[thm]{Lemma}
    \newtheorem{prop}[thm]{Proposition}
    
    \newtheorem{notation}[thm]{Notation}

    \theoremstyle{definition}
    \newtheorem{df}[thm]{Definition}
    \theoremstyle{remark}
    \newtheorem{rem}[thm]{Remark}
\newcommand{\la}{\langle}\newcommand{\ra}{\rangle}

    \newcommand{\ol}{\overline}
    \newcommand{\mf}{\mathfrak}

\begin{document}

\title{Fields of Moduli and Fields of Definition of Curves}
\author{Bonnie Sakura Huggins}
\degreesemester{Fall} \degreeyear{2005} \degree{Doctor of
Philosophy} \chair{Professor Bjorn Poonen}
\othermembers{Professor Kenneth Ribet \\
  Professor Michael Jansson }
\numberofmembers{3} \prevdegrees{B.S. (University of California, Los
Angeles) 1998} \field{Mathematics} \campus{Berkeley}
    \maketitle
    \approvalpage
    \copyrightpage
     \begin{abstract}
     The field of moduli $K_X$ of a curve $X$ over a field $K$ is the intersection
     over all fields of definition of $X_{\ol{K}}$, where $X_{\ol{K}}$ is the base
     extension of $X$ to a curve over an algebraic closure $\ol{K}$ of $K$.
       It has the property that
     for all $\sigma\in\Aut(\ol{K})$, we have $X_{\ol{K}}\cong\el{\sigma}{X_{\ol{K}}}$ if and only if $\sigma|_{K_X}$ is equal
     to the identity.  In this thesis we determine conditions that guarantee that a
     hyperelliptic or plane curve over a field of characteristic not
     equal to $2$ can be defined over its field of moduli.  We also
     give new examples of curves not definable over their fields of
     moduli.

     In Chapter~1, we define the notion of ``field of moduli," we
     compare our definition with definitions given by others, and we show in what ways they
      are equivalent.

     In Chapter~2, we list all of the finite subgroups of the two
     and three dimensional projective general linear groups.  We
     will use these classifications to prove our main results.

     In Chapter~3, we discuss isomorphisms of plane and
     hyperelliptic curves.  We will later use the classification of
     automorphism groups of hyperelliptic and plane curves to
     determine whether a plane curve or a
     hyperelliptic curve can be defined over its field of moduli.

     In Chapter~4, we give our main result in the case of
     hyperelliptic curves.  We show that hyperelliptic curves with
     certain automorphism groups can always be defined over their fields of
     moduli.

     In Chapter~5, we list examples of hyperelliptic curves  not
     definable over their fields of moduli.

     In Chapter~6, we give our main result in the case of plane
     curves.  We show that plane curves with certain automorphism
     groups can always be defined over their fields of moduli.

     In Chapter~7, we give new examples of plane curves not
     definable over their fields of moduli.

    \abstractsignature
    \end{abstract}
    \begin{frontmatter}
\begin{dedication}\null\vfil{\large\begin{center}To my family
\end{center}}\vfil\null\end{dedication}

\tableofcontents 
\begin{acknowledgements} I would like to thank my advisor Bjorn
Poonen to whom I am very grateful
  for many helpful suggestions, advice, and
inspiration. I would also like to thank my friends Rajan Mehta and
Bj{\o}rn Kjos-Hanssen for making my stay in Berkeley pleasant and
fun.

\end{acknowledgements}
\addtocontents{toc}{\contentsline {chapter}{Introduction}{\thepage}}
\chapter*{Introduction}
    \input{intro}
\end{frontmatter}

    \chapter{Fields of moduli}

\input{ch1}
    \chapter{Finite subgroups of projective general linear groups}

\input{ch2}
    \chapter{Isomorphisms of hyperelliptic and plane curves}
    \input{ch3}
    \chapter{Hyperelliptic curves definable over their fields of
    moduli}
    \input{ch4}
    \chapter{Hyperelliptic curves not definable over their fields of
    moduli}
    \input{ch5}
    \chapter{Plane curves definable over their fields of
    moduli}
    \input{ch6}
    \chapter{Plane curves not definable over their fields of moduli}

\input{ch7}

\bibliographystyle{amsplain}
\bibliography{fom}
\end{document}

%% file: intro.tex
The field of moduli $K_X$ of a curve $X$ over a field $K$ is the
intersection
     over all fields of definition of $X_{\ol{K}}$, where $X_{\ol{K}}$ is the base
     extension of $X$ to a curve over an algebraic closure $\ol{K}$ of $K$. It has the property that for all
$\sigma\in\Aut(\ol{K})$, we have
$X_{\ol{K}}\cong\el{\sigma}{X_{\ol{K}}}$ if and only if
$\sigma|_{K_X}$ is the identity.  It is well known that if the genus
$g$ of a curve $X$ is $0$ or $1$, then it can be defined over its
field of moduli. It also can be deduced from~Theorem~1
of~\cite{Weil:1}, that a curve with no nontrivial automorphisms can
be defined over its field of moduli. However, if $g>1$ and if the
automorphism group of $X$ is nontrivial, the curve $X$ may not be
definable over its field of moduli.
 The first examples of curves not definable over their fields of moduli were given by Shimura on
page~177 of~\cite{Sh:1}.  These curves are hyperelliptic $\C$-curves
with two automorphisms.   In this thesis, we study the definability
of plane and hyperelliptic curves, over fields of characteristic not
equal to $2$, over their fields of moduli.

In the first chapter, we relate alternate definitions of ``fields of
moduli" of curves to one another, give our own definition of ``field
of moduli," and show in what sense it is equivalent to the others.
We also briefly mention the relationship between the field of moduli
of a curve of genus $g$ and the residue field at the corresponding
point on the moduli space of curves of genus $g$.

The definability of a curve over its field of moduli depends heavily
on the structure of its group of automorphisms.  To understand the
automorphism groups and isomorphisms of plane and hyperelliptic
curves we need to understand the structure of the finite subgroups
of the $2$ and $3$-dimensional projective general linear groups.  In
Chapter~2, we list these groups and study some of their properties.
In Chapter~3, we discuss the automorphism groups and isomorphisms of
hyperelliptic and plane curves.

 Recall that a curve of genus $2$ is hyperelliptic.  In~\cite{Cardona:1} it is shown that a
curve of genus $2$ over a
 field of characteristic different from $2$, can
be defined over its field of moduli if its automorphism group has
more than two elements.  By Theorem~22 of~\cite{Cardona:2}, any
curve of genus $2$ over a field of characteristic $2$, can be
defined over its field of moduli.  In Section~1 of~\cite{Shaska:1},
it is conjectured that a hyperelliptic curve over a field of
characteristic $0$ is definable over its field of moduli if its
automorphism group has more than two elements.

In Chapter~4, we prove our main theorem for hyperelliptic curves:
\emph{Let $X$ be a hyperelliptic curve over a field of $K$
characteristic not equal to $2$ and let $\iota$ be the hyperelliptic
involution of $X$.  Then $X$ is definable over its field of moduli
of $\Aut(X)/\la\iota\ra$
 is not cyclic or if $\Aut(X)/\la\iota\ra$ is cyclic of order equal to the
 characteristic of $K$.}  Our main result for hyperelliptic curves relies heavily on a result
of D{\`e}bes and Emsalem from~\cite{Debes:1}, that states that a
curve $X$ can be defined over its field of moduli $K$ if a certain
$K$-model of the curve $X/\Aut(X)$ has a $K$-rational point.

The authors of~\cite{Bujalance} have attempted to classify all
hyperelliptic curves over $\C$ with fields of moduli $\R$ relative
to $\C/\R$ but not definable over $\R$.  Due to some errors in their
paper, some curves are missing from their list and many curves on
their list are, in fact, definable over $\R$. In Chapter~5, we give
the complete list of hyperelliptic $\C$-curves, up to isomorphism,
not definable over their fields of moduli relative to $\C/\R$.
 Every curve $X$ in the list has $\Aut(X)/\la\iota\ra$ cyclic of
  order $n$ for some $n\ge 1$.  Evidently the
results of~\cite{Bujalance}, were not known to the author
of~\cite{Shaska:1} at the time he made his conjecture.  Nor were we
aware of them as we had independently constructed similar examples,
not included in the list of~\cite{Bujalance}, before learning of
their results.

In Chapter~6, we prove our main result for plane curves:  \emph{Let
$X$ be a smooth plane curve over a field $K$ of characteristic $p$
where $p=0$ or $p>2$.  Let $F$ be an algebraic closure of $K$. Then
$X$ is definable over its field of moduli if $\Aut(X)$ is not
$\PGL_3(F)$-conjugate to a diagonal subgroup of $\PGL_3(F)$,
$\G_{18}$, $\G_{36}$, or a semi-direct product of a diagonal group
an a $p$-group.}  See Lemma~\ref{gps} in Chapter~2, for a detailed
description of these groups.

Lastly, in Chapter~7, we give new examples of plane curves not
definable over their fields of moduli.  These curves have diagonal
automorphism groups, automorphism groups given by $\G_{18}$, and
automorphism groups given by $\G_{36}$.

%% file: ch1.tex
The notion of fields of moduli for polarized varieties was first
introduced by Matsusaka in~\cite{Matsusaka} and was later defined by
Shimura in~\cite{Sh:2} for polarized abelian varieties and polarized
abelian varieties with further structures. Their definitions, for
the field of moduli of a polarized (abelian) variety $V$,  depend
heavily on the construction of certain subvarieties of projective
space whose closed points consist of Chow points of certain
projective embeddings of $V$. Both authors define this notion in the
case of characteristic zero and of positive characteristic. Koizumi
in~\cite{Koizumi:1} gives a more general definition of field of
moduli for geometric objects satisfying certain conditions. He
removed artificial assumptions previously thought necessary for the
proof of the existence of fields of moduli for polarized varieties
in positive characteristic. In the case of polarized varieties,
Koizumi's definition agrees with the definitions given by Matsusaka
and Shimura.

  We will focus our attention on the field of moduli
of a curve.  We give a definition of field of moduli that agrees
with the definition given by Koizumi.
\section{Fields of definition}

\begin{df} Let $K$ be a field.  A \emph{variety} over $K$
($K$-variety) is an integral separated scheme of finite type over
$\Spec K$.
\end{df}
\begin{rem} We will sometimes speak of a variety $X$ without mention of a base
field.  In this case it should be understood that $X$ is a
$K$-variety for some field $K$.
\end{rem}

\begin{notation} Let $K$ be a field, let $X$ be a $K$-variety, and let $F$ be an extension
field of $K$.  Let $X_F$ denote the base extension $X\times_{\Spec
K} \Spec F$.
\end{notation}

\begin{df}Let $K\subseteq F\subseteq \ol F$ be fields where $\ol F$ is an
algebraic closure of $F$.  Let $X$ be an $F$-variety.  Then $X$ is
\emph{defined} over $K$ if and only if there is a $K$-variety $X'$
such that $X'_F$ is isomorphic (as an $F$-variety) to $X$.  We say
that $K$ is a \emph{field of definition} of $X$.  We say that $X$ is
\emph{definable} over $K$ if there is a $K$-variety $X'$ such that
$X'_{\ol{F}}$ is isomorphic to $X_{\ol{F}}$.
\end{df}

\begin{df} Let $K\subseteq L\subseteq F$ be fields and let $X$ be a geometrically
integral $K$-variety.  Let $Y$ be a subvariety of $X_F$. We say that
$Y$ is \emph{defined} over $L$ as a subvariety of $X_F$ if there
exists an $L$-subvariety $Y'$ of $X_L$ such that $Y'_F\cong Y$ as
subvarieties of $X_F$. In this case, we say that $L$ is a
\emph{field of definition} of $Y$ as a subvariety of $X_F$.  We will
write ``$L$ is of field of definition of $Y\subseteq X_F$" to
indicate that $L$ is a field of definition of $Y$ as a subvariety of
$X_F$.
\end{df}
\begin{lem}\label{mindef}Let $K\subseteq F$ be fields, let $X$
be a geometrically integral $K$-variety, and let $Y$ be a closed
subvariety of $X_F$. Then there is a unique field of definition $L$
of $Y\subseteq X_F$ with $K\subseteq L\subseteq F$ such that for any
field of definition $L'$ of $Y\subseteq X_F$ with $K\subseteq
L'\subseteq F$ we have $L\subseteq L'$. We call $L$ the minimum
field of definition of $Y$ as a subvariety of $X_F$.
\end{lem}
\begin{proof} See Proposition~3.11 in~\cite{Vojta}.
\end{proof}

\begin{rem} Let $F$ be a field, let $X$ be a
 projective $F$-variety, and let
$X\to\mathbb{P}^n_F$ be an embedding. Let $\mathscr{I}_X$ be the
ideal sheaf of $X$ in $\mathbb{P}^n_F$ and let $I\subset
F[X_0,\ldots, X_n]$ be the corresponding ideal. Then the minimum
field of definition of $X\subseteq \mathbb{P}^n_F$ is the smallest
field $K$ contained in $F$ such that $I$ can be generated by
elements in $K[X_0,\ldots, X_n]$.
\end{rem}

\section{The field of moduli of a curve}
\begin{df} Let $K$ be a field.  A \emph{curve} over $K$ is a smooth,
projective, geometrically integral $K$-variety of dimension $1$.
\end{df}
\begin{rem} We will sometimes, especially in the case of plane curves, say ``smooth curve" even though it is
redundant.
\end{rem}

\begin{notation}Let  $X$ be a curve over a field $K$.  For $\sigma\in\Aut(K)$, the curve
$\el{\sigma}{X}$ is the base extension $X\times_K K$ of
 $X$ by the morphism $\Spec K\xrightarrow{\Spec\sigma}\Spec K$.
 \end{notation}

 \begin{df} Let $X$ be a curve over a field $K$.  Let $\ol{K}$ be an algebraic closure
 of $K$.  The \emph{field of moduli} $K_X$ of
 $X$ is the intersection over all fields of definition of $X_{\ol{K}}$.
 \end{df}

  Let $X$ be a curve over a field $K$, let $\ol{K}$ be an algebraic closure of $K$,
   and let $K_X$ be the field of moduli of $X$.
  Many define the field of moduli of a  curve $X$ as the subfield $\ol{K}^H$ of
   $\ol{K}$ fixed by  \[H:=\{\sigma\in\Aut(\ol{K})\mid X_{\ol{K}}\cong\el{\sigma}{X_{\ol{K}}}\}.\]
  Theorem~\ref{koizumi} in Section~\ref{structures} shows
 that, in fact, $\ol{K}^H$ is a purely inseparable extension of $K_X$.  So $K_X$ has the
  property that if $\sigma\in\Aut(\ol{K})$ then $X_{\ol{K}}\cong\el{\sigma}{X_{\ol{K}}}$ if and only if $\sigma|_{K_X}=id$.

Theorem~\ref{koizumi} is due to Koizumi~\cite{Koizumi:1}.  We will
give a proof here based on the proof given in~\cite{Koizumi:1}, but
with more details.

We need to define a few notions and supply some needed results
first.

\section{Chow forms, Chow points, and Chow fields}
\begin{df}  Let $K$ be a field and let $X$ be a $K$-variety.  The
\emph{$K$-cycles} of  $X$ are the elements of the free $\Z$-module
over the irreducible closed $K$-subvarieties of $X$.  The
\emph{components} of a cycle $\nu:=\sum_{i=1}^n m_i[X_i]$ are the
$X_i$ with nonzero coefficients.  A $K$-cycle is \emph{effective} if
all of its components have positive coefficients.  A $K$-cycle is
called a \emph{$K$-divisor} of $X$ if all of its components
 have codimension $1$.
\end{df}

 \begin{df} Let $K\subseteq F$ be fields and let $X$ be a geometrically integral $K$-variety.
   Then we have a morphism of
schemes \[\varphi_{F/K}\colon X_F\to X.\]  If $Y$ is a
$K$-subvariety of $X$ we define $\varphi_{F/K}^*(Y)$ by
\[\varphi_{F/K}^*(Y)=\sum_y\mbox{length}(\mathcal{O}_{X_F,y})\ol y\]
where $y$ runs over the maximal points of $Y_F$ and $\ol y$ is the
closure of $y$. By linearity we can extend $\varphi_{F/K}^*$ to a
map from the $K$-cycles of $X$ to the $F$-cycles of $X_F$.  An
$F$-cycle $\nu_{(F)}$ on $X_F$ is said to be \emph{$K$-rational} if
there exists a $K$-cycle $\nu$ on $X$ such that
$\varphi_{F/K}^*(\nu)=\nu_{(F)}$.
\end{df}
\begin{rem}  Let $K\subseteq F$ be fields, let $X$ be a geometrically integral
$K$-variety,
  and let $Y$ be a closed
subvariety of $X_F$.  Then $Y$ is a $K$-rational cycle of $X_F$ if
and only if $Y\subseteq X_F$ is defined over $K$.
\end{rem}

\begin{thm}[Chow]Let $K$ be a field and let $X\subseteq\mathbb{P}_K^n$ be a $K$-variety of dimension
$r$.  The set of $(r+1)$-tuples of hyperplanes in $\mathbb{P}_K^n$
which have a common intersection with $X$ are parameterized by an
irreducible hypersurface $H\subset (\mathbb{P}^n_K)^{r+1}$.  The
hypersurface $H$ is given by a multi-homogeneous form $f$ in the set
of $n+1$ variables corresponding to each ${\mathbb P}^n_K$.  The
form $f$ is called the Chow form of $X$ and is unique up to
multiplication by an element of $K^{\times}$.
\end{thm}
\begin{proof} See \S1 of~\cite{Chow} for the original proof  or Chapter~8
of~\cite{Rydh} for a more modern treatment.
\end{proof}
\begin{df}Let $K$ be a field, let $X\subseteq\mathbb{P}^n_K$ be a $K$-variety, and let
$\nu:=\sum_{i=1}^l m_i[X_i]$ be an effective $K$-cycle on $X$.  The
\emph{Chow form} of $\nu$ is the product $\prod_{i=1}^l
f_{X_i}^{m_i}$ where for each $i$, $f_{X_i}$ is the Chow form of
$X_i$.
\end{df}
\begin{df}Let $K$ be a field, let $X\subseteq\mathbb{P}^n_K$ be a $K$-variety, let $\nu$
be an effective $K$-cycle on $X$, and let $f$ be the Chow form of
$\nu$.  The coefficients of $f$ can be seen as a point of projective
space called the \emph{Chow point} of $\nu$.
\end{df}
\begin{lem}\label{chowpt} Let $K$ be a field and let $\nu$ be an effective
$K$-cycle in $\mathbb{P}^n_K$.  Then the Chow point of $\nu$
uniquely determines $\nu$.
\end{lem}
\begin{proof} This follows from Proposition~8.24 on page~67 and
Proposition~8.15 on page~64 of~\cite{Rydh}.
\end{proof}
\begin{df}Let $F$ be a field, let $X\subseteq\mathbb{P}^n$ be an $F$-variety, and let $\nu$
be an effective $F$-cycle on $X$.  The \emph{Chow field} $K$ of
$\nu$ is the field obtained by adjoining to the prime field of $F$
the ratios of the nonzero coefficients of the Chow point of $\nu$.
\end{df}
\begin{lem}\label{samuel}Let $F$ be a field,
let $\nu$ be a $F$-cycle on $\mathbb{P}^n_F$, and let $K$ be the
Chow field of $\nu$. Then $\nu$ is rational over a finite purely
inseparable extension of $K$. Furthermore, $\nu$ is rational over
$K$ if $K$ is perfect or if there exists a geometrically integral
$K$-variety $X\subseteq\mathbb{P}^n_K$ such that $\nu$ is a divisor
of $X_F\subseteq\mathbb{P}^n_F$.
\end{lem}
\begin{proof}See page~47 of~\cite{Samuel}.
\end{proof}
\section{Curves as divisors}
\begin{df} Let $F$ be a field.  A \emph{hypersurface} $H_f\subset \mathbb{P}^n_F$ of degree $d$ over $F$ is
an
$F$-variety given by a homogeneous form of $f$ of degree $d$.
\end{df}
\begin{lem}\label{intersection}  Let $F$ be an algebraically closed field.  Suppose we have
 \[X\subset V\subseteq\mathbb{P}^n,\]
 where $X$ is a curve, and $V$ is a smooth projective variety of dimension $r\ge 3$,
  all over $F$.  Then for all sufficiently
large $d$, there exists a  hypersurface $H$ of degree $d$ in
$\mathbb{P}^n$ such that $H\cap V$ is a smooth projective variety of
dimension $r-1$ containing $X$.  Furthermore, the set of
hypersurfaces of degree $d$ with this property forms an open dense
subset of the projective space parameterizing all hypersurfaces of
degree $d$ which contain $X$.
\end{lem}
\begin{proof} Fix $d\in\Z_{>0}$.
Let $\mathscr{I}_X$ be the ideal sheaf of $X$ in ${\mathbb P}^n$ and
let
   $W:=\Gamma(\mathbb{P}^n,
\mathscr{I}_X(d))$.
    Let $PW$ be the projective space associated with the vector space $W$.  So we can view
     $PW$ as the set of hypersurfaces
$H\subset\mathbb{P}^n$ of degree $d$ which contain $X$.   Any
element in $PW$ is determined by a nonzero global section $f\in W$.
For a closed point $v\in V$, let
\[B_v:=\{H\in PW\mid \text{$v \in H \cap V$ and $H \cap V$ is not smooth of dimension $r-1$ at $v$.}\}
\]
Suppose that $H\in PW$ and suppose that $H\not\in B_v$ for all
closed points $v\in V$. Since the dimension of $V$ is larger that
$2$, by Corollary~7.9 on page 244 of~\cite{Hartshorne},
 $H\cap V$ must be connected, hence irreducible since $H \cap V$
is smooth.  So $H\cap V$ must be a smooth projective variety of
dimension $r-1$ containing $X$.

For a closed point $v\in V$, fix $f_0\in
\Gamma(\mathbb{P}^n,\mathcal{O}_{\mathbb{P}^n}(d))$ such that
$v\not\in H_{f_0}$. Let $\mathcal{O}_{v,V}$ be the local ring of $v$
on $V$ with maximal ideal $\mf{m}_v$.   We define a map of
$F$-vector spaces
\[\psi_v\colon W\to \mathcal{O}_{v,V}/\mf{m}_v^2\]
by letting $\psi_v(f)$ be the image of $f/f_0$ in
$\mathcal{O}_{v,V}/\mf{m}_v^2$.  Then $v\in H_f\cap V$ if and only
if the image of $f/f_0$ is zero in $\mathcal{O}_{v,V}/\mf{m}_v$. If
if $v$ is a nonregular point of $V\cap H_f$ then $\psi_v(f)=0$.
 Note that if $\dim(H_f\cap V)\ne
r-1$ then $V\subseteq H_f$ and we must  have $\psi_v(f)=0$ since the
image of $f/f_0$ is zero in $\mathcal{O}_{v,V}$. So if $f$ is not in
the kernel of $\psi_v$, then either $v$ is a regular point of
$H_f\cap V$ and $\dim(H_f\cap V)=r-1$ or $v\not\in H_f$.  So a
hypersurface $H$ is in $B_v$ for some closed point $v\in V$, if and
only if $H=H_f$ for some $f\in\ker(\psi_v)$.

 Let $I$
be the homogeneous ideal corresponding to $\mathscr{I}_X$ and let
$\{g_1,\ldots, g_s\}$ be a set of generators of degrees $d_{g_1},
\ldots, d_{g_s}$, respectively.  {}From now on suppose that $d$ is
strictly greater that $d_{g_i}$ for all $i$.  Since $v$ is a closed
point and $F$ is algebraically closed, $\mf{m}_v$ is generated by
linear forms in the coordinates.

Suppose that $v\not\in X$ and suppose that
$\alpha\in\mathcal{O}_{v,V}/\mf{m}_v^2$.  Then $\alpha$ is equal to
the image of $p/q$ in $\mathcal{O}_{v,V}/\mf{m}_v^2$, for some
$p,q\in\Gamma(\mathbb{P}^n,\mathcal{O}_{\mathbb{P}^n}(1))$
 with $v\not\in H_q$.  Since $v\not\in X$ there exists $g_i\in \{g_1,\ldots, g_s\}$ such
that $v\not\in H_{g_i}$. Let $m:=d-\mbox{ deg}(g_i)$ and pick
$h_v\in \Gamma(\mathbb{P}^n,\mathcal{O}_{\mathbb{P}^n}(m-1)))$ with
$v\not\in H_{h_v}$.  Let $f_p:=pg_ih_v\in W$ and let
$f_q:=qg_ih_v\in W$.   Since $\ker(\psi_v)$ is independent of our
choice of $f_0$, the image of $\psi_v$ is also independent of our
choice of $f_0$. Since $v\not\in H_{qg_ih_v}$, we may assume that
$f_0=qg_ih_v$. Then $\psi_v(f_p)=\alpha$.    So $\psi_v$ is
surjective if $v\not\in X$.

 Now suppose that $v\in X$.  Let $\mathcal{O}_{v,
X}$ be the local ring of $v$ on $X$ with maximal ideal $\mf{n}_v$.
      Since $X\subset V$, there
exists a natural surjective map of $F$-vector spaces
 \[\rho_v\colon \mathcal{O}_{v,V}/\mf{m}_v^2\to
 \mathcal{O}_{v,X}/\mf{n}_v^2.\]
For
 $1\le i\le s$, choose $h_i\in
\Gamma(\mathbb{P}^n,\mathcal{O}(d-d_{g_i}))$ so that $v\not\in
H_{h_i}$ and  let $g_i'=g_ih_i\in W$.  Then the kernel of $\rho_v$
is generated by the $\psi_v(g_i')$.  The composition $\rho_v\psi_v$
is equal to the map
\[\phi_v\colon W\to \mathcal{O}_{v,X}/\mf{n}_v^2\] defined by
letting $\phi_v(f)$ be equal to the image of $f/f_0$ in
$\mathcal{O}_{v,X}/\mf{n}_v^2$. This is, of course, the zero map
since $X\subset H_f$ for all $f\in W$.  It follows that the sequence
 \[W\to \mathcal{O}_{v,X}/\mf{m}_v^2\to
 \mathcal{O}_{v,V}/\mf{n}_v^2\to 0 \]
   is exact.  Since both $X$ and $V$ are nonsingular,
$\mathcal{O}_{v,X}/\mf{n}_v^2$ has dimension $2$ and
$\mathcal{O}_{v,V}/\mf{m}_v^2$
 has dimension $r+1$.  We see that if $v\in X$, then the
image of $\psi_v$ has dimension $r-1$.

  We have
\[\dim_F\ker(\psi_v)=\left\{
                            \begin{array}{ll}
                              \dim_F W-(r+1), & \hbox{if $v\not\in V\cap X$} \\
                              \dim_F W-(r-1), & \hbox{if $v\in V\cap X$.}
                            \end{array}
                          \right.
\]

 Since $B_v$ is the projective space associated with the vector space $\ker(\psi_v)$ we have
 \[\dim B_v=\left\{
                            \begin{array}{ll}
                              \dim_F W-(r+2), & \hbox{if $v\not\in V\cap X$} \\
                              \dim_F W-r, & \hbox{if $v\in V\cap X$.}
                            \end{array}
                          \right.
\]
 Let $B_{V- X}\subset V\times PW$ be the subsets of $V\times PW$
consisting of all pairs $\la v, H\ra$ such that $v\in V- X$ is a
closed point and $H\in B_v$.  Define $B_X$ similarly, using $X$ in
place of $V- X$. Then $B_{V- X}$ and $B_X$ are the sets of closed
points of some quasi-projective varieties $B'_{V- X}$ and $B'_X$
respectively. The fibres of the projections $p_{V- X}\colon B'_{V-
X}\to (V- X)$ and $p_X\colon B'_X\to X$ have dimension $\dim_F
W-(r+2)$ and $\dim_F W-r$ respectively. So
\[\dim B_{V- X}\le\dim (V- X)+\dim_F W-(r+2)=\dim_F W-2\] and
\[\dim B_X\le \dim X + \dim_F W-r = \dim_F W-(r-1).\]
Let $B:=(B_{V- X}\cup B_X)\subseteq V\times PW$.  So $B$ is the set
of closed points of a projective variety $B'$.
 Then the
local ring at a point of $B'$ has dimension less than or equal to
$\dim_F W-2$.  It follows that each irreducible component of $B'$
has dimension less than or equal to $\dim_F W-2$. So $\dim B'\le
\dim_F W-2$. Consider the projection $p\colon B'\to PW$.
 Since $\dim B'\le
\dim_F W-2$ we must have $\dim p(B')\le \dim_F W-2$. Since $\dim
PW=\dim_F W-1$, $p(B')$ is properly contained in $PW$.  If $H\in PW-
p(B')$ then  $H$  satisfies the requirements of the lemma.

Finally note that since $B'$ is projective, $p\colon B'\to PW$ is
proper, so since $B'$ is closed in $V\times PW$, $p(B')$ is closed
in $PW$.  Therefore, $PW- p(B')$ is an open dense subset of $PW$.
\end{proof}
\begin{cor}\label{intinfinite}  Let $K$ be an infinite subfield of an
algebraically closed field $F$.  Let $X$ be an $F$-curve, let $V$ be
a smooth, projective,  geometrically integral $K$-variety of
dimension $r\ge 3$, and suppose we have
\[X\subset V_F\subseteq\mathbb{P}_F^n\]
 Then for sufficiently large $d$,
there exists a hypersurface $H$ of degree $d$  in $\mathbb{P}^n_K$
such that $(H\cap V)_F$ is a smooth projective variety of dimension
$r-1$ containing $X$.  Furthermore, $H\cap V$ is a geometrically
integral $K$-variety.
\end{cor}
\begin{proof} Fix $d\in\Z_{>0}$ and
let $PW$ be as in Lemma~\ref{intersection}.  Then by
Lemma~\ref{intersection}, for sufficiently large $d$, $PW$ contains
an open dense subset $U$ consisting of hypersurfaces $H$ of degree
$d$ such that $H\cap V$ is smooth of dimension $r-1$.

  Since $K$ is infinite, the
$K$-rational points of $PW$ are Zariski dense in $PW$ so they cannot
all be contained in the complement of $U$. So $U\cap PW(K)$ is
nonempty.  So there exists $H\in PW(K)$ such that $(H\cap V)_F$ is
smooth projective variety of dimension $r-1$ containing $X$.

As shown in the proof of Lemma~\ref{intersection}, $(H\cap V)_F$ is
smooth and connected (and therefore integral).  Since $F$ is
algebraically closed, to prove the last statement we show that
$H\cap V$ is smooth and connected. By Proposition~17.7.4(v) on pages
72-73 of~\cite{Groth}, since $(H\cap V)_F$ is smooth, $H\cap V$ is
smooth. Suppose that $H\cap V$ is not connected. Then $H\cap V=S\cup
T$ for some nonempty closed and open subschemes $S$ and $T$. This
implies that $(H\cap V)_F=S_F\cup T_F$. Since $S_F$ and $T_F$ are
nonempty closed and open subschemes of $(H\cap V)_F$, we get a
contradiction.

\end{proof}
\begin{cor}\label{divisor} Let $F$ be an algebraically closed field
and let $X\subset\mathbb{P}^n_F$ be a curve embedded in projective
space.  Let $K$ be the Chow field of $X$.  Then $K$ is the minimum
field of definition of $X\subseteq\mathbb{P}^n_F$.
\end{cor}
\begin{proof}  If $K$ is not infinite, then it is perfect and so our statement
follows by Lemma~\ref{samuel}.  Assume that $K$ is infinite.  If
$n=1$ this is trivial. If $n=2$ then $X$ is a divisor of
$\mathbb{P}_F^2$ so this follows by Lemma~\ref{samuel}. Assume that
$n>2$.  Then we may apply Corollary~\ref{intinfinite} to
$\mathbb{P}^n_F$ to obtain a smooth, projective, geometrically
integral $K$-variety $V^{n-1}$ of dimension $n-1$ with $X\subset
V^{n-1}_F\subset\mathbb{P}^n_F$.
 We may repeat this procedure,
applying Corollary~\ref{intinfinite} to $V^i$ to obtain a smooth,
projective, geometrically integral $K$-variety $V^{i-1}$ of
dimension $i-1$ with $X\subset V^{i-1}_F\subset\mathbb{P}^n_F$,
until we obtain a smooth, projective, geometrically integral
$K$-variety $V^2$ of dimension $2$ with $X\subset
V^2_F\subset\mathbb{P}^n_F$. Then $X$ is a divisor of $V^2_F$.  So
by Lemma~\ref{samuel}, the minimum field of definition of $X$ as a
subvariety of $\mathbb{P}^n_F$ is $K$.
\end{proof}
\section{$\mathcal{A}$-structures}\label{structures}

\begin{df} Let $X$ be a curve over an algebraically closed field $F$.  Let $d\in\Z_{>0}$
be such that every divisor on $X$ of degree $d$ is very ample.  Then
the pair $\mf{A}:=(X,d)$ is called an
 \emph{$\mathcal{A}$-structure} on $X$.
\end{df}

\begin{df}Let $X$ and $X'$ be two curves over an algebraically closed field $F$.
 We say that two $\mathcal{A}$-structures $\mf{A}=(X,d)$ and
$\mf{A}'=(X',d')$, on $X$ and $X'$ respectively, are
\emph{isomorphic} if and only if $X\cong X'$ and $d=d'$.
\end{df}

\begin{notation} Let $X$ be a curve over a field $K$ and let $D$ be a very ample divisor on
$X$.  Let $B$ be a basis for $\Gamma(X,\mathscr{L}(D))$, and let
$f_{B,D}\colon X\to\mathbb{P}^n$ be an embedding given by $B$. Then,
unless otherwise stated, $f_{B,D}(X)$ is to be viewed as a
subvariety of $\mathbb{P}^n$ and not as an abstract curve.  We will
call any embedding given by a basis for $\Gamma(X,\mathscr{L}(D))$
an embedding \emph{determined} by $D$.
\end{notation}

\begin{prop}\label{kpoints}Let $X$ be a variety defined over a field $K$ and let $\ol{K}$ be
an algebraic closure of $K$. Suppose that $K$ is finitely generated
over its prime subfield.  Let $M_X$ be the intersection over all
fields $L$ with $K\subseteq L\subseteq \ol{K}$ such that
$X(L)\ne\varnothing$.  Then $M_X=K$.
\end{prop}
\begin{proof} See~\cite{Rizov}.
\end{proof}

\begin{lem}\label{astructure} Let $X$ be a curve over an algebraically
closed field $F$ and let $\mf{A}=(X,d)$ be an
$\mathcal{A}$-structure on $X$.   Let ${\mathcal{E}}_{\mf{A}}$ be
the set of all
 non-degenerate embeddings of $X$ determined by divisors of degree $d$ on
 $X$.  Then the set of Chow points of \[\{ f(X)\mid f\in{\mathcal{E}}_{\mf{A}}\}\]
 form a set of closed points $V_{\mf{A}}(F)$ of a geometrically reduced quasi-projective variety $V_{\mf{A}}$
 embedded in a projective space $\mathbb{P}^n$.
  Let $M_{\mf{A}}$ be the minimal field of definition of the closure $\ol{V_{\mf{A}}}\subseteq
  \mathbb{P}^n$.
  Then $V_{\mf{A}}\subseteq\mathbb{P}^n$ is defined over $M_{\mf{A}}$
  and $M_{\mf{A}}$ is
 the intersection of the fields of definition of the curves
 \[\{ f(X)\mid f\in{\mathcal{E}}_{\mf{A}}\}.\]
\end{lem}
\begin{proof} The first statement of the Lemma is proved in
Lemma~4 of~\cite{Matsusaka}.  See the proof of Theorem~2.2
in~\cite{Koizumi:1} for a simpler construction of $\ol{V_{\mf{A}}}$
 or page~104 of~\cite{Sh:2} for a similar
construction in the case of an abelian variety.  (Note that in the
terminology of~\cite{Koizumi:1, Matsusaka, Sh:2}, all varieties by
definition are geometrically reduced.)  In each construction it is
shown that $\ol{V_{\mf{A}}}\subseteq\mathbb{P}^n$ is defined over
all fields of definition for the curves
\[\{ f(X)\mid f\in{\mathcal{E}}_{\mf{A}}\}.\]  In the proof of
Lemma~4 of~\cite{Matsusaka} it is shown that
$V_{\mf{A}}\subseteq\mathbb{P}^n$ is defined  over all fields of
definition for the curves
\[\{ f(X)\mid f\in{\mathcal{E}}_{\mf{A}}\}.\]  In Theorem~3
of~\cite{Matsusaka}, it is shown that
$V_{\mf{A}}\subseteq\mathbb{P}^n$ is defined over the minimum field
of definition  $M_{\mf{A}}$ of
$\ol{V_{\mf{A}}}\subseteq\mathbb{P}^n$.

By Corollary~\ref{divisor}, if $f\in {\mathcal{E}}_{\mf{A}}$ then
the minimum field of definition of $f(X)$ is the Chow field of
$f(X)$. Let $\ol{M_{\mf{A}}}\subseteq F$ be the algebraic closure of
$M_{\mf{A}}$. Since $M_{\mf{A}}$ is finitely generated over its
prime field, by Proposition~\ref{kpoints}, $M_{\mf{A}}$ is the
intersection of all fields $L$ with $M_{\mf{A}}\subseteq
L\subseteq\ol{M_{\mf{A}}}$ such that $V_{\mf{A}}(L)\ne\varnothing$.
For each $L$ with $M_{\mf{A}}\subseteq L\subseteq F$, we have
$V_{\mf{A}}(L)\ne\varnothing$ if and only if $L$ is a field of
definition of $f(X)$ for some $f\in{\mathcal{E}}_{\mf{A}}$.  Since
$\ol{M_{\mf{A}}}\subseteq F$, it follows that
 the intersection of the fields of definition of the curves
 \[\{ f(X)\mid f\in{\mathcal{E}}_{\mf{A}}\}\] is equal to $M_{\mf{A}}$ .
\end{proof}
\begin{rem}\label{isomastruct} Let $X$ and $X'$ be two isomorphic curves over an algebraically
closed field $F$ and let $\mf{A}=(X,d)$ and $\mf{A}'=(X',d)$ be
$\mathcal{A}$-structures on $X$ and $X'$ respectively.
 It is easy to see, following the
notation of  Lemma~\ref{astructure}, that
$V_{\mf{A}}(F)=V_{\mf{A}'}(F)$ and $M_{\mf{A}}=M_{\mf{A}'}$.
\end{rem}
\begin{lem}\label{twoastruct}  Let $X$ be a curve over an algebraically closed field
$F$, Let $d_1, d_2\in Z_{>0}$, and suppose that $\mf{A}_1=(X,d_1)$
and $\mf{A}_2=(X,d_2)$ are two $\mathcal{A}$-structures on $X$.
Then, following the notation of
 Lemma~\ref{astructure}, $M_{\mf{A}_1}=M_{\mf{A}_2}$.
\end{lem}
\begin{proof} This is proven on page~52
 of~\cite{Koizumi:1}
\end{proof}
\begin{thm}[Koizumi]\label{koizumi} Let $X$ be a curve over a field
$K$, let $\ol{K}$ be an algebraic closure of $K$, let $d\in\Z_{>0}$,
and suppose that $\mf{A}=(X_{\ol{K}},d)$ is an
$\mathcal{A}$-structure on $X_{\ol{K}}$. Then, following the
notation of Lemma~\ref{astructure}, $M_{\mf{A}}$ is the field of
moduli of $X$. Furthermore, the subfield $\ol{K}^H$ of $\ol{K}$
fixed by
\[H:=\{\sigma\in\Aut(\ol{K})\mid X_{\ol{K}}\cong \el{\sigma}{X_{\ol{K}}}\}\]  is a
purely inseparable extension of $M_{\mf{A}}$.
\end{thm}
\begin{proof} Let $K_X$ be the field of moduli of $X$.  It is immediate from
Lemma~\ref{astructure} that $K_X\subseteq M_{\mf{A}}$.  Let $L$ be a
field of definition of $X_{\ol{K}}$ and let $X'$ be an $L$-curve
such that $X'_{\ol{K}}\cong X_{\ol{K}}$.  Choose an effective
divisor $D$ on $X'$ of degree $d'$ large enough so that
$(X_{\ol{K}}, d')$ is an $\mathcal{A}$-structure $\mf{A}'$ on
$X_{\ol{K}}$.  Let $B$ be a basis for $\Gamma(X', \mathscr{L}(D))$.
Then $f_{B,D}(X)$ is defined over $L$, so $M_{\mf{A}'}\subseteq L$.
By Lemma~\ref{twoastruct}, $M_{\mf{A}'}=M_{\mf{A}}$.
    This implies that $M_{\mf{A}}$ is contained in every field of definition of $X_{\ol{K}}$.
     Thus
$K_X=M_{\mf{A}}$.

Suppose that $\sigma\in\Aut(\ol{K})$ and let $V_{\mf{A}}$ be as in
Lemma~\ref{astructure}.  It is clear that
$\mf{A}^{\sigma}:=(\el{\sigma}{X_{\ol{K}}},d)$ is an
$\mathcal{A}$-structure on $\el{\sigma}{X_{\ol{K}}}$ and if
$X_{\ol{K}}\cong\el{\sigma}{X_{\ol{K}}}$, then
$V_{\mf{A}}=\el{\sigma}{V_{\mf{A}}}$ and so
$\sigma|_{M_{\mf{A}}}=Id$.

Conversely, suppose that $\sigma\in\Aut(\ol{K})$ and
$\sigma|_{M_{\mf{A}}}=Id$.  Then for all $P\in V_{\mf{A}}(\ol{K})$
we have $P^{\sigma}\in V_{\mf{A}}(\ol{K})$.  Let $f\colon
X_{\ol{K}}\to\mathbb{P}^n$ be a non-degenerate embedding of
$X_{\ol{K}}$ determined by a divisor of degree $d$ and let $P$ be
the Chow point of $f(X_{\ol{K}})$.  Then the Chow point of
$\el{\sigma}{f(X_{\ol{K}})}\subset\mathbb{P}^n$ is $P^{\sigma}$.
Since $P^{\sigma}\in V_{\mf{A}}(\ol{K})$, by Lemma~\ref{chowpt},
$\el{\sigma}{f(X_{\ol{K}})}$ is equal to $g(X_{\ol{K}})$ for some
embedding $g$ of $X$, determined by a divisor of $X_{\ol{K}}$ of
degree $d$. So $X_{\ol{K}}\cong\el{\sigma}{X_{\ol{K}}}$.  It follows
that $\sigma|_{M_{\mf{A}}}=Id$ if and only if
$X_{\ol{K}}\cong\el{\sigma}{X_{\ol{K}}}$.

Thus $\ol{K}^H$ is a purely inseparable extension of $K_X$.
\end{proof}

\begin{cor}\label{finiteext}Let $X$ be a curve over a field $K$ and let $K_X$ be
 its field of moduli.  Then $X$ is definable over a finite separable
 extension of $K_X$.
\end{cor}
\begin{proof}  Let $\ol{K}$ be an algebraic closure of $K$ and choose $d\in\Z_{>0}$ so that $\mf{A}:=(X_{\ol{K}},d)$
is an $\mathcal{A}$-structure on $X_{\ol{K}}$.  Let $V_{\mf{A}}$ be
the quasi-projective variety corresponding to $\mf{A}$.  Let
$K_X^s\subseteq\ol{K}$ be the separable closure of $K_X$.  By
Lemma~\ref{seppoint} below, $V_{\mf{A}}(K_X^s)\ne\varnothing$.  So
there exists a finite separable extension field $L$ of $K_X$ such
that $V_{\mf{A}}(L)\ne\varnothing$.  Then $X$ is definable over $L$.
\end{proof}
Since the following lemma is well known, we omit the proof.
\begin{lem}\label{seppoint}  Let $K$ be a separably closed field and
let $X$ be a geometrically reduced $K$-variety.  Then $X(K)$ is
Zariski dense in $X$.
\end{lem}

 The proof of the
following lemma, although slightly different from the proof of
Proposition~3.2 in~\cite{Sekiguchi:3}, is based on an idea given in
the proof of Proposition~3.2 of~\cite{Sekiguchi:3}. Let $K\subseteq
F$ be fields where $F$ is a purely inseparable extension of $K$. Let
$V$ be a polarized abelian variety over $F$. Proposition~3.2
of~\cite{Sekiguchi:3} states that $V$ is definable over $K$ if $K$
contains the field of moduli of $V$.

\begin{lem}\label{insep}Let $X$ be a curve over a field $F$ and suppose that
$F$ is a purely inseparable extension of a field $K$.  Suppose that
$X$ is definable over a finite separable extension of $K$.  Then $X$
is definable over $K$.
\end{lem}
\begin{proof} Let $L$ be a finite separable extension of $K$
contained in an algebraic closure $\ol{F}$ of $F$.  Suppose there
exists a curve $X'$ over $L$ such that $X'_{\ol{F}}\cong
X_{\ol{F}}$.  Then there exists a finite extension field $F'$ of
$F$, containing $L$ such that $X'_{F'}\cong X_{F'}$ over $F'$.
Replacing $F'$ with a larger field if necessary, we may assume that
$F'=F_1F_2$ where $F_1$ and $F_2$ are subfields of $F'$ with
$F_1\cap F_2=K$, where $F'/F_1$ is Galois, and where $F'/F_2$ is
purely inseparable.  So $F_1/K$ is purely inseparable and $F_2/K$ is
Galois.  It follows that $F\subseteq F_1$ and $L\subseteq F_2$.
Furthermore, we have $\Gal(F'/F_1)\cong\Gal(F_2/K)$. We have the
following picture:
\begin{diagram}
         &        &  F'      &         & \\
         & \ruLine^{\text{Galois}}&        & \luLine^{\text{p. ins.}} &\\
       F_1\supseteq F &        &        &          & F_2\supseteq L\\
         & \rdLine_{\text{p. ins.}}&        & \ldLine_{\text{Galois}} &\\
         &        & K   &          &\\
\end{diagram}

 Let
$X_1:=X_{F_1}$, let $X_2:=X'_{F_2}$, and let $F_1(X_1)$ and
$F_2(X_2)$ be their function fields respectively.  Let $F'(X_{F'})$
be the function field of $X_{F'}$. Then
\[F'(X_{F'})=F_2(X_2)\otimes_{F_2} F'=F_1(X_1)\otimes_{F_1} F',\]
  $F'(X_{F'})/F_1(X_1)$ is Galois with
\[\Gal(F'(X_{F'})/F_1(X_1))\cong\Gal(F'/F_1),\] and $F'(X_{F'})/F_2(X_2)$ is
purely inseparable.   Let $M:=F_1(X_1)\cap F_2(X_2)$. Then
$F_2(X_2)/M$ is Galois, with Galois group isomorphic to
$\Gal(F'(X_{F'})/F_1(X_1))$, and $F_1(X_1)/M$ is purely inseparable.
We have the following picture:
\begin{diagram}
         &        & F'(X_{F'})      &         & \\
         & \ruLine^{\text{Galois}}&        & \luLine^{\text{p. ins.}} &\\
       F_1(X_1) &        &        &          & F_2(X_2)\\
         & \rdLine_{\text{p. ins.}}&        & \ldLine_{\text{Galois}} &\\
         &        & M   &          &\\
         &        &             \dLine            &          &\\
         &        &                   K      &          &\\
\end{diagram}
and we have
\[\Gal(F_2(X_2)/M)\cong\Gal(F'(X_{F'})/F_1(X_1))\cong\Gal(F'/F_1)\cong\Gal(F_2/K).\]
An isomorphism $\Gal(F'(X_{F'})/F_1(X_1))\to \Gal(F_2/K)$ is given
by \[\sigma\mapsto \sigma|_{F_2}\] and an isomorphism
$\Gal(F'(X_{F'})/F_1(X_1))\to \Gal(F_2(X_2)/M)$ is given by
\[\sigma\mapsto \sigma|_{F_2(X_2)}.\]
 Since
\[\sigma|_{F_2}=(\sigma|_{F_2(X_2)})|_{F_2},\] we have
\[\Gal(F_2(X_2)/M)|_{F_2}=\Gal(F'(X_{F'})/F_1(X_1))|_{F_2}=\Gal(F_2/K).\]
Since $\Gal(F_2(X_2)/M)\subseteq \Aut(F_2(X_2)/K)$, the sequence
\[1\to\Aut(F_2(X_2)/F_2)\to\Aut(F_2(X_2)/K)\to\Gal(F_2/K)\to 1\]
is split exact.  So by Theorem~\ref{weil}, $X$ is definable over
$K$.
\end{proof}
\begin{cor} Let $X$ be a curve over a field $F$.  Suppose that $F$
is a purely inseparable extension of a field $K$ that contains the
field of moduli of $X$.  Then $X$ is definable over $K$.
\end{cor}
\begin{proof}  Let $K_X$ be the field of moduli of $X$.  By Corollary~\ref{finiteext}, $X$ is definable over a
finite separable extension $L$ of $K_X$.  So $X$ is definable over
the compositum $LK$ which is a finite separable extension of $K$.
Then by Lemma~\ref{insep}, $X$ is definable over $K$.
\end{proof}

\section{Fields of moduli relative to Galois extensions}
We introduce another definition of ``field of moduli" that is
commonly used and is defined relative to a given Galois extension.

\begin{df} Let $X$ be a curve over a field $F$ and let $K$ be a subfield of $F$
such that $F/K$ is Galois. The \emph{field of moduli of $X$
relative} to the extension $F/K$ is defined as the fixed field
$F^{H}$ of
\[H:=\{\sigma\in \Gal(F/K)\mid X\cong \el{\sigma}{X}\mbox{ over } F\}.\]
\end{df}
\begin{prop}\label{closed} Let $X$ be a curve over a field $F$, let $K$ be a subfield of $F$
such that $F/K$ is Galois, let
\[H:=\{\sigma\in \Gal(F/K)\mid X\cong \el{\sigma}{X}\mbox{ over } F\},\] and let $K_m$ be the field of moduli
of $X$ relative to $F/K$. Then the subgroup $H$ is a closed subgroup
of $\Gal(F/K)$ for the Krull topology.  That is,
\[H=\Gal(F/K_m).\]
 The field of $K_m$ is contained in
each field of definition between $K$ and $F$ (in particular, $K_m$ is a finite extension of $K$).  Hence if the
field of moduli is a field of definition, it is the smallest field of definition between $F$ and $K$.  Finally,
the field of moduli of $X$ relative to the extension $F/K_m$ is $K_m$.
\end{prop}
\begin{proof}
See Proposition 2.1 in~\cite{Debes:1}.
\end{proof}

Let $X$ be a curve over a field $F$ and let $K$ be a subfield of $F$
such that $F/K$ is Galois.   The following theorem gives necessary
and sufficient conditions for $L$ to be a field of definition for
$X$.

\begin{thm}[Weil]\label{weil}Let $X$ be a curve over a field $F$ and let
$K$ be a subfield of $F$ such that $F/K$ is Galois.  Let
$\Gamma=\Gal(F/K)$ and suppose for all $\sigma\in\Gamma$ there
exists an $F$-isomorphism $f_{\sigma}\colon X\to\el{\sigma}{X}$ such
that
\[f^{\sigma}_{\tau}f_{\sigma}=f_{\sigma\tau},\;\;\mbox{for all $\sigma,\tau\in\Gamma$.}\]
Then there exist a $K$-curve $X'$ and an isomorphism
\[f\colon X\to X'_F\] defined over $F$ such that
\[f_{\sigma}=(f^{-1})^{\sigma}f,\mbox{ for all }\sigma\in\Gamma.\]
\end{thm}
\begin{proof} See the proof of Theorem 1 of~\cite{Weil:1}.
\end{proof}
\begin{rem} Let $X$ be a curve over a field $F$ and let
$K$ be a subfield of $F$ such that $F/K$ is Galois, and  $F(X)$ be
the function field of $X$.   If $K$ is the field of moduli of $X$
relative to $F/K$ we get an exact sequence
\[1\to\Aut(F(X)/F)\to\Aut(F(X)/K)\to\Gal(F/K)\to 1.\]
The conditions of Theorem~\ref{weil} are satisfied precisely when
this sequence is split exact.  That is, when there exists a subgroup
$H\subseteq\Aut(F(X)/K)$ isomorphic to $\Gal(F/K)$ such that
$H|_F=\Gal(F/K)$.
\end{rem}

The following four results of D\`ebes, Emsalem , and Douai will be
of use to us.
 They rely on the notions of a cover and the field of moduli of a cover,
for which we refer the reader to \S2.4 in~\cite{Debes:2}.
\begin{thm}\label{bk} Let $F/K$ be a Galois extension and
 $X$ be a curve of genus larger than $1$ defined over $F$ with $K$ as field
of moduli.   Then there exists a $K$-model $B$ of the curve
$X/\Aut(X)$ such that the cover $X\to B_F$ with $K$-base $B$ is of
field of moduli $K$.
\end{thm}
\begin{proof}See Theorem 3.1 in~\cite{Debes:1}.  The authors make the additional assumption that
the characteristic of $K$ does not divide $|\Aut(X)|$ but do not use
it in their proof.
\end{proof}
\begin{cor}\label{finitefield}  Suppose that $K$ is a finite field
and that $F$ is algebraically closed.  Then $X$ can be defined over
$K$.
\end{cor}
\begin{proof} It suffices to show that the cover $X\to B_F$ with $K$-base $B$
can be defined over $K$, since a field of definition of the cover is
automatically a field of definition of $X$. By Theorem~\ref{bk}, the
field of moduli of the cover $X\to B_F$ with $K$-base $B$ is $K$. If
$K$ is a finite field then $\Gal(F/K)$ is a projective profinite
group.   In this case, by Corollary 3.3 of~\cite{Debes:2} the cover
$X\to B_F$ can be defined over $K$.
\end{proof}
\begin{cor}\label{qpt}Suppose that $F$ is algebraically closed and that $X$ is
a hyperelliptic curve.  If $B$ has a $K$-rational point, then $K$ is
a field of definition of $X$.
\end{cor}
\begin{proof}It suffices to show that the cover $X\to B_F$ with $K$-base $B$
can be defined over $K$, since a field of definition of the cover is
automatically a field of definition of $X$. By Theorem~\ref{bk}, the
field of moduli of the cover $X\to B_F$ with $K$-base $B$ is $K$. By
Corollary~\ref{finitefield}, we may assume that $K$ is infinite.
 Since $B\cong_K \mathbb {P}^1_K$, $B$ has a rational point off the
branch point set of $X\to B_F$. Then by Corollary 3.4 and
\textsection~2.9 of~\cite{Debes:2}, the cover can be defined over
$K$.
\end{proof}
\begin{cor}\label{pqpt} Suppose that $F$ is algebraically closed and that $|\Aut(X)|$ is prime to
the characteristic of $F$.  If $B$ has a $K$-rational point, then
$K$ is a field of definition of $X$.
\end{cor}
\begin{proof} This is Corollary~4.3(c) of~\cite{Debes:1}.
\end{proof}
The curve $B$ of Theorem~\ref{bk} and Corollary~\ref{qpt} is called
the canonical model of $X/\Aut(X)$ over the field of moduli of $X$.

Let $X$ be a curve over a field $K$ and let $K_X$ be the field of
moduli of $X$.  We now show the relationship between $K_X$ and the
fields of moduli of $X$ relative to Galois extensions of $K$.

\begin{thm}\label{relmod} Let $X$ be a curve over a field $K$ and let $K_X$ be the
field of moduli of $X$.  Then $X$ is definable over $K_X$ if and
only if given any algebraically closed field $F\supseteq K$, and any
subfield $L\subseteq F$ with $F/L$ Galois, $X_{F}$ can be defined
over its field of moduli relative to the extension $F/L$.
\end{thm}
\begin{proof}  Suppose that given any algebraically closed field $F\supseteq K$, and any
subfield $L\subseteq F$ with $F/L$ Galois, $X_{F}$ can be defined
over its field of moduli relative to the extension $F/L$.  By
Corollary~\ref{finiteext}, $X$ is definable over a finite separable
extension field $M$ of $K_X$.  Let $\ol{M}$ be an algebraic closure
of $M$.  Without loss of generality we may assume that $X$ is a
curve over $M$.  Then there exists a field $K_m$ with $K_X\subset
K_m\subset\ol{M}$ such that $\ol {M}/K_m$ is separable and $K_m/K_X$
is purely inseparable.
  Then the field $K_m$ is the
field of moduli of $X_{\ol{M}}$ relative to the extension
$\ol{M}/K_m$. By assumption, $X_{\ol{M}}$ is definable over $K_m$.
So by Lemma~\ref{insep}, $X_{\ol{M}}$ is definable over $K_X$.  Thus
$X$ is definable over $K_X$.

The other direction follows immediately from the definition of
$K_X$.
\end{proof}

\section{Moduli spaces}

 Let $X$ be a curve of genus $g$ ($\ge 2)$ over a field $K$, let
 $K_X$ be its field of moduli, and
let $\mathcal{M}_g$ be the coarse
 moduli space of curves of genus $g$ viewed as a scheme over the prime field of
 $K$.
The curve $X$ gives a morphism $\Spec K\to \mathcal{M}_g$ whose
image $x$ is a closed point of $\mathcal{M}_g$.   Let
$\textbf{K}(x)$ be the residue field at
 $x$.  In this section we show the relationship between $\textbf{K}(x)$ and
$K_X$.
\begin{thm}[Baily, 1962] Following the above notation, suppose the
characteristic of $K$ is $0$.  Then $\textbf{K}(x)=K_X$.

\end{thm}\noindent(cf.Baily~\cite{Baily:1})

\begin{thm}[Sekiguchi, 1985]Following the above notation, $K_X$ is a purely inseparable extension of
$\textbf{K}(x)$.
\end{thm}\noindent (cf.Sekiguchi~\cite{Sekiguchi:1})

\begin{thm}[Sekiguchi, 1985, 1988]There exist
both hyperelliptic and non-hyperelliptic curves whose fields of
moduli are nontrivial purely inseparable extensions of the residue
fields of the corresponding points on their moduli spaces.
\end{thm}\noindent (cf.Sekiguchi~\cite{Sekiguchi:1, Sekiguchi:2})

%% file: ch2.tex
\section{Notation and terminology}
Let $F$ be a field.  Throughout this thesis $\GL_n(F)$ denotes the
group of nonsingular $n \times n$ matrices over the field $F$ and
$\SL_n(F)$ denotes the subgroup of $\GL_n(F)$ consisting of elements
of determinant $1$.

Let $p$ be a prime number and let $\F_q$ be a finite field with
$q:=p^r$ elements, where $r>0$. Throughout this thesis, let
$\GU_n(\F_q)$ be the subgroup of $\GL_n(\F_{q^2})$ consisting of
matrices $M$ such that $M^{-1}$ is the transpose of the matrix
obtained from $M$ via the automorphism $c\mapsto c^q$ of $\F_{q^2}$
and let $\SU_n(\F_q):=\GU_n(\F_q)\cap\SL_n(\F_{q^2})$.

For any group $\G\subseteq\GL_n(F)$, $\mbox{P}\G$ denotes $\G/Z(\G)$
where $Z(\G)$ is the center of $\G$.

 We will use a matrix with round brackets to denote an element of
$\GL_n(F)$ and a matrix with square brackets to denote the image in
$\PGL_n(F)$ of an element of $\GL_n(F)$.  For example,
\[\begin{pmatrix}
a_{11}&\cdots&a_{1n}\\
\cdot&\cdots&\cdot\\
a_{n1}&\cdots&a_{nn}\\
\end{pmatrix}\] denotes a matrix in $\GL_n(F)$ and
\[\begin{bmatrix}
a_{11}&\cdots&a_{1n}\\
\cdot&\cdots&\cdot\\
a_{n1}&\cdots&a_{nn}\\
\end{bmatrix}\]  denotes its image in $\PGL_n(F)$.

\section{Finite subgroups of the 2-dimensional projective general linear groups}
\indent\par Throughout this section let $F$ be an algebraically
closed field of characteristic $p$ with $p=0$ or
 $p>2$.
 \begin{lem}\label{subgroups}
Any finite subgroup $\G$ of $\PGL_2(F)$ is conjugate to one of the
following groups:

\begin{enumerate}
\item[Case] I: when $p=0$ or $|\G|$ is relatively prime to $p$.
\begin{enumerate}

\item  $\G_{C_n}:=\left\{\begin{bmatrix}

  \zeta^r & 0 \\
0 & 1 \\
\end{bmatrix}\colon r=0,1,\ldots,n-1\right\}\cong C_n,~n\ge 1$

\item $\G_{D_{2n}}:=\left\{\begin{bmatrix}

  \zeta^r & 0 \\
0 & 1 \\
\end{bmatrix},
\begin{bmatrix}

  0 & \zeta^r \\
1 & 0 \\
\end{bmatrix}\colon r=0,1,\ldots,n-1\right\}\cong D_{2n},~{n>1}$

\item  $\G_{A_4}:=\left\{\begin{bmatrix}
\pm 1 & 0 \\
0 & 1 \\
\end{bmatrix},
\begin{bmatrix}
0 & \pm 1 \\
1 & 0 \\
\end{bmatrix},
\begin{bmatrix}
  i^{\nu} & i^{\nu} \\
1 & -1 \\
\end{bmatrix},
\begin{bmatrix}
i^{\nu} & -i^{\nu} \\
1 & 1 \\
\end{bmatrix},
\begin{bmatrix}
1 & i^{\nu} \\
1 & -i^{\nu} \\
\end{bmatrix},\right.$

$\left. \begin{bmatrix}
-1 & -i^{\nu} \\
1 & -i^{\nu} \\
\end{bmatrix}\colon\nu=1,3\right\}\cong A_4$

\item  $\G_{S_4}:=\left\{\begin{bmatrix}
  i^{\nu} & 0 \\
0 & 1 \\
\end{bmatrix},
\begin{bmatrix}

  0 & i^{\nu} \\
1 & 0 \\
\end{bmatrix},
\begin{bmatrix}
i^{\nu} & -i^{\nu+\nu'} \\
1 &  i^{\nu'} \\
\end{bmatrix}\colon \nu, \nu'=0,1,2,3\right\}{\;\cong S_4}$

\item  $\G_{A_5}:=\left\{\begin{bmatrix}
  \epsilon^r & 0 \\
0 & 1 \\
\end{bmatrix},
\begin{bmatrix}

  0 & \epsilon^r \\
-1 & 0 \\
\end{bmatrix},
\begin{bmatrix}
  \epsilon^r\omega & \epsilon^{r-s} \\
1 & -\epsilon^{-s}\omega \\
\end{bmatrix},
\begin{bmatrix}
  \epsilon^r\overline{\omega} & \epsilon^{r-s} \\
1 & -\epsilon^{-s}\overline{\omega} \\
\end{bmatrix}\colon\right.$

$\left. r, s = 0,1,2,3,4\right\}\cong A_5$
\end{enumerate}
\noindent where $\omega:=\frac{-1+\sqrt{5}}{2}$,
$\overline{\omega}:=\frac{-1-\sqrt{5}}{2}$, $\zeta$ is a primitive
$n^{th}$ root of unity, $\epsilon$ is a primitive $5^{th}$ root of
unity, and $i$ is a primitive $4^{th}$ root of unity.
 \item [Case] II: when $|\G|$ is divisible by $p$.
\begin{enumerate}
  \item[(f)]
$\G_{\beta,A}:=\left\{ \begin{bmatrix}
                      \beta^k &a\\
                      0&1\end{bmatrix}\colon a\in A, \; k\in\Z\right\}$, where $A$ is a finite
                      additive subgroup of $F$ containing $1$ and $\beta$ is a root of unity such that $\beta A=A$
\item[(g)] $\PSL_2(\mathbb{F}_{q})$
  \item[(h)] $\PGL_2(\F_{q})$
\end{enumerate}
where $\F_{q}$ is the finite field with $q:=p^r$ elements, where
$r>0$.
\end{enumerate}
\end{lem}
\begin{proof}
See \S\S71-74 in~\cite{Weber:2} and Chapter 3 in~\cite{Suzuki:1}.
\end{proof}
\begin{rem} It
can be directly verified that $\G_{A_4}$ and $\G_{S_4}$ are
subgroups of $\PGL_2(F)$ when the characteristic of $F$ is $3$.
Indeed, in this case $\G_{A_4}$ is $\PGL_2(F)$ conjugate to
$\PSL_2(\F_3)$ and $\G_{S_4}$ is $\PGL_2(F)$ conjugate to
$\PGL_2(\F_3)$.  So the result of Lemma~\ref{normalizersub}(b) is
still valid in characteristic $3$.
\end{rem}

\begin{lem}\label{normalizersub} Let $N(\G)$ be the normalizer of $\G$ in
$\PGL_2(F)$.  Then
\begin{enumerate}
\item [(a)] $N(\G_{C_n})=\left\{\begin{bmatrix}
\alpha & 0 \\
0 & 1 \\
\end{bmatrix},
\begin{bmatrix}
0 & \alpha \\
1 & 0 \\
\end{bmatrix}: \alpha\in F^{\times}\right\}$ if $n>1$,
\item[(b)] $N(\G_{D_4})=\G_{S_4}$, $N(\G_{D_{2n}})=\G_{D_{4n}}$ if $n>2$,
 \item[(c)]$N(\G_{A_4})=\G_{S_4}$,
\item[(d)] $N(\G_{S_4})=\G_{S_4}$, \item[(e)] $N(\G_{A_5})=\G_{A_5}$,
\item[(g)]$N(\PSL_2(\mathbb{F}_{q}))=\PGL_2(\mathbb{F}_{q}),\mbox{ and }$ \item[(h)]
$N(\PGL_2(\mathbb{F}_{q}))=\PGL_2(\mathbb{F}_{q})$.
\end{enumerate}
\end{lem}
\begin{proof} \hfill
\begin{enumerate}
\item[(a)] See \S71 in~\cite{Weber:2}. \item[(b)]
 Since $\G_{D_4}$ is a normal subgroup of $\G_{S_4}$, $\G_{S_4}\subseteq N(\G_{D_4})$. Conjugation of $\G_{D_4}$ by
$\G_{S_4}$ gives a homomorphism $\G_{S_4}\to\Aut(D_4)\cong S_3$.  A
computation shows that the centralizer $Z$ of $\G_{D_4}$ in
$\PGL_2(F)$ is $\G_{D_4}$.  The kernel of this homomorphism is
$Z\cap \G_{S_4}=Z$.  Since $\G_{S_4}/Z\cong S_3$, every automorphism
of $\G_{D_4}$ is given by conjugation by an element of $\G_{S_4}$.
Let $U\in N(\G_{D_4})$. Then $UV\in Z=\G_{D_4}$ for some $V\in
\G_{S_4}$, so $U\in \G_{S_4}$.

For $n>2$, see \S71 in~\cite{Weber:2}.
\item[(c)] Since $\G_{D_4}$ is a
characteristic subgroup of $\G_{A_4}$, $N(\G_{A_4})\subseteq
N(\G_{D_4}){=\G_{S_4}}$. As $\G_{A_4}$ is normal in $\G_{S_4}$, we
get $N(\G_{A_4})=\G_{S_4}$.
\item[(d)] Since $\G_{A_4}$ is a characteristic subgroup of
$\G_{S_4}$, $N(\G_{S_4})\subseteq N(\G_{A_4}){=\G_{S_4}}$. Thus
$N(\G_{S_4})=\G_{S_4}$. \item[(e)] Conjugation of $\G_{A_5}$ by
$N(\G_{A_5})$ gives a homomorphism $N(\G_{A_5})\to\Aut(A_5)$. The
kernel of this homomorphism is the centralizer of $\G_{A_5}$ in
$N(\G_{A_5})$, which is just the centralizer $Z$ of $\G_{A_5}$ in
$\PGL_2(F)$. A computation shows that $Z$ is just the identity.
Since $\Aut(A_5)$ is finite, $N(\G_{A_5})$ is a finite subgroup of
$\PGL_2(F)$.  Since $\G_{A_5}\subseteq N(\G_{A_5})$, by
Lemma~\ref{subgroups} we must have $N(\G_{A_5})=\G_{A_5}$.

\item[(g)] We first show that $N(\PSL_2(\mathbb{F}_{q}))$ is finite.  Conjugation of
$\PSL_2(\mathbb{F}_{q})$ by $N(\PSL_2(\mathbb{F}_{q}))$ gives a
homomorphism
$N(\PSL_2(\mathbb{F}_{q}))\to\Aut(\PSL_2(\mathbb{F}_{q}))$. The
kernel of this homomorphism is the centralizer $Z$ of
$\PSL_2(\mathbb{F}_{q})$ in $\PGL_2(F)$.  A computation shows that
$Z$ is just the identity.  Since $\Aut(\PSL_2(\mathbb{F}_{q}))$ is
finite, so is $N(\PSL_2(\mathbb{F}_{q}))$.  By Lemma~\ref{subgroups}
any finite subgroup of $\PGL_2(F)$ containing
$\PSL_2(\mathbb{F}_{q})$ must be isomorphic to either
$\PGL_2(\mathbb{F}_{q'})$ or $\PSL_2(\mathbb{F}_{q'})$ for some
$q'$. Since $\SL_2(\mathbb{F}_{q})$ is normal in
$\GL_2(\mathbb{F}_{q})$, $\PSL_2(\mathbb{F}_{q})$ is a normal
subgroup of $\PGL_2(\mathbb{F}_{q})$.  So
$\PGL_2(\mathbb{F}_{q})\subseteq N(\PSL_2(\mathbb{F}_{q}))$, in
particular $\PSL_2(\mathbb{F}_{q})$ is strictly contained in
$N(\PSL_2(\mathbb{F}_{q}))$.  By the corollary on page 80
of~\cite{Suzuki:1}, $\PSL_2(\mathbb{F}_{q'})$ is simple for $q'>3$.
It follows that $N(\PSL_2(\mathbb{F}_{q}))\ne
\PSL_2(\mathbb{F}_{q})$ for $q\ge 3$.  By Theorem 9.9 on page 78
of~\cite{Suzuki:1}, the only nontrivial normal subgroup of
$\PGL_2(\mathbb{F}_{q'})$ is $\PSL_2(\mathbb{F}_{q'})$ if $q'>3$.
Therefore $N(\PSL_2(\mathbb{F}_{q}))=\PGL_2(\mathbb{F}_{q})$.
\item[(h)] Clear from the proof of the previous case.
\end{enumerate}
\end{proof}

\begin{lem}\label{propersubgroups}  Let $\G$ be one of the subgroups listed in Lemma~\ref{subgroups} and let
$\G'\subset\PGL_2(F)$ be a finite subgroup which properly contains
$\G$.
\begin{enumerate}
    \item[(a)] If $\G=\G_{D_4}$, then $\G'=\G_{A_4}$ , $\G'\cong\G_{S_4}$, $\G'\cong
    \G_{A_5}$, $\G'\cong D_{4n}$ for some $n>1$, $\G'\cong\PSL_2(\F_q)$ with $q>3$,
    or $\G'\cong\PGL_2(\F_q)$ for some finite field $\F_q$.
    \item[(b)] If $\G=\G_{D_6}$, then $\G'\cong S_4$, $\G'\cong A_5$, $\G'=\G_{D_{6n}}$ for some $n>1$,
    $\G'\cong\PSL_2(\F_q)$,
    or $\G'\cong\PGL_2(\F_q)$ for some finite field $\F_q$ where $3|q-1$.
    \item[(c)] If $\G=\G_{D_8}$, then $\G'\in\{\G_{S_4}, \G_{D_{4n}}\colon n>1\}$,
    $\G'\cong\PSL_2(\F_q)$,
    or $\G'\cong\PGL_2(\F_q)$ for some finite field $\F_q$ where $4|q-1$.
    \item[(d)] If $\G=\G_{D_{10}}$, then $\G'\cong A_5$, $\G'=\G_{D_{10n}}$ for some
    $n>1$,
    $\G'\cong\PSL_2(\F_q)$,
    or $\G'\cong\PGL_2(\F_q)$ for some finite field $\F_q$ where $5|q-1$.
    \item[(e)] If $\G=\G_{D_{2n}}$ with $n>5$,
then $\G'=\G_{D_{2n'}}$ for some $n'>n$ with $n|n'$,
$\G'\cong\PSL_2(\F_q)$,
    or $\G'\cong\PGL_2(\F_q)$ for some finite field $\F_q$ where $n|q-1$.
    \item[(f)] If $\G=\G_{A_4}$, then $\G'=\G_{S_4}$, $\G'\cong
    A_5$,
    $\G'\cong\PSL_2(\F_q)$,
    or $\G'\cong\PGL_2(\F_q)$ for some finite field $\F_q$.
   \end{enumerate}
\end{lem}
\begin{proof}\hfill \begin{enumerate}
    \item[(a)]   By Lemma~\ref{subgroups}, $\G'$ must be isomorphic to $A_4$,
    $S_4$, $A_5$,
     $D_{4n}$ for some $n>1$, $\PSL_2(\F_q)$, or $\PGL_2(\F_q)$ for some finite field $\F_q$.
      Suppose that $\G'$ is isomorphic to $A_4$.  (Note that $\PSL_2(\F_3)\cong A_4$.)  The group $A_4$ has a unique subgroup isomorphic
      to $D_4$.  So $\G'\subseteq N(\G_{D_4})$.
    Lemma~\ref{normalizersub}, $N(\G_{D_4})=\G_{S_4}$.  Since
    $\G_{A_4}$ is the unique subgroup of
    $\G_{S_4}$ isomorphic to $A_4$, the result follows.

    \item[(b)]  Since $D_6$ is not a subgroup of $A_4$, by Lemma~\ref{subgroups},
     $\G'\cong S_4$, $\G'\cong A_5$,  $\G'\cong\G_{D_{6n}}$ for some $n>1$,
    $\G'\cong\PSL_2(\F_q)$,
    or $\G'\cong\PGL_2(\F_q)$ for some finite field $\F_q$.  In the last two cases note
    that $\PSL_2(\F_q)$ and $\PGL_2(\F_q)$ have a diagonal element of order $3$ if and only if $3|q-1$.

    If $\G'\cong D_{6n}$ for some $n>1$, then
     $\G'$ has an element of order $3n$ that commutes with the elements of order $3$ in $\G_{D_6}$.  A computation
     shows that $\G'=\G_{D_{6n}}$.
    \item[(c)]  Since $D_8$ is not a subgroup of $A_4$, $A_5$, or a cyclic group, by Lemma~\ref{subgroups},
     $\G'\cong S_4$, $\G'\cong D_{8n}$ for some $n>1$, $\G'\cong\PSL_2(\F_q)$,
    or $\G'\cong\PGL_2(\F_q)$ for some finite field $\F_q$.  In the last two cases note
    that $\PSL_2(\F_q)$ and $\PGL_2(\F_q)$ have a diagonal element of order $4$ only if $4|q-1$.

    Suppose that $\G'\cong S_4$.  It can be shown that $S_4$ has $3$
    subgroups isomorphic to $D_8$ and that each one contains a
     subgroup isomorphic to $D_4$ which is normal in $S_4$.  So
     $\G'\subseteq N(\G_{D_4})$.  By Lemma~\ref{normalizersub},
     $N(\G_{D_4})=\G_{S_4}$.  So $\G'=\G_{S_4}$.

     If $\G'\cong D_{8n}$ for some $n>1$, then
     $\G'$ has an element of order $4n$ that commutes with the elements of order $4$ in $\G_{D_8}$.  A computation
     shows that $\G'=\G_{D_{8n}}$.
    \item[(d)] By Lemma~\ref{subgroups}, $\G'\cong A_5$, $\G'\cong\G_{D_{10n}}$ for some
    $n>1$,
    $\G'\cong\PSL_2(\F_q)$,
    or $\G'\cong\PGL_2(\F_q)$ for some finite field $\F_q$.  In the last two cases note
    that $\PSL_2(\F_q)$ and $\PGL_2(\F_q)$ have a diagonal element of order $5$ if and only if
    $5|q-1$.

    If $\G'\cong D_{10n}$ for some $n>1$, then
     $\G'$ has an element of order $5n$ that commutes with the elements of order $5$ in $\G_{D_{10}}$.  A computation
     shows that $\G'=\G_{D_{10n}}$.
    \item[(e)] By Lemma~\ref{subgroups}, $\G'\cong\G_{D_{2n'}}$ for some $n'>n$ with $n|n'$,
$\G'\cong\PSL_2(\F_q)$,
    or $\G'\cong\PGL_2(\F_q)$ for some finite field $\F_q$.   In the last two cases note
    that $\PSL_2(\F_q)$ and $\PGL_2(\F_q)$ have an diagonal element of order $n$ only if
    $n|q-1$.

    If $\G'\cong\G_{D_{2n'}}$ for some $n'>n$ with $n|n'$, then
     $\G'$ has an element of order $n'$ that commutes with the elements of order $n$ in $\G_{D_{2n}}$.  A computation
     shows that $\G'=\G_{D_{2n'}}$.
    \item[(f)] By Lemma~\ref{subgroups}, $\G'\cong S_4$, $\G'\cong A_5$,
    $\G'\cong\PSL_2(\F_q)$,
    or $\G'\cong\PGL_2(\F_q)$ for some finite field $\F_q$.

    Suppose that $\G'\cong S_4$.  Since
    $S_4$ has a unique subgroup isomorphic to $A_4$ we have
    $\G'\subseteq N(\G_{A_4})$.  By Lemma~\ref{normalizersub},
     $N(\G_{A_4})=\G_{S_4}$.  So $\G'=\G_{S_4}$.
\end{enumerate}
\end{proof}

\section{Finite subgroups of the 3-dimensional projective general linear groups}
\indent\par Throughout this section let $F$ be an algebraically
closed field of characteristic $p$ with $p=0$ or
 $p>2$.  Before we begin classifying the finite subgroups of
 $\PGL_3(F)$ we need to prove a few lemmas.
\begin{lem}\label{gpisom}
Let $K$ be a field of characteristic $p>0$.
 Let $\G$ be a finite subgroup of $\PGL_n(K)$. Then $\G$ is
isomorphic to a finite subgroup of $\PGL_n(\F_q)$ (and
$\PSL_n(\F_q)$) for some finite field $\F_q$ of order $q=p^r$ for
some $r>0$.
\end{lem}
\begin{proof}
Since $\G$ is finite,  $\G\subseteq\PGL_n(A)$ for some finitely
generated $\F_{p}$-algebra $A$ where $\F_{p}$ is a finite field of
with $p$ elements. Let $V=\Spec(A)$. Then $V$ gives an affine
variety over $\F_{p}$. Let $\ol{\F_{p}}$ be an algebraic closure of
$\F_{p}$ and let $P\in V(\ol{\F_{p}})$.  Let $M\in\G$ and write
\[M:=\begin{bmatrix}
f_{11}&\ldots&f_{1n}\\
.&\ldots&.\\
f_{n1}&\ldots&f_{nn}\\
\end{bmatrix},\] where $f_{ij}\in A$ for all $i$, $j$ and where the
$f_{ij}(P)$ are not identically zero for all $P\in V(\ol{\F_{p}})$.
Let $M(P)$ denote the element
\[\begin{bmatrix}
f_{11}(P)&\ldots&f_{1n}(P)\\
.&\ldots&.\\
f_{n1}(P)&\ldots&f_{nn}(P)\\
\end{bmatrix}\in\PGL_n(\ol{\F_{p}}).\]   The map $\psi_P\colon\G\to\PGL_n(\ol{\F_{p}})$ given by
$M\mapsto M(P)$ is a homomorphism and is injective if and only if
$M(P)\ne Id$ for all $M\in\G-\{Id\}$.  We show that there exists
$P\in V(\ol{\F_{p}})$ such that $\psi_P$ gives an isomorphism of
$\G$ onto a finite subgroup of $\PGL_n(\F_q)$ for some finite field
$\F_q$.

 Note that if $M\in\G-\{Id\}$ is the image of a diagonal
matrix then $M$ is not in the kernel of $\psi_P$ for any $P$ since
it is the image of a matrix in $\GL_n(\ol{\F_{p}})$.

Let \[\mf{H}:=\{M\in\G\colon M \mbox{ is not the image of a diagonal
matrix}\}.\] For each $M\in\mf{H}$ fix a lift
\[M':=\begin{pmatrix}
f_{11}&\ldots&f_{1n}\\
.&\ldots&.\\
f_{n1}&\ldots&f_{nn}\\
\end{pmatrix}\in \GL_n(A),\] where the
$f_{ij}(P)$ are not identically zero for all $P\in V(\ol{\F_{p}})$.
For each $M'$, choose $f_{ij}\ne 0$ with $i\ne j$, let $f^{(M)}:=
f_{ij}$, and let $f=\prod_{M\in\mf{H}} f^{(M)}$. Then since $f\ne
0$, there must exist $Q\in V(\ol{\F_{p}})$ such that $f(Q)\ne 0$. It
follows that the map $\psi_Q$ gives an isomorphism from $\G$ onto a
finite subgroup of $\PGL_n(\F_{q})$ for some finite field $\F_{q}$.
(Note that $\PGL_n(\ol{\F_p})=\PSL_n(\ol{\F_p})$.)
\end{proof}

Given a finite group $G$ of order $r$ and a field $L$,
let $L[G]$ be the algebra of $G$ over $L$. Recall that a left
$L[G]$-module corresponds to a representation of $G$ over $L$ and
that a simple $L[G]$-module corresponds to an irreducible
representation of $G$ over $L$.  Let $S_L$ be the set of isomorphism
classes of simple $L[G]$-modules.

\begin{lem}\label{finextgl} Let $\tilde{\G}$ be
a finite irreducible subgroup of $\GL_n(F)$ of order $r$.  Let
$E\subset F$ be the field obtained by adjoining the $r^{th}$ roots
of unity to the prime field of $F$.  Then $\G$ is $\GL_n(F)$
conjugate to a finite subgroup of $\GL_n(E)$.
\end{lem}
\begin{proof}  It is shown in~\cite{Brauer:2} that an irreducible
representation of a finite group $G$ of order $r$ can be written in
the field of the $r^{th}$ roots of unity.  Since any finite
irreducible subgroup of $\GL_n(F)$ is the image of an irreducible
representation $G\to\GL_n(F)$, our result follows.
\end{proof}
The next two corollaries follow easily from Lemma~\ref{finext}.
\begin{cor} \label{finext} Let $\G$ be a finite irreducible subgroup of $\PGL_n(F)=\PSL_n(F)$.
Then $\G$ is $\PGL_n(F)$-conjugate to a finite subgroup of
$\PGL_n(E)$ and $\PSL_n(E)$ where $E\subset F$ is a finite extension
of the prime subfield of $F$.
\end{cor}
\begin{proof}  This is clear by Lemma~\ref{finextgl}.
\end{proof}
\begin{cor}\label{prtop} Let $\G$ be a finite subgroup of $\PGL_n(F)$ of order $g$.  Let $p$ be the
characteristic of $F$.  Suppose that $p$ is zero or $p\nmid g$. Then
$\G$ is $\PGL_n(F)$-conjugate to a finite subgroup of $\PGL_n(E)$
where $E\subset F$ is a finite extension of the prime subfield of
$F$.
\end{cor}
\begin{proof}  If $p\nmid r$, then by Theorem~1.2 of~\cite{Lang}, the group algebra $F[G]$ is
semisimple.  So every $F[G]$-module is a direct sum of irreducible
modules.  Our result follows easily.
\end{proof}

\begin{lem}\label{repbijection} Let $p$ be a prime.  Let $G$ be a finite group of order $r$.
Let $K$ be a field of characteristic $0$. Assume further that $K$ is
complete with respect to a discrete valuation $\nu$ with valuation
ring $A$, maximal ideal $\mf{m}$, and residue field $k:=A/\mf{m}$ of
characteristic $p>0$. Any representation $G\to\GL_n(K)$ is
isomorphic to a representation $G\to\GL_n(A)\subset\GL_n(K)$.
Furthermore if $p\nmid r$ then the operation of reduction
$\pmod{\mf{m}}$ defines a bijection from $S_K$ onto $S_k$.
\end{lem}
\begin{proof} This follows from Proposition~43 and the remark on page~128
of~\cite{Serre:1}.
\end{proof}
\begin{cor}\label{glbij} Two finite subgroups $\G_1$, $\G_2\subseteq\GL_n(k)$ of order prime to $p$ are
conjugate if and only if there exists two $\GL_n(K)$ conjugate
subgroups $\G_1'$, $\G_2'\subseteq\GL_n(A)$ with $\G_i'\to\G_i$
under the reduction $\pmod{\mf{m}}$ map.
\end{cor}
\begin{proof} This is clear by Lemma~\ref{repbijection}.
\end{proof}

An element of $\PGL_3(F)$ which is of the form
\[\begin{bmatrix}
a & b & 0\\
c & d & 0\\
0 & 0 & 1\\
\end{bmatrix}\]
is called \emph{intransitive}.  A subgroup of $\PGL_3(F)$ consisting
entirely of intransitive elements is also called intransitive.  If a
subgroup $\mf{I}$ is intransitive then it is isomorphic to a
subgroup of $\GL_2(F)$, and there is a natural map from $\mf{I}$
onto a subgroup $\ol{\mf{I}}\subset \PGL_2(F)$.

\begin{lem}\label{gps}
Any finite subgroup $\G$ of $\PGL_3(F)$ is conjugate to one of the
following groups:
\begin{enumerate}
\item[Case] I: when $p=0$ or $|\G|$ is relatively prime to $p$.
\begin{enumerate}
    \item an intransitive group $\mf{I}$ whose image $\ol {\mf{I}}$ in $\PGL_2(F)$
    is equal to one of the groups in
    Lemma~\ref{subgroups} Case~I,
    \item a group generated by
    \[T:=\begin{bmatrix}
    0 & 1 & 0\\
    0 & 0 & 1\\
    1 & 0 & 0\\
    \end{bmatrix}\] and $\mf{H}$, where $\mf{H}$ is
    a finite group generated by the image in
    $\PGL_3(F)$ of diagonal
    matrices.  Such a group will be called a group of type $\mf{C}$.
     Let
     \[S:=\begin{bmatrix}
     1&0&0\\
     0&\omega&0\\
     0&0&\omega^2\\
     \end{bmatrix}\]where $\omega^3=1$.  The group of order $9$ generated by $S$ and $T$ will be
     called $\G_9$.

    \item a group generated by
    \[R:=\begin{bmatrix}
    1&0&0\\
    0&0&1\\
    0&1&0\\
    \end{bmatrix}\] and a group of type $\mf{C}$.  Such a group will be called a group of type
    $\mf{D}$.  Following the notation of~(b), the group of order $18$ generated by $S$, $T$, and
    $R$ will be called $\G_{18}$.
     \item the group $\G_{36}$ of order $36$
    generated by $\G_{18}$ and
    \[V:=\begin{bmatrix}
    1 & 1 & 1\\
    1 & \omega & \omega^2\\
    1 & \omega^2 & \omega\\
    \end{bmatrix},\]
where $\omega$ is a primitive cube root of unity.
    \item the group $\G_{72}$ of order $72$ generated by $\G_{36}$ and $UVU^{-1}$, where
    \[U:=\begin{bmatrix}
    1 & 0 & 0\\
    0 & 1 & 0\\
    0 & 0 & \omega\\
    \end{bmatrix},\] where $\omega$ is a primitive cube root of unity.
    \item the group $\G_{216}$ of order $216$ generated by $\G_{72}$ and $U$ of (e),
    \item the group $\G_{60}$ of order $60$ generated by
  \[E_1:=\begin{bmatrix}
1 & 0 & 0\\
0 & \epsilon^4 & 0\\
0 & 0 & \epsilon\\
\end{bmatrix}, E_2:=\begin{bmatrix}
1 & 0 & 0\\
0 & 0 & 1\\
0 & 1 & 0\\
\end{bmatrix}, \mbox{ and }E_3:=\begin{bmatrix}
1 & 1 & 1\\
2 & \epsilon^2+\epsilon^{-2} & \epsilon+\epsilon^{-1}\\
2 & \epsilon+\epsilon^{-1} & \epsilon^2+\epsilon^{-2}\\
\end{bmatrix}\]
and where $\epsilon$ is a primitive $5^{th}$ root of unity.
  \item[(h)] the group $\G_{360}$ of order $360$ generated by $E_1, E_2, E_3$ of
  (g) and
  \[E_4:=\begin{bmatrix}
1 & \lambda_1 & \lambda_1\\
2\lambda_2 & \epsilon^2+\epsilon^{-2} & \epsilon+\epsilon^{-1}\\
2\lambda_2 & \epsilon+\epsilon^{-1} & \epsilon^2+\epsilon^{-2}\\
\end{bmatrix},\]
where $\lambda_1=\frac{1}{4}(-1+\sqrt{-15})$ and
$\lambda_2=\frac{1}{4}(-1-\sqrt{-15})$.
  \item[(i)] the group $\G_{168}$ of order $168$ generated by
  \[F_1:=\begin{bmatrix}
\beta & 0 & 0\\
0 & \beta^2 & 0\\
0 & 0 & \beta^4\\
\end{bmatrix}, F_2:=\begin{bmatrix}
0 & 1 & 0\\
0 & 0 & 1\\
1 & 0 & 0\\
\end{bmatrix}, \mbox{ and }F_3:=\begin{bmatrix}
a & b & c\\
b & c & a\\
c & a & b\\
\end{bmatrix},\]
where $\beta$ is a primitive $7^{th}$ root of unity,
$a=\beta^4-\beta^3$, $b=\beta^2-\beta^5$, and $c=\beta-\beta^6$.
\end{enumerate}
\item [Case] II: when $|\G|$ is divisible by $p$.
\begin{enumerate}
  \item[(j)]  $\PSL_3(\mathbb{F}_q)$
  \item[(k)]  $\PGL_3(\mathbb{F}_q)$
  \item[(l)] $\PSU_3(\mathbb{F}_q)$
  \item[(m)]  $\PGU_3(\mathbb{F}_q)$
  \item[(n)]  $\G_{\PSL_2(q)}$ which is defined as the image of $\PSL_2(\mathbb{F}_q)$ under the injective map
  \[\begin{bmatrix}
  a&b\\
  c&d\\
  \end{bmatrix}\mapsto
  \begin{bmatrix}
  a^2& ab & b^2\\
  2ac & ad+bc&2bd\\
  c^2 &cd&d^2\\
  \end{bmatrix}.\]

  \item[(o)] $\G_{\PGL_2(q)}$ which is defined as the image of $\PGL_2(\mathbb{F}_q)$ under the map from~(n).
  \item[(p)]  $\G_{A_6}:=\la T, V, B\ra$, $\G_{A_7}:=\la T, V, B, W\ra$, or the group
  $N(\G_{A_6}):= \la T, V, B, U\ra$ containing $\G_{A_6}$ with index $2$,
   and with
  \[T:=\begin{bmatrix}
  -1&0&0\\
  0&-1&0\\
  0&0&1\\
  \end{bmatrix},\]
  \[V:=\begin{bmatrix}
  1&0&0\\
  0&0&1\\
  0&1&0\\
  \end{bmatrix},\]
   \[B:=\begin{bmatrix}
  1&-\alpha-1&-\alpha+1\\
  -\alpha+1&2&2\alpha+2\\
  -\alpha-1&3\alpha+2&2\\
  \end{bmatrix},\]
  \[W:=
\begin{bmatrix}
1 & 0 & 0\\
0 & \omega^2 & 0\\
0 & 0 & \omega\\
\end{bmatrix},\] and

  \[U:=\begin{bmatrix}
1&0&0\\
0&\alpha&\alpha\\
0&\alpha&-\alpha\\
\end{bmatrix}.\]

where $\alpha^2=2$ and $\omega=3\alpha+2$.  In all cases
char$(F)=5$.  The group $\G_{A_6}$ is independent of our choice of
$\alpha$.   That is, the generators obtained by replacing $\alpha$
with $-\alpha$ also generate $\G_{A_6}$.
\item[(q)] $\G_{60}$ of~(g), $\G_{168}$ of~(i), a group of type $\mf{C}$ from~(b), or a group
of type $\mf{D}$ from~(c).  In all cases char$(F)=3$.  For
convenience,
 in characteristic $3$, we will refer to $\G_{60}$ as $\G_{60}^3$, $\G_{168}$ as $\G_{168}^3$, a group
of type $\mf{C}$ as  a group of type $\mf{C}^3$, and a group of type
$\mf{D}$ as  a group of type $\mf{D}^3$.
\item [(r)] a semidirect product $\mf{P}_{A_i}\rtimes \mf{I}$
where, for $i\in\{1,2\}$, $\mf{P}_{A_i}$ is an elementary abelian
$p$-group consisting entirely of elements of the form
\begin{equation}\label{a1}\tag{$A_1$}
\begin{bmatrix}
1&0&0\\
0&1&0\\
\alpha&\beta&1\\
\end{bmatrix}
\end{equation}
or
\begin{equation}\label{a2}\tag{$A_2$}
\begin{bmatrix}
1&0&\alpha\\
0&1&\beta\\
0&0&1\\
\end{bmatrix}
\end{equation}
and where $\mf{I}$ is an intransitive group whose image
$\ol{\mf{I}}\subset\PGL_2(F)$ is one of the groups listed in
Lemma~\ref{subgroups} that is not equal to a group given in
Lemma~\ref{subgroups}(a) or~(f).  If the order of $\mf{I}$ is prime
to $p$ then $\mf{P}_{A_i}$ is nontrivial.

\item [(s)] A semidirect product $\mathfrak{P}\rtimes \mathfrak{D}$ where $\mathfrak{P}$ is a nontrivial $p$-group
consisting entirely of elements of the form
\[\begin{bmatrix}
1&0&0\\
\alpha&1&0\\
\beta&\gamma&1\\
\end{bmatrix}\] and where $\mathfrak{D}$ is a finite diagonal subgroup of $\PGL_3(F)$

\end{enumerate}
where $\F_q$ is the finite field with $q:=p^r$ elements, with $r>0$.
\end{enumerate}
\end{lem}
\begin{proof}\hfill
\begin{enumerate}
  \item [Case] I:
The finite subgroups of $\PGL_3(F)$, where $p=0$, are classified in
Chapter~VII of~\cite{Miller:1}.   Using Corollary~\ref{glbij}, we
obtain a classification for the finite subgroups of $\PGL_3(F)$ in
any characteristic where $p\nmid |\G|$.
\begin{enumerate}
    \item [(a)]  See pages~206-207 and page~236 of~\cite{Miller:1}.
    \item [(b)] See Lemma~108 on page~230 and page~236 of~\cite{Miller:1}.
    \item [(c)] By~\S112 on page~236 of~\cite{Miller:1}, any group not conjugate to any of the other
    groups in this Lemma is conjugate to a group $\G'$ generated by
    a group of type $\mf{C}$, and an element $Q$ of the form
    \[\begin{bmatrix}
    1 & 0 & 0\\
    0 & 0 & a\\
    0 & b & 0\\
    \end{bmatrix}.\]   We need to show that $\G'$ is $\PGL_3(F)$ conjugate to a group generated by a group of type
    $\mf{C}$ and $R$.  Observe that
    \[Z:=(QTQ^{-1}T)(QT^2Q^{-1}T^2)\in \G'.\]
    A computation shows that
    \[Z=\begin{bmatrix}
    1 & 0 & 0\\
    0 & a^3 & 0\\
    0 & 0 & b^3\\
    \end{bmatrix}.\]
Then
    \[Z^{-1}Q^3=\begin{bmatrix}
    1 & 0 & 0\\
    0 & 0 & b/a\\
    0 & a/b & 0\\
    \end{bmatrix}.\]
    Let
    \[Y:=\begin{bmatrix}
    1 & 0 & 0\\
    0 & b & 0\\
    0 & 0 & a\\
    \end{bmatrix}.\]
    A computation shows that
    $YTY^{-1}=Q^2TQ^{-1}TQ^{-1}T$.  So $YTY^{-1}$ is in $\G'$.  Since $(YTY^{-1})T^{-1}$ is the
    image of a diagonal matrix we must have $YTY^{-1}=DT$ for some $D\in\mf{H}$.
    Observe that $YRY^{-1}=Z^{-1}Q^3$.  Note that $Z^{-1}Q^2\in\mf{H}$.  Since $\mf{H}$ is
    generated by the images of diagonal matrices and since diagonal
    matrices commute we have $Y^{-1}\mf{H}Y=\mf{H}$.  So we have \[\G'=\la\mf{H}, T, Q\ra=
    \la\mf{H}, DT, (Z^{-1}Q^2)Q\ra=\la \mf{H}, YTY^{-1},
    YRY^{-1}\ra.\]
     It follows that $\G=Y^{-1}\G'Y$ is of the desired form.
    \item [(d)]and (e)  See pages 236-239 of~\cite{Miller:1}.
    \item [(f)]  See pages 236-239 of~\cite{Miller:1} and page~217 of~\cite{Mitchell:1}.
    \item [(g)]  See pages 250-252 of~\cite{Miller:1} and page~224 of~\cite{Mitchell:1}.
    \item [(h)]  See pages 250-252 of~\cite{Miller:1} and page~225 of~\cite{Mitchell:1}.
    \item [(i)]  See pages 250-251 of~\cite{Miller:1} and \S131 of~\cite{Weber:2}.
    \end{enumerate}
\item[Case] II:  The subgroups of $\PSL_3(\F_q)$, for a finite field
$\F_q$,
 are classified
in~\cite{Bloom:1}.   By Lemma~\ref{gpisom}, any finite subgroup of
$\PGL_3(F)$ is isomorphic to a subgroup of $\PSL_3(\F_q)$ for some
finite field $\F_q$. By Corollary~\ref{finext}, any finite
irreducible subgroup of $\PGL_3(F)$ is $\PGL_3(F)$-conjugate to a
subgroup of $\PSL_3(\F_q)$.  Note that the groups listed in (j)-(q)
are all irreducible.
\begin{enumerate}
  \item [(j)]  See part~(1) of Theorem~1.1 of~\cite{Bloom:1}.
  \item [(k)]  See part~(3) of Theorem~1.1 and Lemma~6.1 of~\cite{Bloom:1}.
  \item [(l)]  See part~(2) of Theorem~1.1 of~\cite{Bloom:1}.
  \item [(m)]  See part~(4) of Theorem~1.1 and Lemma~6.2 of~\cite{Bloom:1}.
  \item [(n)] and~(o)  See part~(5) of Theorem~1.1 and Lemma~6.3
  of~\cite{Bloom:1}.  Note that the map given in Lemma~6.3 of~\cite{Bloom:1} is
  defined incorrectly.  The map we give comes from the symmetric square representation
  of Chapter~1~\S6 of~\cite{Serre:1}.
  \item [(p)]   Suppose that the characteristic of $F$ is $5$.  Let
  $\G'_{A_6}:=W\G_{A_6}W^{-1}=\la T', S', V'\ra$, where
    $T':=WTW^{-1}$, $V':=WVW^{-1}$, and $B':=WBW^{-1}$.  We have $T'=T$,
\[V':=
\begin{bmatrix}
1&0&0\\
0&0&\omega\\
0&\omega^2&0\\
\end{bmatrix},\]
and
\[B':=
\begin{bmatrix}
2&4&1\\
1&4&2\\
4&2&4\\
\end{bmatrix}.\]
 By Lemma~6.6
of~\cite{Bloom:1}, any subgroup of $\PSL_3(F)$ isomorphic to $A_6$
is conjugate to $\G'_{A_6}$.

  If $M\in\PGL_3(\F_{25})\subset\PGL_3(F)$, let
$M^{\sigma}$ be the element of $\PGL_3(\F_{25})$ obtained by
applying the map $a\mapsto a^5$ to the entries of $M$.  We need to
show that $\G_{A_6}=\G_{A_6}^{\sigma}$.

Write
\[X:=(V'T')^2=\begin{bmatrix}
1&0&0\\
0&-1&0\\
0&0&-1\\
\end{bmatrix}.\]
So $X\in \G'_{A_6}$.  Let
\[R:=
\begin{bmatrix}
1&0&0\\
0&0&1\\
0&1&0\\
\end{bmatrix}.\]  Then $RT'R^{-1}=TX$, $RXR^{-1}=X$, $RV'R^{-1}=(V')^{\sigma}$, and
$RB'R^{-1}=(B')^{-1}$.  It follows that
$(\G'_{A_6})^{\sigma}=R\G'_{A_6}R^{-1}$.  Note that
\[W(W^{-1})^{\sigma}R=V'\in \G'_{A_6}.\]
So
\[(W(W^{-1})^{\sigma}R)\G'_{A_6}(W(W^{-1})^{\sigma}R)^{-1}=\G'_{A_6}\]
$\Leftrightarrow$
\[(W^{-1})^{\sigma}(R\G'_{A_6}R^{-1})W^{\sigma}=W^{-1}\G'_{A_6}W \]
$\Leftrightarrow$
\[(W^{-1}\G'_{A_6}W)^{\sigma}=W^{-1}\G'_{A_6}W\]
$\Leftrightarrow$
\[\G_{A_6}=\G_{A_6}^{\sigma}.\]

By Lemma~6.6 of~\cite{Bloom:1}, $\G'_{A_6}$ has index $2$ in
$N(\G'_{A_6})$ and every subgroup isomorphic to $N(\G'_{A_6})$ is
conjugate to $N(\G'_{A_6})$ and is generated by $WUW^{-1}$ and
$\G'_{A_6}$.

By Lemma~6.6 of~\cite{Bloom:1}, every subgroup isomorphic to $A_7$
in $\PSL_3(F)$, is conjugate to $\G_{A_7}$.

  \item [(q)]  By parts~(6) and~(7) of Theorem~1.1
  of~\cite{Bloom:1}, there exists one conjugacy class each of groups isomorphic
  to $\G_{60}$ and $\G_{168}$.  One can verify directly that the generators for $\G^3_{60}$ and
  $\G^3_{168}$ satisfy the necessary relations as generators for the groups.  For the
  last two types of groups, see parts~(1) and~(2) of Theorem~7.1 of~\cite{Bloom:1}.

  \item [(r)] and~(s) Let $\G$ be a finite subgroup of $\PGL_3(F)$ that is not $\PGL_3(F)$-conjugate
  to any of the groups listed earlier in this Lemma.
For any subgroup $\mf{H}\subseteq\PGL_3(F)$ let $\psi$ be the
isomorphism given by $M\mapsto (M^{-1})^{Tr}$, where
 $(M^{-1})^{Tr}$ is the image in $\PGL_3(F)$ of the inverse transpose
of any lift of $M\in\GL_3(F)$.
  It follows from
   the argument given in \S7 of~\cite{Bloom:1} and from Theorem~7.1 of~\cite{Bloom:1}, that up
   to conjugation and/or isomorphism by $\psi$, the group $\G$
   consists entirely of elements of the form
\begin{equation}\label{I}\tag{I}
\begin{bmatrix}
a&b&0\\
c&d&0\\
e&f&1\\
\end{bmatrix}
\end{equation}  and that $\G$ contains an elementary abelian subgroup $\mf{P}_A$
consists entirely of elements of the form
\[\begin{bmatrix}
1&0&0\\
0&1&0\\
\alpha&\beta&1\\
\end{bmatrix}.\]
 The subgroup $\mf{P}_A$ is the kernel of the
homomorphism $\G\to \tilde{\G}'\subset\GL_2(F)$ given by
\[\begin{bmatrix}
a&b&0\\
c&d&0\\
e&f&1\\
\end{bmatrix}\mapsto \begin{pmatrix}
a & b\\
c&d\\
\end{pmatrix}.\]  After conjugating by an element of the form (I) we
may assume that the image $\G'$ of $\tilde{\G}'$ in $\PGL_2(F)$ is
one of the subgroups listed in Lemma~\ref{subgroups}.

First suppose that $\G'$ is one of the groups given in
Lemma~\ref{subgroups}(a) or~(f).  Conjugating $\G$ by the element
\[\begin{bmatrix}
0&1&0\\
1&0&0\\
0&0&1\\
\end{bmatrix}\]
we obtain a group that consists of elements of the form
\[\begin{bmatrix}
\zeta_1&0&0\\
\alpha&\zeta_2&0\\
\beta&\gamma&1\\
\end{bmatrix}.\]Let $\mf{P}$ be the Sylow-$p$ group of $\G$.   The group $\mf{P}$ consists entirely of
elements of the form
\[\begin{bmatrix}
  1&0&0\\
  \alpha&1&0\\
  \beta&\gamma&1\\
  \end{bmatrix}.\]
 Define $\varphi\colon\G\to\PGL_3(F)$ by
\[\begin{bmatrix}
\zeta_1&0&0\\
  \alpha&\zeta_2&0\\
  \beta&\gamma&1\\
  \end{bmatrix}\mapsto
  \begin{bmatrix}
\zeta_1&0&0\\
  0&\zeta_2&0\\
  0&0&1\\
  \end{bmatrix}\] and let $\mf{H}$ be the image of $\varphi$.  Then the
  sequence
  \[1\to\mf{P}\to\G\to\mf{H}\to 1\] is split exact by Hall's
  Theorem.  Let $\phi\colon\mf{H}\to\G$ be a section that is a
  homomorphism.  It can be deduced from Theorem~$8$ of~\S96 of~\cite{Miller:1} that
   $\phi(\mf{H})$ is conjugate to
  subgroup $\mf{D}$, generated by the images of diagonal matrices, via a lower triangular
  matrix $M$.  The group $M\mf{P}M^{-1}$ consists entirely of
  elements of the form
  \[\begin{bmatrix}
  1&0&0\\
  \alpha&1&0\\
  \beta&\gamma&1\\
  \end{bmatrix}.\]  So we obtain the a group of the type given
  in~(s).

{}From now on we will assume that $\G'$ is one of the groups given
in Lemma~\ref{subgroups} not equal to a group given by
Lemma~\ref{subgroups}(a) or~(f).

If the order of $\G'$ is prime to $p$, then by Hall's Theorem the
sequence
\[1\to\mf{P}_A\to\G\to\tilde{\G}'\to 1\] is split exact.
 Let $\psi\colon \tilde{\G}'\to \G$ be a section that is a homomorphism.
  Then by part~(a) of this lemma, $\psi(\tilde{\G}')$ is $\PGL_3(F)$-conjugate to
  and intransitive group $\mf{I}$.   Let $M\in\PGL_3(F)$ and suppose that
  $M\psi(\tilde{\G}')M^{-1}=\mf{I}$.  Since $\G'$ is one of the groups given
in Lemma~\ref{subgroups} not equal to a group given by
Lemma~\ref{subgroups}(a) or~(f), and since $\psi(\tilde{\G}')$
consists of elements of the form~\eqref{I}, it can be verified that
$\psi(\tilde{\G}')$ fixes exactly one point of $\mathbb{P}^2(F)$,
the point $P:=[0\colon 0\colon 1]$.  It can be verified that $P$ is
the only point of $\mathbb{P}^2(F)$ fixed by $\mf{I}$.  It follows
that $M$ must also fix $P$.  So $M$ must be of the form~\eqref{I}.
Then $M\mf{P}_AM^{-1}$ is of the form~\eqref{a1}.  So $\G$ is
$\PGL_3(F)$ conjugate to a group of the desired form.

If $p$ divides the order of $\G'$ then $\G'$ is one of the groups
given in Lemma~\ref{subgroups}(g) or~(h), then by Theorem~3.4(5)
of~\cite{Bloom:1}, $\tilde{\G}'$ contains the element
\[M':=\begin{pmatrix}
-1&0\\
0&-1\\
\end{pmatrix}.\]
After conjugating $\G$ by an element of the form~\eqref{I}, we may
assume that
\[M:=\begin{bmatrix}
-1&0&0\\
0&-1&0\\
0&0&1\\
\end{bmatrix}\in\G.\]  Then for any
\[H:=\begin{bmatrix}
a&b&0\\
c&d&0\\
e&f&1\\
\end{bmatrix}\in\G,\]
the element
\[L:=(M^{-1}HMH^{-1})^{p-1}
\]
\[= \begin{bmatrix}
1&0&0\\
0&1&0\\
e&f&1\\
\end{bmatrix}\in\G.\]  Then $H=WL$ where
\[W:=\begin{bmatrix}
a&b&0\\
c&d&0\\
0&0&1\\
\end{bmatrix}\in\G.\]   So $\G$ is of the form $\mf{P}_{A_1}\rtimes\mf{I}$ and it is easily seen that
 the transpose inverse of
$\G$ is of the form $\mf{P}_{A_2}\rtimes\mf{I}$.
\end{enumerate}
\end{enumerate}
\end{proof}

\begin{lem}\label{normalizergps} Let $\G$ be one of the groups listed in Lemma~\ref{gps} and let
$N(\G)$ be the normalizer of $\G$ in $\PGL_3(F)$.  Then, using the
notation of Lemma~\ref{gps},
\begin{enumerate}
  \item[(b)] suppose that $\G$ is a group of type $\mf{C}$.  So $\G$ is generated by the element $T$ and
   a diagonal
  subgroup $\mf{H}$.  Let $m=3^l$ be the largest power of $3$ such that $\mf{H}$ has an element
    of order $m$ and let $\zeta_m$ be a
    primitive $m^{th}$ root of unity.
    Let \[S_{m,2}:=\begin{bmatrix}
      \zeta_m & 0 & 0\\
      0  & \zeta_m^{-1} & 0\\
      0  & 0  & 1\\
    \end{bmatrix}\] and let
    \[S_{m,1}:=\begin{bmatrix}
      \zeta_m & 0 & 0\\
      0  & 1 & 0\\
      0  & 0  & 1\\
    \end{bmatrix}.\]Then
  \begin{enumerate}
    \item[i.] if $\G=\G_9$ then $N(\G)=\G_{216}=\G_9\rtimes\la US^2, V\ra$.
    \item[ii.] if $S_{m,1}\in\mf{H}$ and if $\mf{H}\ne\la S_{3,2}\ra$, then $N(\G)\subseteq\la \G, R, S_{3m,2}\ra$.
    \item[iii.] if $S_{m,1}\not\in\mf{H}$ and if $\mf{H}\ne\la S_{3,2}\ra$, then $N(\G)\subseteq\la \G, R, S_{m,1}\ra$.
  \end{enumerate}
  \item[(c)] suppose that $\G$ is a group of type $\mf{D}$.  So $\G$ is generated by the elements $T$, $R$, and
   a diagonal
  subgroup $\mf{H}$. Then
  \begin{enumerate}
    \item[i.] if $\G=\G_{18}$ then $N(\G)=\G_{216}$.
    \item[ii.] if $S_{m,1}\in\mf{H}$ and if $\mf{H}\ne\la S_{3,2}\ra$, then $N(\G)\subseteq\la \G, S_{3m,2}\ra$.
    \item[iii.] if $S_{m,1}\not\in\mf{H}$ and if $\mf{H}\ne\la S_{3,2}\ra$, then $N(\G)\subseteq\la \G,
    S_{m,1}\ra$.
  \end{enumerate}
  \item[(d)] $N(\G_{36})=\G_{72}$.
  \item[(e)] $N(\G_{72})=\G_{216}=\G_{72}\rtimes\la U\ra$.
  \item[(f)] $N(\G_{216})=\G_{216}$.
  \item[(g)] $N(\G_{60})=\G_{60}$.
  \item[(h)] $N(\G_{360})=\G_{360}$.
  \item[(i)] $N(\G_{168})=\G_{168}$.
  \item[(j)] $N(\PSL_3(\F_q))=\PGL_3(\F_q)$.  Note that if $q\not\equiv 1\pmod{3}$, then
  $\PSL_3(\F_q)=\PGL_3(\F_q)$.  Otherwise, $\PGL_3(\F_q)=\la\PSL_3(\F_q), M\ra$ where
  \[M:=\begin{bmatrix}
  \alpha&0&0\\
  0&1&0\\
  0&0&1\\
  \end{bmatrix}\] for some $\alpha$ where $\alpha\in (\F_q)^3-\F_q$.
  \item[(k)] $N(\PGL_3(\F_q))=\PGL_3(\F_q)$.
  \item[(l)]$N(\PSU_3(\F_q))=\PGU_3(\F_q)$.   Note that if $q\equiv 0\pmod{3}$, then
$\PSU_3(\F_q))=\PGU_3(\F_q)$.  Otherwise
$\PGU_3(\F_q)=\la\PSU_3(\F_q), M\ra$
  where
  \[M:=\begin{bmatrix}
  \alpha&0&0\\
  0&1&0\\
  0&0&1\\
  \end{bmatrix}\] for some $\alpha$
   where $\alpha\in(\F_q)^3-\F_q$ if $q\equiv 1\pmod{3}$ and $\alpha\in(\F_{q^2})^{q-1}-(\F_{q^2})^{3(q-1)}$
if $q\equiv 2\pmod{3}$.
  \item[(m)] $N(\PGU_3(\F_q))=\PGU_3(\F_q)$.
  \item[(n)] $N(\G_{\PSL_2(q)})=\G_{\PGL_2(q)}=\la \G_{\PSL_2(q)}, M\ra$
  where
  \[M:=\begin{bmatrix}
  \alpha^2&0&0\\
  0&\alpha&0\\
  0&0&1\\
  \end{bmatrix}\] for any $\alpha\in\F_q$ that is not a square.
  \item[(o)] $N(\G_{\PGL_2(q)})=\G_{\PGL_2(q)}$.
  \item[(p)] $N(N(\G_{A_6}))=N(\G_{A_6})$ and
  $N(\G_{A_7})=\G_{A_7}$.
  \item[(q)] $N(\G^3_{60})=\G^3_{60}$, $N(\G^3_{168})=\G^3_{168}$. If $\G$ is of type $\mf{C}^3$
 then $N(\G)=\G\rtimes\mf{K}$
    where $\mf{K}\le\la  R\ra$.  If $\G$ is of type $\mf{D}^3$
 then $N(\G)=\G$.

\end{enumerate}
\end{lem}
\begin{proof}\hfill
\begin{enumerate}
               \item[(b)] It can be deduced by Lemma~1.5 of~\cite{Conlon:1} that any group that contains
    a diagonal element of order
    $m$ contains the element $S_{m,2}$.
               \begin{enumerate}
                 \item[i.] Suppose that $\G=\G_9$.  By \S115 of~\cite{Miller:1},
                 $N(\G_9)=\G_{216}$ and $\G_{216}$ is generated by
                 $\G_9$, $U$, and $V$.  Note that $\la\G_9, U, V\ra=\la\G_9, US^2, V\ra$.
                  One can verify that
                  $\la US^2, V\ra\cap\G_9=Id$.
                   \item[ii.] Suppose that
                 $S_{m,1}\in\mf{H}$ and $\mf{H}\ne\la S_{3,2}\ra$.  By the Lemma of~\S108
 of~\cite{Miller:1}, $N(\G)$ is a group generated by $R$, $T$, and the
 images in $\PGL_3(F)$ of diagonal matrices.  Let $D$
 be the image in $\PGL_3(F)$ of a diagonal matrix.  A
 computation shows that if $DTD^{-1}\in\G$ then
 $D^3\in\G$.  So if $D\in N(\G)$, then $D=D'M$, where $D'\in\mf{H}$
  and $M$ is the image of a diagonal matrix of order a power of $3$.
   By Lemma~1.5 of~\cite{Conlon:1} we may assume that $M$ is a power of $M=S_{3m,1}$ or
   $S_{3m,2}$.  A computation shows that $S_{3m,1}\not\in
   N(\G)$ and that $S_{3m,2}\in N(\G)$.
   It follows
 that $N(\G)\subseteq \la \G, S_{3m,2},
 R\ra$.
                 \item[iii.] Suppose that $S_{m,1}\not\in\mf{H}$ and  $\mf{H}\ne\la S_{3,2}\ra$.
                    The proof is
                 the same as for $N(\G)$ in (ii), except that $M$ must be a
                 power of $S_{m,1}$ or $S_{3m,2}$ and the final
                 computation shows that $S_{3m,2}\not\in
   N(\G)$ and that $S_{m,1}\in N(\G)$.

               \end{enumerate}
               \item[(c)]
               \begin{enumerate}
                 \item[i]Suppose that $\G=\G_{18}$.
                  By \S115 of~\cite{Miller:1}, $N(\G_{18})=\G_{216}$.
                 \item[ii.]and~iii.  By assumption $\G$ is generated by $R$ and a group $\mf{G}'$ of type $\mf{C}$.
                   It is easily verified
                 that $\mf{G}'$ is a characteristic subgroup of
                 $\G$.  It follows that
                 $N(\G)\subseteq N(\mf{G}')$.  The
                 result follows from~(b).
               \end{enumerate}
               \item[(d)] By \S115 of~\cite{Miller:1}, $N(\G_{36})\subseteq
    \G_{216}$.  Since $UVU^{-1}\not\in \G_{36}$, $N(\G_{36})\ne \G_{216}$.  Since $\G_{36}$ has index
$2$ in $\G_{72}$ it is normal in $\G_{72}$.  Since $\G_{72}$ has
prime index in $\G_{216}$, we must have $N(\G_{36})=\G_{72}$.
               \item[(e)] By \S115 of~\cite{Miller:1},
$N(\G_{72})\subseteq
    \G_{216}$.  A computation shows that $\G_{72}$ is closed under conjugation by $\G_{216}$,
     so $N(\G_{72})=\G_{216}$.  Since $U \in \G_{216} - \G_{72}$ has order $3$, $\G_{216}=\G_{72}\rtimes\la
     U\ra$.
    \item [(f)]  By \S115 of~\cite{Miller:1}, $N(\G_{216})=
    \G_{216}$.
    \item [(g)]  By \S124 of~\cite{Miller:1}, $N(\G_{60})=
    \G_{60}$.
    \item [(h)]   By \S124 of~\cite{Miller:1}, $N(\G_{360})=
    \G_{360}$.

    \item [(i)]  By \S124 of~\cite{Miller:1}, $N(\G_{168})=
    \G_{168}$.
    \item [(j)]and~(l)  Clear by Lemmas~6.1 and ~6.2 of~\cite{Bloom:1}.
    \item [(k)]and~(m) Let $\G$ be equal to $\PGL_3(\F_q)$ or
    $\PGU_3(\F_q)$ and let $\mf{H}$ be equal to $\PSL_3(\F_q)$ or
    $\PSU_3(\F_q)$ respectively.  By Theorem~1.1 of~\cite{Bloom:1},
    $\mf{H}$ is a normal subgroup of $\G$ of index $1$ or $3$.
       We show that $\mf{H}$ is a characteristic
     subgroup of $\G$, so in particular $N(\G)\subseteq N(\mf{H})$.
       Then our result follows from~(j) and~(l).

       Let $\varphi$ be an automorphism of $\G$ and let
       $\mf{H}'=\varphi(\mf{H})$.  Then since
       $\mf{H}$ and $\mf{H}'$ are normal subgroups of $\G$, the group
       $\mf{H}\cap\mf{H}'$ is a normal subgroup of
       $\mf{H}$.  By Theorem~5.14 of~\cite{Bloom:1}, $\mf{H}$ is a simple
     group.  So $\mf{H}\cap\mf{H}'=\{Id\}$ or
       $\mf{H}$.  Suppose that
       $\mf{H}\cap\mf{H}'=\{Id\}$.  So $|\mf{H}\mf{H}'|=|\mf{H}||\mf{H}'|$.
        By page~208
       of~\cite{Mitchell:1}, $|\mf{H}'|=|\mf{H}|>3$.  (Mitchell calls $\PSL_3(\F_q)$
       and $\PSU_3(\F_q)$, LF$(3,q)$ and HO$(3, q^2)$ respectively.)
        Since $\mf{H}\mf{H}'\subseteq\G$, $|\mf{H}||\mf{H}'|\le |\G|=3|\mf{H}|$.
        This is a contradiction.  Therefore $\mf{H}=\mf{H}'$
        and so $\mf{H}$ is a characteristic subgroup of $\G$.
    \item [(n)]  This is clear by Lemma~6.3 of~\cite{Bloom:1}.
    \item [(o)]  This follows from the fact that $\PSL_2(\F_q)$ is a
    characteristic subgroup of $\PGL_2(\F_q)$ and part~(n) of this Lemma.
    \item[(p)] By Lemma~6.6 of~\cite{Bloom:1},
    $N(N(\G_{A_6}))=N(\G_{A_6})$.  By Lemma~6.7 of~\cite{Bloom:1},
    $N(\G_{A_7})=\G_{A_7}$.
    \item [(q)] By Lemma~6.4 of~\cite{Bloom:1},
    $N(\G^3_{60})=\G^3_{60}$.  By Lemma~6.5 of~\cite{Bloom:1},
    $N(\G^3_{168})=\G^3_{168}$.  The last two statement follow from
    simple computations very similar to the ones in the proofs of
    part~(b) and~(c).
             \end{enumerate}
\end{proof}

\begin{lem}\label{diagelt}  Let $\mf{I}\subset\PGL_3(F)$ be a finite intransitive group.   Suppose that the image $\ol{\mf{I}}$ of
$\mf{I}$ in $\PGL_2(F)$ is one of the groups listed in
Lemma~\ref{subgroups} which is not equal to one of the groups given
by Lemma~\ref{subgroups}(a) or~(f).  Then $\mf{I}$ has an element of
the form
 \[\begin{bmatrix}
 \zeta_1&0&0\\
 0&\zeta_2&0\\
 0&0&1\\
 \end{bmatrix}\] where both $\zeta_1$ and $\zeta_2$ are roots of unity
 not equal to $1$.
\end{lem}
\begin{proof}  If $\ol{\mf{I}}$  is one of the groups
given in Lemma~\ref{subgroups}(g) or~(h), then by Theorem~3.4(5)
of~\cite{Bloom:1}, $\mf{I}$ contains the element
\[\begin{bmatrix}
-1&0&0\\
0&-1&0\\
0&0&1\\
\end{bmatrix}.\]  In the other cases, after looking at the generators
for the groups listed in Lemma~\ref{subgroups}, one can deduce that
$\ol{\mf{I}}$ contains a dihedral subgroup generated by
\[R:=\begin{bmatrix}
\zeta&0\\
0&1\\
\end{bmatrix}\] and
\[S:=\begin{bmatrix}
0&\zeta\\
1&0\\
\end{bmatrix}\] where $\zeta$ is a primitive $n^{th}$ root of unity
for some $n>1$.  Then $\mf{I}$ contains the element
\[\tilde{S}:=\begin{bmatrix}
0&\zeta \omega &0 \\
\omega&0&0\\
0&0&1\\
\end{bmatrix}\] where  $\omega$ is a some root of unity.  Then
\[\tilde{S}^2=\begin{bmatrix}
\zeta \omega^2&0 &0 \\
0&\zeta\omega^2&0\\
0&0&1\\
\end{bmatrix}\] is of the desired form unless $\zeta=1/\omega^2$.  Suppose that
$\zeta=1/\omega^2$.  Then $\mf{I}$ contains the element
\[\tilde{R}:=\begin{bmatrix}
 \omega^{-2}\omega'&0 &0 \\
0&\omega'&0\\
0&0&1\\
\end{bmatrix}\] and the element
\[(\tilde{R}\tilde{S})^2=\begin{bmatrix}
(\omega')^2/\omega^2&0&0\\
0&(\omega')^2/\omega^2&0\\
0&0&1\\
\end{bmatrix}\]where  $\omega'$ is a some root of unity.  One of
$\tilde{R}$ or $(\tilde{R}\tilde{S})^2$ is of the desired form, for
if $(\tilde{R}\tilde{S})^2$ is not then $\omega^2=(\omega')^2$ and
so
\[\tilde{R}=\begin{bmatrix}
 (\omega')^{-1}&0 &0 \\
0&\omega'&0\\
0&0&1\\
\end{bmatrix}\] is of the desired form since $(\omega')^{-2}=\zeta$ implies that $\omega'$ is a primitive root
of unity for some $n>2$.
\end{proof}

%% file: ch3.tex
\section{Isomorphisms of hyperelliptic curves}

 Throughout this section let $K$ be a perfect field
of characteristic not equal to $2$, let $F$ be an algebraic closure
of $K$, and let $X$ be a hyperelliptic curve over $F$.  In
particular, $X$ admits a degree-$2$ morphism to ${\mathbb P}_F^1$
and the genus of $X$ is at least $2$.  Each element of Aut($X$)
  induces an automorphism of $\mathbb{P}^1_F$ fixing the branch points.
The number of branch points is $\ge 3$ (in fact $\ge 6$), so
$\Aut(X)$ is finite.  We get a homomorphism $\Aut(X) \to
\Aut({\mathbb P}^1_F) = \PGL_2(F)$ with kernel generated by the
hyperelliptic involution $\iota$. Let $\G\subset \PGL_2(F)$ be the
image of this homomorphism. Replacing the original map $X \to
{\mathbb P}^1_F$ by its composition with an automorphism $g \in
\Aut({\mathbb P}^1_F) = \PGL_2(F)$ has the effect of changing $\G$
to $g\G g^{-1}$, so we may assume that $\G$ is one of the groups
listed in Lemma~\ref{subgroups}.
  Fix an equation $y^2=f(x)$ for $X$
  where $f\in F[x]$ and $\mbox{disc}(f)\ne 0$.  So the function field $F(X)$ equals $F(x,y)$.

\begin{prop}\label{isom}  Let $X'$ be a
hyperelliptic curve over $F$ given by $y^2=f'(x)$, where $f'(x)$ is
another squarefree polynomial in $F[x]$.  Every isomorphism
$\varphi\colon X\to X'$ is given by an expression of the form:

\[ (x,y)\mapsto \left(\frac{ax+b}{cx+d},
\frac{ey}{(cx+d)^{g+1}}\right),\]

\noindent for some $M=\begin{pmatrix}

  a & b \\
c & d \\
\end{pmatrix}\in \GL_2(F)$ and $e\in F^{\times}$.
The pair $(M,e)$ is unique up to replacement by $(\lambda M, e
\lambda^{g+1})$ for $\lambda\in F^{\times}$. If $\varphi'\colon
X'\to X''$ is another isomorphism, given by $(M',e')$, then the
composition $\varphi'\varphi$ is given by $(M'M,e'e)$.

\end{prop}
\begin{proof}
See Proposition 2.1 in~\cite{Poonen:1}.
\end{proof}

Throughout the rest of this section assume that  $K$ is the field of
moduli of $X$ relative to the extension $F/K$ and let
$\Gamma=\Gal(F/K)$.
\begin{lem}\label{isomnormalizer}  Suppose
$\sigma\in\Gamma$ and suppose that the isomorphism $\varphi\colon
X\to \el{\sigma}{X}$ is given by $(M,e)$.  Let $\overline{M}$ be the
image of $M$ in $\PGL_2(F)$.  If $\G\ne \G_{\beta, A}$ then
$\overline{M}$ is in the normalizer $N(\G)$ of $\G$ in $\PGL_2(F)$.
If $\G=\G_{\beta, A}$ then $M$ is an upper triangular matrix.
\end{lem}

\begin{proof}
Since Aut$(\el{\sigma}{X})=\{\psi^{\sigma}\mid\psi\in\Aut(X)\}$, the
group of automorphisms of $\mathbb{P}^1$ induced by
Aut($\el{\sigma}{X})$ is $\G^{\sigma}:=\{U^{\sigma}\mid U\in \G\}$.

Let $\psi$ be an automorphism of $X$ given by $(V,v)$. Since $\psi$
is an automorphism, $V\in \GL_2(F)$ is a lift of some element $\ol
V\in \G$. Then $\varphi\psi\varphi^{-1}$ is an automorphism of
$\el{\sigma}{X}$ given by $(MVM^{-1}, v)$.  We have
$\ol{MVM^{-1}}{=\ol M\;\ol V\;\ol{M}^{-1}}{\in \G^{\sigma}}$. It
follows that $\overline{M}\G\overline{M}^{-1}=\G^{\sigma}.$   If
$\G\ne \G_{\beta, A}$, by Lemma~\ref{subgroups}, $\G^\sigma=\G$. So
$\overline{M} \in N(\G)$.  If $\G=\G_{\beta, A}$, then since
$\G^{\sigma}$ has an elementary abelian subgroup of the same form as
$\G$, a simple computation shows that $M$ is an upper triangular
matrix.
\end{proof}
\begin{lem}\label{ingsig}
Suppose that for every $\tau\in\Gamma$ there exists an isomorphism
$\varphi_{\tau}\colon X\to \el{\tau}{X}$ given by $(M_{\tau},e)$
where $\overline{M}_{\tau}\in \G^{\tau}$. Then $X$ can be defined
over $K$. Furthermore, $X$ is given by an equation of the form
$z^2=h(x)$ where $h\in K[x]$.
\end{lem}
\begin{proof}  Let $P_1,\ldots,P_n$ be the hyperelliptic branch points of $X\to\mathbb{ P}^1$.  Let ${\tau\in\Gamma}$.
The isomorphism $\varphi_{\tau}\colon X\to \el{\tau}{X}$ induces an
isomorphism on the canonical images $\mathbb{P}^1\to\mathbb{P}^1$
which is given by $\overline{M}_{\tau}$.  Write
$\tau(\infty)=\infty$.  The hypothesis $\overline{M}_\tau \in
\G^\tau$ implies that $\overline{M}_\tau$ maps $\{\tau(P_1), \ldots,
\tau(P_n)\}$ to itself; since it also maps $\{P_1,\dots,P_n\}$ to
$\{\tau(P_1), \ldots, \tau(P_n)\}$, we get $\{\tau(P_1), \ldots,
\tau(P_n)\} = \{P_1,\dots,P_n\}$.   So
\[h(x):=\prod_{P_j\ne\infty} (x-P_j)\in K[x].\]
 It follows that $X$ can be defined over $K$.
\end{proof}
\begin{cor}\label{def}Suppose that $N(\G)=\G$ and $\G\ne \G_{\beta, A}$. Then $X$ can be defined over $K$.
\end{cor}
\begin{proof}By Lemma~\ref{subgroups}, $\G^{\sigma}=\G$ for all $\sigma\in\Gamma$.  Let $\tau\in\Gamma$.  By Lemma~\ref{isomnormalizer}, any isomorphism $X\to\el{\tau}{X}$ is
given by $(M,e)$ where $\ol M\in N(\G)=\G=\G^{\tau}$.
\end{proof}

\begin{lem}\label{lift} Suppose there exists an automorphism $\psi$ of
$\mathbb{P}^1$ such that for all $\sigma\in\Gamma$ the automorphism
$(\psi^{-1})^{\sigma}\psi$ lifts to an isomorphism
$\varphi_{\sigma}\colon X\to\el{\sigma}{X}$.  Then $\psi$ lifts to
an isomorphism
\[X\to Y_F,\]
 where $Y$ is a $K$-model of $X$ is given by an equation of the form $z^2=h(x)$ with $h\in K[x]$.
\end{lem}
\begin{proof} Let $P_1,\ldots,P_n$ be the hyperelliptic branch points of $X\to\mathbb{ P}^1$.  By assumption $X$ is the smooth
projective model of $y^2=f(x)$, where $f(x)\in F[x]$.  Let
\[f(x):=\lambda\prod_{P_j\ne\infty}(x-P_j),\]
where $\lambda\in F^{\times}$.

 Suppose $\sigma\in\Gamma$.  Writing $\sigma(\infty)=\infty$, we have
by assumption
\[\{\sigma(P_1),\ldots,\sigma(P_n)\}=\{(\psi^{-1})^{\sigma}\psi(P_1),\ldots,(\psi^{-1})^{\sigma}\psi(P_n)\}.\]
So
\[\{\sigma(\psi(P_1)),\ldots,\sigma(\psi(P_n))\}\]
\[=\{\psi^{\sigma}(\sigma(P_1)),\ldots, \psi^{\sigma}(\sigma(P_n))\}\]
\[=\{\psi^{\sigma}((\psi^{-1})^{\sigma}\psi(P_1)),\ldots, \psi^{\sigma}((\psi^{-1})^{\sigma}\psi(P_n))\}\]
\[=\{\psi(P_1),\ldots,\psi(P_n)\}.\]
 Then
\[h(x):=\prod_{\psi(P_j)\ne\infty} (x-\psi(P_j))\in K[x].\]
Let $Y$ be the hyperelliptic curve over $K$ given by $z^2=h(x)$.
  The hyperelliptic branch points of $Y_F\to\mathbb{P}^1$ are $\{\psi(P_1),\ldots,\psi(P_n)\}$.  Suppose
that $\psi$ is given by
\[x\mapsto \frac{ax+b}{cx+d}.\]
Let
\[r(x):=h\left(\frac{ax+b}{cx+d}\right)(cx+d)^n\in F[x].\]
If $P$ is a zero of $r$, then either $P$ is a zero of $f$ with
$\psi(P)\ne\infty$, or $c\ne 0$ and $P=-\frac{d}{c}$. If $c\ne 0$
and $-\frac{d}{c}$ is a zero of $r$, then the degree of $h$ is
$n-1$.  In this case $\infty$ is a branch point of
$Y_F\to\mathbb{P}^1$ so we must have $\psi(Q)=\infty$ for some
hyperelliptic branch point of $Q$ of $X\to\mathbb{P}^1$.  Then since
$\psi(-\frac{d}{c})=\infty$, $Q=-\frac{d}{c}$ is a zero of $f$. So
$r|f$.

Conversely, let $Q$ be a zero of $f$.  Clearly, if
$\psi(Q)\ne\infty$ then $Q$ is a zero of $r$.  Suppose that
$\psi(Q)=\infty$.  Then the degree of $h$ is $n-1$ and $cx+d$ must
divide $r$.  Also, since $\psi(Q)=\infty$
 we must have $c\ne 0$ and $Q=-\frac{d}{c}$.  So $Q$ is a zero of $r$.  So $f|r$ and we must have $f=\lambda' r$
 for some $\lambda'\in F^{\times}$.  By Proposition~\ref{isom}, it follows that there exists $e\in F^{\times}$ such that
 $(M,e)$ gives an isomorphism $X\to Y_F$ where

 \[M:=\begin{pmatrix}
 a&b\\
 c&d\\
 \end{pmatrix}.\]

\end{proof}

\section{Isomorphisms of plane curves}

Throughout this section let $F$ be an algebraically closed field of
characteristic not equal to $2$.

\begin{thm}\label{planeisom}  Let $X$ and $X'$ be smooth plane curves
of degree $d>3$ over $F$.  Then any isomorphism from $X$ to $X'$ is
induced by a linear transformation of $\mathbb{P}^2$.
\end{thm}
\begin{proof}
See \cite{Chang:1}.
\end{proof}

\begin{lem}\label{isnorm} Let $K$ be a subfield of $F$ such that $F/K$ is
Galois and let $\Gamma:=\Gal(F/K)$. Let $X$ be a smooth plane curve
over $F$ of degree $d> 3$. Let $\G\subset \PGL_3(F)$ be the
automorphism group of $X$ and let $N(\G)$ be the normalizer of $\G$
in $\PGL_3(F)$.  Let $\sigma\in\Gamma$ and suppose that there exists
an isomorphism $\varphi\colon X\to\el{\sigma}{X}$.  If
$\G^{\sigma}=\G$ then $\varphi$ is given by an element $M$ in the
normalizer $N(\G)$ of $\G$.
\end{lem}
\begin{proof}  Let $\sigma\in\Gamma$ and suppose that $\varphi\colon X\to\el{\sigma}{X}$ is an isomorphism.
  By Theorem~\ref{planeisom}, $\varphi$ is given by $M\in\PGL_3(F)$.  Since
  $\Aut(\el{\sigma}{X})=\G^{\sigma}$, we have
   $M\G M^{-1}=\G^{\sigma}=\G$,
  we must have $M\in N(\G)$.
\end{proof}

\begin{lem}\label{conjcurve}   Let $X$ be a smooth plane curve over $F$ and suppose that
the automorphism group $\Aut(X)$ of $X$ is given by $\G$, one of the
groups listed in Lemma~\ref{gps}.  Let $K$ be a subfield of $F$ such
that $F/K$ is Galois and let $\Gamma:=\Gal(F/K)$.   Let
$\sigma\in\Gamma$ and suppose that $M\in\PGL_3(F)$ gives an
isomorphism $X\to\el{\sigma}{X}$. Then $M\in N(\G)$ unless
\begin{enumerate}
    \item [(a)]  $\G=\mf{I}$ is intransitive where $\ol{\mf{I}}$ is one of the groups
     in Lemma~\ref{subgroups}.   If $\ol{\mf{I}}$ is not a group given in Lemma~\ref{subgroups}(a) or~(f),
     then $M$ is an element
    of an intransitive group $\mf{I}'$ containing $\mf{I}$ with $\ol{\mf{I}'}$ equal to the
    normalizer of $\ol {\mf{I}}$ in $\PGL_2(F)$.
    \item [(h)] $\G=\G_{360}$.   Then $M\in \G_{360}$ if and only if $\sigma(\sqrt{-15})=\sqrt{-15}$.  If
    $\sigma(\sqrt{-15})=-\sqrt{-15}$, then $M=M'A$
     where $A\in \G_{360}$ and
\[M':=\begin{bmatrix}
\lambda_2 & 0 &0\\
0& 1&0\\
0&0&1\\
\end{bmatrix}.\]
    \item [(r)] $\G=\mf{P}_{A}\rtimes\mf{I}$ is a group given in
    Lemma~\ref{gps}(r), where $\mf{P}_A$ is an elementary abelian $p$-group and $\mf{I}$ is an intransitive group.
 Then $M=M'A$, where $A\in\G$ and $M'$ is an
element
    of an intransitive group $\mf{I}'$ containing $\mf{I}$ with $\ol{\mf{I}'}$ equal to the
    normalizer of $\ol {\mf{I}}$ in $\PGL_2(F)$.
    \item [(s)] $\G=\mf{P}\rtimes\mf{D}$ is a group given in Lemma~\ref{gps}(s),
    where is $\mf{P}$ a nontrivial $p$-group
consisting entirely of elements of the form
\[\begin{bmatrix}
1&0&0\\
\alpha&1&0\\
\beta&\gamma&1\\
\end{bmatrix}\] and where $\mathfrak{D}$ is a finite diagonal subgroup of $\PGL_3(F)$.

\end{enumerate}
\end{lem}
\begin{proof} In all but the four cases listed, it is easily
verified that $\G^{\sigma}=\G$ and so by Lemma~\ref{isnorm}, $M$ is
in $N(\G)$.  We now consider the other four cases.
\begin{enumerate}
    \item [(a)]  Any intransitive group, whose image in $\PGL_2(F)$ is given by one of the groups of
    Lemma~\ref{subgroups} that is not equal to a group given in~(a) or~(f), fixes exactly one point of
    $\mathbb{P}^2(F)$, the point $P:=[0\colon 0\colon 1]$.
      Since $\mf{I}^{\sigma}=M\mf{I}M^{-1}$ is also intransitive and since $\ol{\mf{I}^{\sigma}}$
       is a group of the same type as $\mf{I}$, $M\mf{I}M^{-1}$ must
      also fix
      only the point $P$.  It follows that $M$ must also fix $P$, so $M$ must be  of the form
      \[\begin{bmatrix}
      a & b &0\\
      c & d &0\\
      e & f &1\\
      \end{bmatrix}.\]
      Let $\psi$ be the inverse transpose isomorphism of $\PGL_3(F)$
      defined in the proof of Lemma~\ref{gps}(r).  Then both $\psi(\mf{I})$ and $\psi(\mf{I}^{\sigma})$
      fix only the point $P$.  So
      \[\psi(M)=\begin{bmatrix}
      d&-c&cf-de\\
      -b&a&-af+be\\
      0&0&ad-bc\\
      \end{bmatrix}\] must also fix $P$.  So $cf-de=-af+be=0$.   This
      implies that
      $(e,f)$ is a scalar multiple of both $(c,d)$ and $(a,b)$.  Since
      $M$ is invertible $ad-bc\ne 0$.  So we must have
      $(e,f)=(0,0)$.
Therefore $M$ must be of the form
      \[\begin{bmatrix}
      a & b &0\\
      c & d &0\\
      0 & 0 &1\\
      \end{bmatrix}.\]

      Since $\ol{\mf{I}^{\sigma}}=\ol {\mf{I}}$, the element $\ol M$ must be in the normalizer of $\ol {\mf{I}}$
      in $\PGL_2(F)$.

    \item [(h)]  If $\sigma(\sqrt{-15})=\sqrt{-15}$, then $\G_{360}^{\sigma}=\G_{360}$ and $M$ must
    be in the normalizer $N(\G_{360})$ of $\G_{360}$.  By
    Lemma~\ref{normalizergps}, $N(\G_{360})=\G_{360}$.

    Suppose that $\sigma(\sqrt{-15})=-\sqrt{-15}$.  Note that $\G_{360}=\la E_1, E_2, E_3, E_4\ra$
    and $\G_{360}^{\sigma}=\la E_1, E_2, E_3, E_4^{\sigma}\ra$.  A computation shows that
    \[M'E_1(M')^{-1}=E_1,\]
    \[ M'E_2(M')^{-1}=E_2,\]
    \[M'E_3(M')^{-1}=E_4^{\sigma},\] and \[ M'E_4(M')^{-1}=E_3.\]
    So $M'\G_{360}(M')^{-1}=\G_{360}^{\sigma}$.
    Then since
    \[\G_{360}^{\sigma}=M\G_{360}M^{-1}=M'\G_{360}(M')^{-1},\]
    we must have $(M')^{-1}M\in N(\G_{360})=\G_{360}$.  It follows that $M=M'A$ for some $A\in \G_{360}$.

    \item [(r)]  Suppose that $\G=\mf{P}_{A}\rtimes\mf{I}$ is a
    group given by Lemma~\ref{gps}(r), where $\mf{P}_A$ is an elementary abelian $p$-group
    consisting of elements of the form~($A_i$) with $i\in\{1,2\}$ and where $\mf{I}$ is an intransitive group whose
    image $\ol{\mf{I}}$ is one of the groups listed in Lemma~\ref{subgroups} not equal
    to a group given in Lemma~\ref{subgroups}(a) or~(f).
     Let $\psi$ be the inverse transpose isomorphism of $\PGL_3(F)$
    defined in proof of Lemma~\ref{gps}(r).  There exists an
    intransitive element $B\in\PGL_3(F)$ so that
    $B\psi(\G)B^{-1}:=\mf{P}_{A}'\rtimes\mf{I}$, is a group given by Lemma~\ref{gps}(r)
    where $\mf{P}_{A}'$ is an elementary abelian $p$-group
    consisting of elements of the form~($A_j$) with
    $j\in \{1,2\}-\{i\}$.  Let $\phi$ be the automorphism of $\PGL_3(F)$ defined by $D\mapsto B\psi(D)B^{-1}$.
    Note that the image of an intransitive element under $\phi$ is
    intransitive and that $\phi(\mf{P}_{A})=\mf{P}_{A}'$ and that
     $\phi(C)^{\sigma}=\phi(C^{\sigma})$ for all $C\in\PGL_3(F)$.
     Since we will show that $M\G M^{-1}=\G^{\sigma}$ implies that $M=M'A$, with $A\in\G$ and where $M'$ is an
element
    of an intransitive group $\mf{I}'$ containing $\mf{I}$ with $\ol{\mf{I}'}$ equal to the
    normalizer of $\ol {\mf{I}}$, we
    may assume without loss of generality that $i=2$.

    For convenience, we will
    write an element $g\in\G$ as $(v, V)$ where $v$ is a
    $2\times 1$ matrix and where
    $V$ is in $\GL_2(F)$.  That is, if
    \[g:=\begin{bmatrix}
    a&b&\alpha\\
    c&d&\beta\\
    0&0&1\\
    \end{bmatrix}\in\G,\] then we will write $g=(v, V)$ where
     $v=(\alpha, \beta)$ and
    \[V:=\begin{pmatrix}
    a&b\\
    c&d\\
    \end{pmatrix}\in\GL_2(F).\]
Let $(v,V)$, $(v',V')$ be in $\G$.  Then \[(v,V)(v',V')=(v+Vv',
VV').\]

Now suppose that $M\in\PGL_3(F)$ gives an isomorphism
$X\to\el{\sigma}{X}$.  Let $[X_0\colon X_1\colon X_2]$ be coordinate
functions on $\mathbb{P}^2$.  The unique line fixed by $\G$ is
$\{X_2=0\}$.  So the unique line fixed by $\G^{\sigma}$ is
$\{X_2=0\}$.  The unique line fixed by $M\G M^{-1}$ must be the
image of the line $\{X_2=0\}$ under $M$.  Since $\G^{\sigma}=M \G
M^{-1}$, $M$ must fix the line $\{X_2=0\}$.  So we can write
$M:=(u,U)$.  We have
$\G^{\sigma}=\mf{P}_{A}^{\sigma}\rtimes\mf{I}^{\sigma}$.  Let
$g:=(v,V)$ be in $\G$. Then $MgM^{-1}$
\[=(u,U)(v,V)(-U^{-1}u, U^{-1})=(u+Uv-UVU^{-1}u,UVU^{-1})\in\G^{\sigma}.\]  We see
that \[\mf{P}_{A}^{\sigma}=\{(Uv,Id)\colon (v,Id)\in\mf{P}_{A}\}\]
and \[\mf{I}^{\sigma}=\{(0, UVU^{-1})\colon (0,V)\in\mf{I}\}.\]
   Then for all $(0,V)\in\G$ we
must have
\[M(0,V)M^{-1}(0,(UVU^{-1})^{-1})=(u-UVU^{-1}u,Id)\in\G^{\sigma}.\] By Lemma~\ref{diagelt},
$\G^{\sigma}$ has an element of the form $(0,W)$, where
\[W:=\begin{pmatrix}
\zeta_1 &0\\
0&\zeta_2\\
\end{pmatrix}\] and where both $\zeta_1$ and $\zeta_2$ roots of unity are not equal to
 $1$.  Then there exists $n>1$ with $p\nmid n$ such that $\zeta_1$ and $\zeta_2$
 are zeros of $\sum_{i=0}^{n-1} x^i$.  Choose
$l\in\Z_{>0}$ so that $ln\equiv 1\pmod p$.  Since $UV'U^{-1}=W$ for
some $(0,V')\in\G$, the element
\[\left(\prod_{j=0}^{n-1} (u-W^ju,Id)\right)^l=(u,Id)\in\G^{\sigma}.\]
So $(U^{-1}u, Id)$ is in $\G$.  Since
$(0,U)\mf{I}(0,U)^{-1}=\mf{I}^{\sigma}$ by part~(a) of this lemma we
must have $\ol{(0,U)}\in N(\ol{\mf{I}})$, so $(0,U)$ in an element
of an intransitive group $\mf{I}'$ where $\ol{\mf{I}'}$ is the
normalizer of $\ol{\mf{I}}$ in $\PGL_2(F)$.
 Then $M=(u,U)=(0,U)(U^{-1}u,Id)$ is of the desired form.
 \item[(s)]  In this case it may not be true that $\G^{\sigma}=\G$
 so $M$ may not be in the normalizer of $\G$.

\end{enumerate}

\end{proof}

%% file: ch4.tex
Let $K$ be a perfect field, let $F$ be an algebraic closure of $K$,
and let $\Gamma=\Gal(F/K)$. Let $X$ be a hyperelliptic curve over
$F$ and let $B$ be the canonical $K$-model of $X/\Aut(X)$ given in
Theorem~\ref{bk}.  In the proof of Theorem~\ref{bk}, D\`{e}bes and
Emsalem show  the canonical model exists by using the following
argument.  For all $\sigma\in\Gamma$ there exists an isomorphism
$\varphi_{\sigma}\colon X\rightarrow \el{\sigma}{X}$ defined over
$F$.  Each induces an isomorphism $\tilde{\varphi_{\sigma}}\colon
X/\Aut(X)\rightarrow \el{\sigma}{X}/\Aut(\el{\sigma}{X})$ that makes
the following diagram commute:
\[
\begin{CD}
{X} @> {\varphi_{\sigma}} >> {\el{\sigma}{X}} \\
@ V{\rho} VV @ VV {\rho^{\sigma}} V \\
{X/\Aut(X)} @>> {\tilde{\varphi_{\sigma}}} >
{\el{\sigma}{X}/\Aut(\el{\sigma}{X})} \\
\end{CD}
\]
Composing $\tilde{\varphi_{\sigma}}$ with the canonical isomorphism
\[i_{\sigma}\colon \el{\sigma}{X}/\Aut(\el{\sigma}{X})\rightarrow
\el{\sigma}{(X/\Aut(X))}\] we obtain an isomorphism
\[\overline{\varphi_{\sigma}}\colon  X/\Aut(X)\rightarrow
\el{\sigma}{(X/\Aut(X))}.\] The family
$\{\overline{\varphi_{\tau}}\}_{\tau\in\Gamma}$ satisfy Weil's
cocycle condition
$\overline{\varphi_{\tau}}^{\sigma}\,\overline{\varphi_{\sigma}}=\overline{\varphi_{\sigma\tau}}$
given in Theorem~\ref{weil}. This shows that $B$ exists. 

Let $F(B_F)$ be the function field of $B_F$. Since
$B_F\cong\mathbb{P}^1$, $F(B_F)=F(t)$ for some element $t$. We use
$t$ as a coordinate on $B_F$.  Suppose $\sigma\in\Gamma$ and suppose
that $\overline{\varphi_{\sigma}}$ is given by
\[t\mapsto \frac{at+b}{ct+d}.\] Define $\sigma^*\in\Aut(F(t)/K )$ by
\[\sigma^*(t)=\frac{at+b}{ct+d},\; {\sigma^*(\alpha)}=\sigma(\alpha),\; \alpha\in F.\]  One can verify that
$(\sigma\tau)^*(w)=\sigma^*(\tau^*(w))$ for all $w\in F(t)$.  So we
get a homomorphism $\Gamma\to\Aut(F(B_F)/K)$,
$\sigma\mapsto\sigma^*$. The curve $B$ is the variety over $K$
corresponding to the fixed field of
$\Gamma^*=\{\sigma^*\}_{\sigma\in\Gamma}$.
 The following lemma and corollary will be of use.
\begin{lem}\label{odddiv}Let $C$ be
a curve of genus $0$ over $K$ and suppose that $C$ has a divisor $D$
rational over $K$ of odd degree. Then $C(K)\ne\varnothing$.
\end{lem}

\begin{proof} Let $\omega$ be a canonical divisor on $C$.
Since $\Deg(\omega)=-2$, we can take a linear combination of $D$ and
$\omega$ to obtain a divisor $D'$ of degree $1$. Since
$\Deg(\omega-D')<0$, by the Riemann-Roch theorem $l(D')>0$.  So
there exists an effective divisor $D''$ linearly equivalent to $D'$
rational over $K$.  Since $D''$ is effective and of degree 1 it
consists of a point in $C(K)$.
\end{proof}

\begin{cor}\label{odd}Let $L/K$ be a separable field extension of odd degree.  Let $C$ be
a curve of genus $0$ defined over $K$ and suppose that $C(L)\ne\varnothing$. Then $C(K)\ne\varnothing$.
\end{cor}

\begin{proof}Let $P\in C(L)$ and let $n=[L:K]$.  Let $\tau_1, \ldots, \tau_n$ be the distinct embeddings of $L$ into an algebraic closure of $L$.  Then
$D=\Sigma\tau_i(P)$ is a divisor of degree $n$ defined over $K$.  By Lemma~\ref{odddiv}, $C(K)\ne\varnothing$.

\end{proof}
\section{The main result for hyperelliptic curves}
\begin{thm}\label{mainhyp} Let $K$ be a perfect field of characteristic not equal
to $2$ and let $F$ be an algebraic closure of $K$.
  Let
$X$ be a hyperelliptic curve over $F$ and let $\G=\Aut(X)/\langle
\iota\rangle$ where $\iota$ is the hyperelliptic involution of $X$.
Suppose that $\G$ is not cyclic or that $\G$ is cyclic of order
divisible by the characteristic of $F$. Then $X$ can be defined over
 its field of moduli relative to the extension $F/K$.
\end{thm}

\begin{proof}
Let $\Gamma=\Gal(F/K)$.  By Proposition~\ref{closed} we may assume
that $K$ is the field of moduli of $X$. By Proposition~\ref{isom} we
may assume that $\G$ is given by one of the groups in
Lemma~\ref{subgroups}. Fix an equation $y^2=f(x)$ for $X$
  where $f\in F[x]$ and $\mbox{disc}(f)\ne 0$.  So the function field $F(X)$ equals
  $F(x,y)$.
There are eight cases.
\begin{enumerate}

\item[(b1)]  $\G\cong D_4$.  The element $t:=x^2+x^{-2}$ is fixed by $\G_{D_4}$ and is a rational function of
degree $4$ in $x$.  So the function field of $X/\Aut(X)$ equals
$F(t)$.  We use $t$ as a coordinate on
 $X/\Aut(X)$.  The map $\rho\colon X\rightarrow X/\Aut(X)$ is given by $(x,y)\mapsto (x^2+x^{-2})$.  Let
$\sigma\in\Gamma$.  By Lemmas~\ref{isomnormalizer}
and~\ref{normalizersub}, $\varphi_{\sigma}\colon X\to
\el{\sigma}{X}$ is given by $(M,e)$ where $\overline{M}\in
\G_{S_4}$.  A computation shows that $\sigma^*(t)$ is one of the
following:
\begin{enumerate}
  \item [i.] $t$
  \item [ii.] $-t$
  \item [iii.] $\frac{2t+12}{t-2} $
  \item [iv.] $\frac{2t-12}{-t-2}$
  \item [v.] $\frac{2t-12}{t+2}$
  \item [vi.] $\frac{2t+12}{-t+2}$.
\end{enumerate}

Since $\overline{\varphi}_{\tau} \colon  X/\Aut(X) \to \el{\tau}{(X/\Aut(X))}$ is defined over $K$ for all
$\tau\in\Gamma$, we have
$\overline{\varphi_{\tau}}\,\overline{\varphi_{\sigma}}=\overline{\varphi_{\sigma\tau}}$ for all
$\tau\in\Gamma$.  The fractional linear transformations i through vi form a group under composition isomorphic
to $S_3$. The map $\tau \mapsto \tau^*|_{K(t)}$ defines a homomorphism from $\Gamma$ to this group. The kernel
of this homomorphism is $\Lambda:=\{\tau\in\Gamma\mid \tau^*(t)=t\}$.  So $|\Gamma/\Lambda|= 1, 2, 3, $ or $6$.
\begin{enumerate}
  \item [Case 1:] $|\Gamma/\Lambda|=1$.  In this case the fixed field of
$\Gamma^*$ is $K(t)$ and
  $B=\mathbb{P}^1_K$.
  \item [Case 2:] $|\Gamma/\Lambda|=2$.  Let $\sigma$ be a representative of the
  nontrivial coset.  There are three cases.
\begin{enumerate}
  \item  $\sigma^*(t)=-t$.  Then $t=0$ corresponds to a point $P\in B(K)$.
  \item $\sigma^*(t)=\frac{2t+12}{t-2}$.  Then $t=6$ corresponds to a point $P\in B(K)$.
  \item $\sigma^*(t)=\frac{2t-12}{-t-2}$.  Then $t=-6$ corresponds to a point ${P\in B(K)}$.
\end{enumerate}
 \item [Case 3:] $|\Gamma/\Lambda|=3$.  Since the fixed field of
$\Lambda^*$ is
  $F^{\Lambda}(t)$, $B$ has a $F^\Lambda$-rational point.
By Corollary~\ref{odd}, since $[F^{\Lambda}:K]$ is odd, $B$ has a
$K$-rational point.
  \item [Case 4:] $|\Gamma/\Lambda|=6$.  Let $\Pi$ be a
  subgroup of $\Gamma$ containing $\Lambda$ such that
$\Pi/\Lambda$ is a subgroup of $\Gamma/\Lambda$ of order
  $2$. By
Case 2, $B$ has a $F^{\Pi}$ rational point. Since $[F^{\Pi}:K]=3$ is
odd, by Corollary~\ref{odd},
 $B$ has a $K$-rational point.

\end{enumerate}
\item[(b2)]
$\G\cong D_{2n}$, $n>2$.  The function field of $X/\Aut(X)$ equals
the subfield of $F(X)$ fixed by $\G_{D_{2n}}$ acting by fractional
linear transformations. Then $t:=x^n+x^{-n}$ is fixed by
$\G_{D_{2n}}$ and is a rational function of degree $2n$ in $x$, so
the function field of $X/\Aut(X)$ equals $F(t)$. Therefore we use
$t$ as coordinate on $X/\Aut(X)$.  The map $\rho\colon X\rightarrow
X/\Aut(X)$ is given by $(x,y)\mapsto (x^n+x^{-n})$. Let
$\sigma\in\Gamma$.  By Lemmas~\ref{isomnormalizer}
and~\ref{normalizersub}, $\varphi_{\sigma}\colon X\to
\el{\sigma}{X}$ is given by $(M,e)$ where $\overline{M}\in D_{4n}$.
Then the map $\rho^{\sigma}\varphi_{\sigma}\colon X\rightarrow
\el{\sigma}{X}/\Aut(\el{\sigma}{X})$ is given by $(x,y)\mapsto \pm
(x^n+x^{-n})$.  So $\sigma^*(t)=\pm t$.  The curve $B$ corresponds
to the fixed field of $F(t)$ under $\Gamma^*$. Then $t=0$
corresponds to a point $P\in B(K)$.

\item[(c)]  $\G\cong A_4$.  The element $t':=x^2+x^{-2}$ is fixed by the normal subgroup $\G_{D_4}$.  From (c), we
see that the element
\[t:=\frac{1}{4}t'\left(\frac{2t'-12}{t'+2}\right)\left(\frac{2t'+12}{-t'+2}\right)
=\frac{x^{12} - 33x^8 - 33x^4 + 1}{-x^{10} + 2x^6 - x^2}\] is fixed
by $\G_{A_4}$ and is a rational function of degree $12$ in $x$.  So
the function field of $X/\Aut(X)$ equals $F(t)$.  We use $t$ as
coordinate on $X/\Aut(X)$.  The map $\rho\colon X\rightarrow
X/\Aut(X)$ is given by
\[(x,y)\mapsto (x^{12} - 33x^8 - 33x^4 + 1)/(-x^{10} + 2x^6 - x^2).\]Let $\sigma\in\Gamma$. By Lemmas~\ref{isomnormalizer} and~\ref{normalizersub},
$\varphi_{\sigma}\colon X\to \el{\sigma}{X}$ is given by $(M,e)$
where $\overline{M}\in \G_{S_4}$.    A computation shows that
$\sigma^*(t)=\pm t$.  Then $t=0$ corresponds to a point $P\in B(K)$.

\item[(d)]  $\G\cong S_4$.  By Lemma~\ref{normalizersub}, $N(\G)=\G$.  So by Corollary~\ref{def}, $X$ can be defined
over $K$.
 \item[(e)] $\G\cong A_5$. By Lemma~\ref{normalizersub}, $N(\G)=\G$.  So by Corollary~\ref{def}, $X$ can be
defined over $K$.

\item[(f)]  $\G=\G_{\beta, A}$. Let $d$ be the order of $\beta$ and let $t=g(x):=\prod_{\alpha\in A}(x-\alpha)^d$.
Then $t$ is a rational function of degree $|\G|$ fixed by
$\G_{\beta, A}$ acting by fractional linear transformations.  So the
function field of $X/\Aut(X)$ equals $F(t)$.  We use $t$ as a
coordinate function of $X/\Aut(X)$.  Let $\sigma\in\Gamma$. By
Lemma~\ref{isomnormalizer}, $\varphi_\sigma\colon X\to
\el{\sigma}{X}$ is given by $(M,e)$ where $M$ is an upper diagonal
matrix.  So $\sigma^*(t)=g^{\sigma}(ax+b)$ for some $a\ne 0$ and
$b$.  Let $P$ be the point of $X/\Aut(X)$ corresponding to
$x=\infty$. Then since $g^{\sigma}(a\infty+b)=g(\infty)$, $P$
corresponds to a point in $B(K)$.

\item[(g)] $\G=\PSL_2(\mathbb{F}_{q})$.    It can be deduced from Theorem 6.21 on page 409
of~\cite{Suzuki:1} that $\PSL_2(\mathbb{F}_q)$ is generated by the
image in $\PGL_2(F)$ of the following matrices
\[\left\{\begin{pmatrix}
           0&-1\\
           1&0\\
           \end{pmatrix},\begin{pmatrix}
           1&a\\
           0&1\\
           \end{pmatrix}\colon a\in\F_{p^r}\right\}.\]
Let \[g(x)=\frac{((x^q
-x)^{q-1}+1)^{\frac{q+1}{2}}}{(x^q-x)^{\frac{q^2-q}{2}} }.\] One can
verify that $g(-1/x)=g(x)$ and $g(x+a)=g(x)$ for all $a\in\F_{p^r}$.
Since $g$ is a rational function of $x$ of degree
$\frac{q^3-q}{2}=|\PSL_2(\F_q)|$, the function field of $X/\Aut(X)$
is $F(t)$ where $t=g(x)$.  We use $t$ as a coordinate function on
$X/\Aut(X)$.  The map $\rho\colon X\rightarrow X/\Aut(X)$ is given
by
\[(x,y)\mapsto \frac{((x^q
-x)^{q-1}+1)^{\frac{q+1}{2}}}{(x^q-x)^{\frac{q^2-q}{2}} }.\] Let
$\sigma\in\Gamma$. By Lemmas~\ref{isomnormalizer}
and~\ref{normalizersub}, $\varphi_{\sigma}\colon X\to
\el{\sigma}{X}$ is given by $(M,e)$ where $\overline{M}\in
\PGL_2(\mathbb{F}_q)$.   A computation shows that $\sigma^*(t)=\pm
t$.  Then $t=0$ corresponds to a point $P\in B(K)$.

 \item[(h)] $\G= \PGL_2(\mathbb{F}_{q})$.  By Lemma~\ref{normalizersub}, $N(\G)=\G$.  So by
Corollary~\ref{def}, $X$ can be defined over $K$.

\end{enumerate}
\end{proof}
\begin{thm}\label{mainhyp2} Let $K$ be a field of characteristic not equal
to $2$, let $X$ be a hyperelliptic curve over $K$ and let
$\G=\Aut(X)/\langle \iota\rangle$ where $\iota$ is the hyperelliptic
involution of $X$. Suppose that $\G$ is not cyclic or that $\G$ is
cyclic of order divisible by the characteristic of $F$. Then $X$ is
definable over its field of moduli.
\end{thm}
\begin{proof}This follows from Theorem~\ref{mainhyp} and
Theorem~\ref{relmod}.
\end{proof}
\section{A note about defining equations}
Let $K$ be a perfect field of characteristic not equal to $2$ and
let $F$ be an algebraic closure of $K$.  Let $X$ be a hyperelliptic
$X$ curve over $F$ with field of moduli $K_m$ and let $\iota$ be the
hyperelliptic involution of $X$. By Theorem~\ref{mainhyp}, $X$ is
definable over $K_m$.  In this section we show that $X$ is given by
an equation of the form $y^2=f(x)$, with $f(x)\in K_m[x]$.

\begin{lem}\label{eveng}Let $X$ be a hyperelliptic curve of even genus $g$ defined over a field $K$
 of characteristic not equal to $2$.
Then $X$ is $K$-isomorphic to a hyperelliptic curve given by an equation of the form $y^2=f(x)$ with $f\in
K[x]$.
\end{lem}
\begin{proof}
The canonical morphism factors as
\[X\xrightarrow{\rho} C\xrightarrow{i}\mathbb{P}^{g-1}_K\] where $\Deg(\rho)=2$ and $C$ is a $K$-curve of genus $0$.
On $\mathbb{P}^{g-1}_K$ we have the invertible sheaf $\mathcal{O}(1)$ and $\mathcal{L}:=i^*\mathcal{O}(1)$ is an
invertible sheaf on $C$.  Since $\rho^*(\mathcal{L})$ is a canonical divisor $\omega$ on $X$, we have $2\Deg
\mathcal{L}=\Deg\omega=2g-2$. So $\Deg \mathcal{L}=g-1$ is odd. By Lemma~\ref{odd}, $C$ has a $K$-rational
point. So $C\cong_K\mathbb{P}^1_K$. So the function field of $K(X)$ of $X$ is a quadratic extension of $K(x)$.
By Kummer theory $K(X)=K(x,y)$ where $y^2=f(x)$ for some $f\in K[x]$.
\end{proof}

\begin{prop}\label{canform}Let $X$ be as in Theorem~\ref{mainhyp}. Then $X$ has a $K_m$-model given by an equation of the
form $z^2=h(x)$ with $h\in K_m[x]$, where $K_m$ is the field of
moduli of $X$ relative to the extension $F/K$.
\end{prop}
\begin{proof}
Throughout this proof we use the same notation and make the same
assumptions as in the proof of Theorem~\ref{mainhyp}. So, for
example, we assume that $K_m=K$.

 In cases (b1), (b2), (c), and (g) we construct a family of
automorphisms $\{\psi_{\sigma}\}_{\sigma\in\Gamma}$ of
$\mathbb{P}^1$ that satisfy Weil's cocycle condition of
Theorem~\ref{weil} and such that each automorphism $\psi_{\sigma}$
lifts to an isomorphism $X\to\el{\sigma}{X}$. Let $C$ be the
$K$-model of $\mathbb{P}^1$ corresponding to
$\{\psi_{\sigma}\}_{\sigma\in\Gamma}$.  We then show that $C$ has a
$K$-rational point, so $C$ is isomorphic over $K$ to
$\mathbb{P}^1_K$.  Then by Theorem~\ref{weil} there exists an
automorphism $\psi$ of $\mathbb{P}^1$ such that for all
$\sigma\in\Gamma$ we have $(\psi^{-1})^{\sigma}\psi=\psi_{\sigma}$.
By Lemma~\ref{lift}, $X$ has a $K$-model $Y$ given by $z^2=h(x)$.

If $\G\ne \G_{\beta, A}$ then for each $\sigma\in\Gamma$ we have
$\sigma^*(t)\in K(t)$. So the map $\sigma\mapsto\sigma^*|_{K(t)}$
defines a homomorphism from $\Gamma$ to a group of fractional linear
transformations of $K(t)$ whose kernel is
$\Lambda:=\{\sigma\in\Gamma\mid\sigma^*(t)=t\}$.  If $\G\ne
\G_{\beta, A}$ let $L=F^{\Lambda}$, and if $\G=\G_{\beta, A}$ let
$L=F$. By Lemma~\ref{ingsig}, we may assume that $f\in L[x]$.

There are eight cases.
\begin{enumerate}
\item[(b1)] $\G\cong D_{2n}$, $n>2$.  In this case $[L:K]$ is equal to $1$ or $2$.  If $[L:K]=1$ there is nothing
to prove so assume $[L:K]=2$.  So $L=K(\sqrt c)$ for some $c\in K$.

 Let $\tau\in\Gamma$. By
Lemma~\ref{isomnormalizer}, if $\varphi_{\tau}\colon
X\to\el{\tau}{X}$ is an isomorphism given by $(M_{\tau},e_{\tau})$,
then the image of $M_{\tau}$ in $\PGL_2(F)$ is an element of
$N(\G_{D_{2n}})=\G_{D_{4n}}$.  If $\tau|_L=Id$, then $\tau^*(t)=t$
so $\ol M_{\tau}\in \G_{D_{2n}}$. If $\tau|_L\ne Id$, then
$\tau^*(t)\ne t$ so $\ol M_{\tau}\in \G_{D_{4n}}-\G_{D_{2n}}$.
Composing $\varphi_{\tau}$ with any automorphism of $X$ gives
another isomorphism $X\to\el{\tau}{X}$, so any element in the same
coset as $\ol M_{\tau}$ in $\G_{D_{4n}}/\G_{D_{2n}}$ will lift to an
isomorphism $X\to\el{\tau}{X}$.

  Choose $\gamma\in F$
satisfying $\gamma^n=\sqrt c$. For $\sigma\in\Gamma$ define
$\psi_{\sigma}$ by $(x)\mapsto (\frac{\sigma(\gamma)}{\gamma}x)$. If
$\sigma|_L=Id$, then $\frac{\sigma(\gamma)}{\gamma}$ is an $n^{th}$
root of unity.  If $\sigma|_L\ne Id$, then
$\frac{\sigma(\gamma)}{\gamma}$ is a $2n^{th}$ root of unity that is
not an $n^{th}$ root of unity.  So $\psi_{\sigma}$ lifts to a map
$X\to\el{\sigma}{X}$ for all $\sigma\in\Gamma$. The family
$\{\psi_{\sigma}\}_{\sigma\in\Gamma}$ satisfies Weil's cocycle
condition and the point $x=0$ corresponds to a point $P\in C$.

\item[(b2)]  $\G\cong D_4$.  Let $\sigma\in\Gamma$.  In this case $[L:K]$ is equal to $1, 2, 3, $ or $6$. We first
show that if $2$ divides $[L:K]$, then we may assume that there
exists an element $\sigma\in\Gamma$ such that $\sigma^*(t)=-t$.  So
suppose that $2$ divides $[L:K]$ and let $\sigma\in\Gamma$ be such
that $\sigma|_L$ has order $2$. As shown in the proof of
Theorem~\ref{mainhyp}(c), $\sigma^*(t)$ is equal to one of the
following:
\begin{enumerate}
  \item[i.]  $-t$
  \item[ii.] $\frac{2t+12}{t-2}$
  \item[iii.] $\frac{2t-12}{-t-2}$.
\end{enumerate}
If $\sigma^*(t)=\frac{2t+12}{t-2}$, then there is an isomorphism $(M_{\sigma},e_{\sigma}): X\to\el{\sigma}{X}$
with $\ol{M_{\sigma}}$ equal to the image of
\[\begin{pmatrix}
1 & -1\\
1 & 1\\
\end{pmatrix}\]
in $\PGL_2(F)$.  Let $T$ be the image of
\[ \begin{pmatrix}
i & 1\\
1 & i\\
\end{pmatrix}\]
in $\PGL_2(F)$.
 Since $T\in N(\G_{D_4})$ and since $T^{\sigma}\ol{M_{\sigma}}T^{-1}$ is given by the image in $\PGL_2(F)$ of
 either
 \[ \begin{pmatrix}
-i & 0\\
0 & 1\\
\end{pmatrix}\]
or
\[ \begin{pmatrix}
0 & i\\
1 & 0\\
\end{pmatrix}\]
 replacing $f(x)$ with $f(\frac{-ix+1}{x-i})$ we may assume that $\sigma^*(t)=-t$.  Similarly, if
 $\sigma^*(t)=\frac{2t-12}{-t-2}$, replacing $f(x)$ with $f(\frac{x+1}{-x+1})$, we may assume that
 $\sigma^*(t)=-t$.
\begin{enumerate}
  \item [Case 1:] $[L:K]=1$.  Clear.
  \item [Case 2:] $[L:K]=2$.  Then $L=K(\sqrt c)$ for some $c\in K$.  Choose $\gamma\in F$ so that $\gamma^2=\sqrt c$.
  For $\sigma\in\Gamma$ define $\psi_{\sigma}$ by
$(x)\mapsto (\frac{\sigma(\gamma)}{\gamma}x)$.  Then as in (b), each
automorphism $\psi_{\sigma}$ lifts to an isomorphism
$X\to\el{\sigma}{X}$ and the family
$\{\psi_{\sigma}\}_{\sigma\in\Gamma}$ satisfies Weil's cocycle
condition. The point $x=0$ corresponds to a point $P\in C$.

\item [Case 3:] $3|[L:K]$.  Let $P_1,\ldots, P_k$ be the hyperelliptic branch points of $X\to\mathbb{P}^1$. The
group $\G_{D_4}$ acting by linear transformation permutes these
points.  By Proposition~2.1 and Lemma~2.2 of~\cite{Brandt:1}, $f(x)$
must be a scalar multiple of a polynomial of one of the forms:
\[g(x), xg(x), (x^2-1)g(x), (x^2+1)g(x), (x^4-1)g(x),\]
\[ x(x^2-1)g(x), x(x^2+1)g(x), \mbox{ or}, x(x^4-1)g(x),\]
for some $g(x)$ where
\[g(x):=\prod_{i=1}^l(x^4+\lambda_ix^2+1),\;\;\lambda_i\in F-\{\pm 2\},\lambda_i\ne\lambda_j\mbox{ if }i\ne j.\]

 Let $\sigma\in\Gamma$,
suppose that $\sigma|_L$ has order $3$, and let
$\varphi_{\sigma}\colon X\to\el{\sigma}{X}$ be an isomorphism given
by $(M_{\sigma},e_{\sigma})$.  By Lemmas~\ref{isomnormalizer}
and~\ref{normalizersub}, $\ol M_{\sigma}\in N(\G_{D_4})=\G_{S_4}$.
Since $(\sigma^*)^3(t)=t$,  $\ol M_{\sigma}^3\in \G_{D_4}$.  Note
that $S_4/D_4\cong S_3$. Replacing $\sigma$ with $\sigma^2$ if
necessary and composing $\varphi_{\sigma}$ with an automorphism of
$X$ if necessary, by Lemma~\ref{normalizersub}, we may assume that
$\ol M_{\sigma}$ is given by the image in $\PGL_2(F)$ of the matrix
\[\begin{pmatrix}
 i & i\\
 1 & -1\\
 \end{pmatrix}.\]

Then, writing $\sigma(\infty)=\infty$, the hyperelliptic branch points $\sigma(P_1)\ldots\sigma(P_K)$ of
$\el{\sigma}{X}\to\mathbb{P}^1$ are
\[\frac{iP_1+i}{P_1-1},\ldots, \frac{iP_k+i}{P_k-1}.\]
      So if $0$ and $\infty$ are branch points of $X\to\mathbb{P}^1$, then $i$ and $-i$ must be branch points
      of $\el{\sigma}{X}\to\mathbb{P}^1$.  In this case, since $\{\sigma^{-1}(i),\sigma^{-1}(-i)\}=\{i,-i\}$,
      $i$ and $-i$ must also be branch points of $X\to\mathbb{P}^1$.  Then $1$ and $-1$ must be branch points of
      both $X$ and $\el{\sigma}{X}$.  It follows that if $x|f(x)$, then $x(x^4-1)|f(x)$.  Similarly, one can
      show that if $x^2\pm 1|f(x)$ then $x(x^4-1)|f(x)$.

So without any loss of generality we may assume that
\[f(x)=g(x)\mbox{ or } x(x^4-1)g(x).\]
     If $f(x)=x(x^4-1)g(x)$, since $\Deg(f)=2g+1$ where $g$ is the genus of $X$, $g$ is even. By Theorem~\ref{mainhyp},
$X$ is definable over $K$.
 By Lemma~\ref{eveng}, $X$ has a $K$-model given by an equation of the form $z^2=h(x)$ where $h\in K[x]$.

Assume that $f(x)=g(x)$ and let $X'$ be the hyperelliptic curve over
$F$ given by $z^2=x(x^4-1)g(x)$.  So the hyperelliptic branch points
of $X'\to\mathbb{P}^1$ are $P_1\ldots P_k,0,\infty,\pm 1,\pm i$. Let
$\tau\in\Gamma$. Suppose that $(M_{\tau},e_{\tau})$ gives an
isomorphism $X\to\el{\tau}{X}$. Suppose that
\[ M_{\tau}:=\begin{pmatrix}
a&b\\
c&d\\
\end{pmatrix}.\]
By Lemmas~\ref{isomnormalizer} and~\ref{normalizersub}, $\ol
M_{\tau}\in N(D_4)=\G_{S_4}$.  The hyperelliptic branch points of
$\el{\sigma}{X}\to\mathbb{P}^1$ are
\[\frac{aP_1+b}{cP_1+d},\ldots, \frac{aP_k+b}{cP_k+d}.\] Since $\infty$ and the roots of $x(x^4-1)$ are permuted by any element of
$\G_{S_4}$, it is easily verified that there exist $\lambda\in
F^{\times}$ such that $(M_{\tau},\lambda e_{\tau})$ gives an
isomorphism $X'\to\el{\tau}{X'}$.  Similarly, given an isomorphism
$X'\to\el{\tau}{X'}$ given by $(M'_{\tau},e_{\tau}')$, there exists
$\lambda'\in F^{\times}$ such that $(M'_{\tau}\,\lambda' e_{\tau}')$
gives an isomorphism $X\to\el{\tau}{X}$.

        By the above
there exists a $K$-model $Y$ of $X'$ given by $w^2=r(x)$ with $r\in
K[x]$. Let $\varphi'\colon X'\to Y\times_K  F$ be an isomorphism
given by $(M,e)$. Suppose that
\[M=\begin{pmatrix}
\alpha & \beta\\
\gamma & \delta\\
\end{pmatrix}.\] Define $\psi\colon \mathbb{P}^1\to\mathbb{P}^1$ by $x\mapsto\frac{\alpha x+\beta}{\gamma x+\delta}$.  Then for all $\sigma$ in
$\Gamma$, $(\psi^{-1})^{\sigma}\psi$ lifts to an isomorphism $X\to\el{\sigma}{X}$.  By Lemma~\ref{lift}, $X$ has
a $K$-model given by an equation of the form $z^2=h(x)$ where $h\in K[x]$.

\end{enumerate}

\item[(c)]  $\G\cong A_4$.  The proof is identical to the proof of Case 1 and Case 2 of~(b2).

\item[(d)]  $\G\cong S_4$.  Clear since $L=K$.

 \item[(e)] $\G\cong A_5$.  Clear since $L=K$.

\item[(f)]  $\G=\G_{\beta, A}$.  By Theorem~\ref{mainhyp}, $X$ is definable over $K$.  Let
$\{\varphi_{\sigma}\}_{\sigma\in\Gamma}$ be a family of isomorphisms
$\varphi_{\sigma}\colon X\to\el{\sigma}{X}$ given by
$\{(M_{\sigma},e_{\sigma})\}_{\sigma\in\Gamma}$ satisfying Weil's
cocycle condition.  Let $Y$ be the corresponding $K$-model of $X$.
By Lemma~\ref{isomnormalizer} the induced automorphisms of
$\mathbb{P}^1$, given by the images in $\PGL_2(F)$ of
$\{M_{\sigma}\}_{\sigma\in\Gamma}$, fix $P=\infty$.  It follows that
the function field $K(Y)$ is a quadratic extension of $K(x)$.  By
Kummer theory, $Y$ is given by an equation of the form $z^2=h(x)$
with $h\in K[x]$.

\item[(g)] $\G=\PSL_2(\mathbb{F}_{q})$. The proof is identical to the proof of~(b1).

 \item[(h)] $\G= \PGL_2(\mathbb{F}_{q})$. Clear since $L=K$.

\end{enumerate}

\end{proof}

%% file: ch5.tex
By Theorem~22 of~\cite{Cardona:2}, a curve of genus $2$ defined over
an algebraic closure $F$ of a perfect field $K$ of characteristic
$2$
 can be defined over its field of moduli relative to $F/K$.   So
by Theorem~\ref{relmod}, any curve of genus $2$ over a field of
characteristic
  $2$ can be defined over its field of moduli.  Let $X$ be a hyperelliptic curve over a field $K$ and let $\iota$ be
the hyperelliptic involution of $X$.   By
  Theorem~\ref{mainhyp2}, and the results of~\cite{Cardona:2}, if
  $X$ is not definable over its field of moduli then the
  characteristic of $K$ is not $2$ and $\Aut(X)/\la\iota\ra$ is a cyclic group of order not equal to the characteristic
  of $K$.

  The first examples of curves not definable over their fields
of moduli were discovered by Shimura.  These curves are
hyperelliptic $\C$-curves with automorphism groups isomorphic to
$\Z/2\Z$ and are given on page~177 of~\cite{Sh:1}.

The authors of~\cite{Bujalance} have attempted to classify all
hyperelliptic $\C$-curves with fields of moduli $\R$ relative to
$\C/\R$ but not definable over $\R$.  Due to some errors in their
paper, some curves are missing from their list and many curves  on
their list are, in fact, definable over $\R$. We give the full and
correct list in this section.

Suppose $n$, $m$, $r$, $s\in\Z_{>0}$.  Assume that $2nr>5$, and that
$sm$ is even.   Assume also that if $n$ is odd, then $r$ is also
odd.  Let $z^c$ be the complex conjugate of $z$ for any $z\in\C$.
Suppose $a_1,\ldots,a_r, b_1,\ldots, b_s \in \C$.  Consider the
polynomials $f(x), g(x)\in\C[x]$ given by
\begin{align*}
  f(x)&:=\prod_{i=1}^r(x^n-a_i)(x^n+1/a_i^c),\\
 g(x)&:=x\prod_{i=1}^s(x^m-b_i)(x^m+1/b_i^c).\\
 \end{align*}\vspace*{-.76in}
\newline   Assume that $f(x)$ and $g(x)$ both have no
repeated zeros and that both are not polynomials in $\R[x]$. Assume
that the map $P\mapsto 1/P$ does not map the zero set of $f$ into
itself and does not map the set of nonzero zeros of $g$ to itself.
For any
  root of unity $\zeta\ne 1$, assume
that  \[\{a_i, -1/a_i^c\}_{i=1}^r\ne\{\zeta
a_i,-\zeta/a_i^c\}_{i=1}^r\] and that
\[\{b_i,-1/b_i^c\}_{i=1}^s\ne\{\zeta
b_i,-\zeta/b_i^c\}_{i=1}^s.\] Lastly, if $n=3$ assume that the map
\[P\mapsto\frac{-(P-\sqrt{3}-1)}{P(\sqrt{3}+1)+1}\]  does not map the
zero set of $f$ into itself, and if $m=3$ assume that $1+\sqrt{3}$
is not a zero of~$g$.

\begin{rem}\label{bmistake}  Here we discuss the errors in~\cite{Bujalance}.

Due to an error in their paper on page~5, the curves isomorphic to
those given by $y^2=f(x)$ with $n$ odd are missing.  We describe
their error here. Suppose $g\in\Z_{>1}$, $n\in\Z_{>0}$ and that
$\gamma:=(g+1)/n$ is an integer.  Define
\[\Lambda:=\la d,x_1,\ldots,x_{\gamma},y\mid d^2x_1x_2\ldots x_{\gamma} y=1,x_i^2=1,y^n=1\ra.\]
So in the notation of page~5 of~\cite{Bujalance}, $\Lambda$ has
presentation~(i).  In the second to last paragraph of page~5
of~\cite{Bujalance} the authors state that there does not exist a
surjective homomorphism $\theta\colon\Lambda\to\Z/4\Z$. This is
false since if $\gamma$ and $n$ are both odd a surjective
homomorphism is given by $\theta(d)=1$, $\theta(x_i)= 2n$,
$\theta(y)=2n-2$.  With an easy adjustment of their proof we can
show the existence of such curves.

Due to an error on page~10, the curves isomorphic to the curves
given by $y^2=g(x)$ with $m$ odd are missing.  For $m$ odd, the
authors give equations of curves that are actually definable over
$\R$. These are curves given by equation~(11) of~\cite{Bujalance} on
page~10:
\begin{align*}
  w^2&=z\prod_{i=1}^s(z^m-d_i)(z^m-(-1)^{m+1}/d_i^c)\\
 &=z\prod_{i=1}^s(z^m-d_i)(z^m-1/d_i^c)\\
 \end{align*}\vspace*{-.76in}
\newline with certain conditions on $d_1,\ldots, d_s\in\C$ and where $sm$ is even.  Let $Z$ be a
hyperelliptic curve given by such an equation.
 Then the map $\mu$ defined by
\[(z,w)\mapsto (1/z,e_Zw/z^{sm+1}) \] where
\[e_Z^2=\prod_1^s d_i^c/d_i\] gives an isomorphism $Z\to\el{c}{Z}$.  Note that
$|e_Z|=1$,
 so we must have $e_Z^ce_Z=1$.
A computation shows that $\mu^c\mu=Id$.  Thus by Theorem~\ref{weil},
$Z$ is definable over $\R$.  This error can be explained as follows.

We follow the notation of~\cite{Bujalance}.  On page~10, the authors
assume that $\omega\in\C$ is an element such that
$\omega^2=\zeta_{m}$, where $\zeta_m$ is a primitive $m^{th}$ root
of unity.   They assume that their curve $Z$ has an equation of the
form $w^2=f(z)$ and they consider anticonformal automorphisms of the
field $\C(w,z)/(w^2-f(z))$ and the anticonformal automorphisms
induced by these on the field of meromorphic functions $\C(z)$ on
the Riemann sphere. In particular, they state that the anticonformal
map $\tau\colon\C(z)\to\C(z)$ with $\tau(i)=-i$ and
$\tau(z)=-1/(\omega z)$ has no fixed points. Note that
$\tau^2(z)=\zeta_mz$ and that $\tau^2$ is conformal.  By ``fixed
point" of $\tau$ they mean an element of the form
$z':=\frac{az+b}{cz+d}$, with $a,b,c,d\in\C$, such that
the determinant of \[\begin{pmatrix} a&b\\
c&d\\
\end{pmatrix}\] is nonzero, and such that $\tau(z')=z'$.  By Proposition~1 of~\cite{Bujalance},
if an induced anticonformal map of the Riemann sphere of order $2$
with fixed points exists, then the corresponding hyperelliptic curve
can be defined over $\R$.  In writing their defining equation they
make the assumption, although not explicitly stated, that
$\omega^m=-1$.  So since $m$ is odd, the element $-1/\omega$ is an
$m^{th}$ root of unity.  Composing $\tau$ with an appropriate power
of $\tau^2$ we get an anticonformal map $\tau^{2k+1}$ of order $2$
defined by $\tau^{2k+1}(z)=1/z$, $\tau^{2k+1}(i)=-i$. Clearly the
point $\frac{i(z-1)}{z+1}$ is a fixed point.

Our equations $f(x)$ and $g(x)$ are obtained in the same method as
in~\cite{Bujalance}.  For $f(x)$, we follow the argument on page~7
but we use the anticonformal automorphism $\tau\colon\C(z)\to\C(z)$
defined by $\tau(i)=-i$ and $\tau(z)=1/(\zeta_{2n} z)$, where
$\zeta_{2n}$ is a primitive $2n^{th}$ root of unity.  For $g(x)$, we
follow the argument on page~10 but we use the anticonformal
automorphism $\tau\colon\C(z)\to\C(z)$ defined by $\tau(i)=-i$ and
$\tau(z)=1/(\zeta_{2m} z)$, where $\zeta_{2m}$ is a primitive
$2m^{th}$ root of unity.  In both cases it can be verified that
$\la\tau\ra$ contains no anticonformal automorphism of order $2$.

\end{rem}

\begin{lem}\label{hypce}Following the notation above, let $X$ be the
hyperelliptic curve over $\C$ given by $y^2=f(x)$ and let $Y$ be the
hyperelliptic curve over $\C$ given by $y^2=g(x)$.  Let $\iota_X$ be
the hyperelliptic involution of $X$ and let $\iota_Y$ be the
hyperelliptic involution of $Y$.  Then
$\Aut(X)\cong\Z/2\Z\times\Z/n\Z$ and $\Aut(Y)\cong\Z/2m\Z$. We have,
\[\Aut(X)/\la\iota_X\ra\cong\Z/n\Z\] and
\[\Aut(Y)/\la\iota_Y\ra\cong\Z/m\Z.\]
Furthermore, an isomorphism $X\to\el{c}{X}$ is given by
\[(x,y)\mapsto (1/\zeta_{2n}x,e_Xy/x^{nr})\] and an isomorphism $Y\to\el{c}{Y}$ is given by
\[(x,y)\mapsto (1/\zeta_{2m}x,e_Yy/x^{sm+1}), \] where
$\zeta_{2n}$ is a primitive $2n^{th}$ root of unity, $\zeta_{2m}$ is
a primitive $2m^{th}$ root of unity,
\[e_X^2=\prod_1^r-a_i^c/a_i,\] and
\[e_Y^2=(1/\zeta_{2m})\prod_1^s -b_i^c/b_i.\]
\end{lem}
\begin{proof}  It is clear that $y^2=f(x)$ and $y^2=g(x)$ give the equations of hyperelliptic curves.
The claims about the automorphism groups follow from Lemma~3.1 and
Corollary~3.2 of~\cite{Bujalance}.  Two very simple computations
show that the maps above give isomorphisms $X\to\el{c}{X}$ and
$Y\to\el{c}{Y}$.
\end{proof}
\begin{prop}\label{hypnotdef} We follow the notation above.  Let $Z$ be a hyperelliptic curve over $\C$ with field
of moduli $\R$ relative to $\C/\R$.  Then $Z$ is not definable over
$\R$ if and only if $Z$ is isomorphic to a curve given by an
equation of the form $y^2=f(x)$ or $y^2=g(x)$ where $f(x)$ and
$g(x)$ satisfy all of the conditions given above.
\end{prop}
\begin{proof}  This follows easily from the arguments given in~\cite{Bujalance}
after fixing the errors of~\cite{Bujalance} mentioned in
Remark~\ref{bmistake}.  We note that it is also quite easy to show
that Weil's cocycle condition of Theorem~\ref{weil} cannot be
satisfied in either case.
\end{proof}

 Note that from Proposition~\ref{hypnotdef}, we can easily construct
curves over $\Q(i)$ with fields of moduli equal to $\Q$ but not
definable over $\Q$.  We simply need to take $f(x)$ or $g(x)$ in
$\Q(i)[x]$.   Let $X$ be such a curve.  Let $\ol\Q$ be the algebraic
closure of $\Q$ contained in $\C$. It is easy to see that the field
of moduli of $X$ is $\Q$.  The curve $X$ cannot be definable over
$\Q$ because that would imply that the curve $X_{\C}$ is definable
over $\R$.

%% file: ch6.tex
\section{Invariant forms}
 Let $F$ be an algebraically closed field of characteristic $p\ne 2$. Let
 $f,g\in F[X_0,\ldots,X_{n-1}]$ be two forms, i.e. two homogeneous
polynomials.  We define an equivalence relation on the set of forms
by $f\sim g$ if $f=\lambda g$ for some $\lambda\in F^{\times}$.  Let
$[f]$ denote the equivalence class of $f$. Suppose
\[M:=\begin{bmatrix}
a_{11}&\ldots&a_{1n}\\
.&\ldots&.\\
a_{n1}&\ldots&a_{nn}\\
\end{bmatrix}\in\PGL_{n}(F).\]     We
 define a right action of $\PGL_n(F)$ on the set
  of equivalence classes of
forms by
\[([f(X_0,\ldots,X_{n-1})])M=[f(a_{11}X_0+\cdots+a_{1n}X_{n-1},\ldots,a_{n1}X_0+\cdots+a_{nn}X_{n-1})].\]
 So $([f])MM'=(([f])M))M'$ for all $M$, $M'\in\PGL_n(F)$.
Let $\G$ be a subgroup of $\PGL_n(F)$.  We say that a form $f$ is
$\G$-invariant if
 $([f])M= [f]$ for all $M\in\G$. If
 $P=[\alpha_0\colon\cdots\colon\alpha_{n-1}]\in\mathbb{P}^{n-1}(F)$,
let $M(P)$ denote the point
\[[a_{11}\alpha_0+\ldots+a_{1n}\alpha_{n-1}\colon\cdots\colon a_{n1}\alpha_0+\ldots+a_{nn}\alpha_{n-1} ].\]
\section{Invariant binary forms}
Let $F$ be an algebraically closed field of characteristic $p\ne 2$.
Let
 $f\in F[X_0,X_1]$ be a binary form, i.e. a homogeneous
polynomial in two variables.  Let $\G$ be a subgroup of $\PGL_2(F)$.
 If $f$ is $\G$-invariant then $\G$, acting by linear transformation on
 $\mathbb{P}^1(F)$, permutes the zero set of $f$, and given a finite $\G$-invariant subset of
 $\mathbb{P}^1(F)$, we can construct a $\G$-invariant form.   We will call a
 $\G$-invariant
 form
 \emph{$\G$-minimal} if its zeros are given by the $\G$-orbit of a single point and if each zero occurs with multiplicity $1$.  It is easy to see that a binary form is
 $\G$-invariant if and only if it is a product of $\G$-minimal forms.
 \begin{lem}\label{minform} Suppose that $\G$ is one of groups given in Lemma~\ref{subgroups}.
  Let $P\in\mathbb{P}^1(F)$.  Then the orbit of $P$ under the action of $\G$ is of size $|\G|$ unless
  $P$ is a zero of one of the following $\G$-minimal forms.
 \begin{enumerate}
    \item[(a)] $\G=\G_{C_{n}}$, $n>1$.
\[X_0,\quad X_1\]
    \item[(b)] $\G=\G_{D_{2n}}$, $n>1$.
     \[X_0X_1,\quad X_0^n-X_1^n,\quad X_0^n+X_1^n\]
    \item[(c)] $\G=\G_{A_4}$.
    \[X_0X_1(X_0^4-X_1^4),\quad X_0^4+2i\sqrt{3}X_0^2X_1^2+X_1^4,\quad X_0^4-2i\sqrt{3}X_0^2X_1^2+X_1^4\]
    \item[(d)] $\G=\G_{S_4}$.
    \[ X_0^{12}-33X_0^8X_1^4-33X_0^4X_1^8+X_1^{12},\quad X_0^8+14X_0^4X_1^4+X_1^8,\quad X_0X_1(X_0^4-X_1^4)\]
    \item [(e)]$\G=\G_{A_5}$.
    \[ X_0X_1(X_0^{10}+11X_0^5X_1^5-X_1^{10}),\]
    \[-(X_0^{20}+X_1^{20})+228(X_0^{15}X_1^5-X_0^5X_1^{15})-494X_0^{10}X_1^{10},\]
    \[(X_0^{30}+X_1^{30})+522(X_0^{25}X_1^5-X_0^5X_1^{25})-10005(X_0^{20}X_1^{10}+X_0^{20}X_1^{10})\]
    \item [(f)] $\G=\G_{\beta,A}$.
    \[ X_1\] and if $\beta\ne 1$
    \[\prod_{a\in A}(X_0-aX_1)\]
    \item [(g)]and~(h) $\G=\PSL_2(\F_q)$ or $\G=\PGL_2(\F_q)$
    \[(X_0^q -X_0X_1^{q-1})^{q-1}  + X_1^{q (q-1)},\quad (X_0^q-X_0X_1^{q-1}) X_1 \]

 \end{enumerate}
 \end{lem}
\begin{proof}  The $\G$-minimal forms listed in parts (b)-(e) are called ``Grundformen" by Weber in~\cite{Weber:2}.
  See \S70
of~\cite{Weber:2} for a discussion of ``Grundformen."  See \S71,
\S72, \S73, and \S76 of~\cite{Weber:2} for the proofs of (a)-(e).

  Let
$P\in\mathbb{P}^1(F)$ and suppose that the orbit of $P$ under the
action of $\G$ contains less than $|\G|$ points. It follows that for
some $M\in \G-\{Id\}$, $M(P)=P$.  If $M'\in \G$, then
$((M')^{-1}MM')((M')^{-1}(P))=(M')^{-1}(P)$. So any element in $\G$
which is conjugate to $M$ fixes a point in the orbit of $P$ under
the action of $\G$.

We now prove part~(f).  A computation shows that any nonidentity
element in $\G_{\beta, A}$ is conjugate to an element in the set
\[\left\{
\begin{bmatrix}
\beta^k&0\\
0&1\\
\end{bmatrix},
\begin{bmatrix}
1&a\\
0&1\\
\end{bmatrix}\colon k\in\{1,\ldots, r-1\}, a\in A-\{0\}\right\},\]
 where $r$ is the multiplicative order of $\beta$.  An
element of the form
\[\begin{bmatrix}
\beta^k&0\\
0&1\\
\end{bmatrix}\] fixes the points $[1\colon 0]$ and $[0\colon 1]$.  An
element in of the form
\[\begin{bmatrix}
1&a\\
0&1\\
\end{bmatrix}\] fixes the point $[1\colon 0]$.  Our result follows.

To prove parts~(g) and~(h), let $\ol{\F_q}$ be an algebraic closure
of $\F_q$, let $x=X_0/X_1$ and note that $\G$ acts by fractional
linear transformation on $\ol{\F_q}(x)$.  As shown in the proof of
Theorem~\ref{mainhyp}, the fixed field $\ol{\F_q}(x)^{\G}$ of
$\ol{\F_q}(x)$ under the action of $\G$ is $\ol{\F_q}(u)$, where
\[u:=\left(\frac{( (x^q -x)^{q-1}  + 1)^{\frac{q+1}{2}}}{ (x^q-x)^{q(\frac{q-1}{2})} }\right)^i \]
with $i=1$ if $\G=\PSL_2(\F_q)$ and $i=2$ if $\G=\PGL_2(\F_q)$.

  The primes that
 ramify in $\ol{\F_q}(x)/\ol{\F_q}(u)$ correspond to the
points of $\mathbb{P}^1(F)$ that have orbits of size less that
$|\G|$. That is, if $\mf{P}$ is a prime of $\ol{\F_q}(u)$ that
ramifies in $\ol{\F_q}(x)$ with ramification index $e$, say
\[\mf{P}=(\mf{Q}_1\ldots\mf{Q}_r)^e,\]
then the $\mf{Q}_i$ correspond to the element of an orbit under the
action of $\G$ that has  $r=|\G|/e$ elements.  By~Theorem~1
of~\cite{Valentini:1}, $\ol{\F_q}(u)$ has $2$ primes that ramify in
$\ol{\F_q}(x)$  with ramification degrees $|\G|/(q(q-1))$ and
$|\G|/(q+1)$. These primes correspond to $u=0$ and $u=\infty$
respectively. Our result follows.
\end{proof}
\section{Useful equations of hyperelliptic curves}

The following lemmas will be used to prove part of our main result
for plane curves.

\begin{lem}\label{fixedfield} Let $F$ be an algebraically closed field of characteristic not equal to $2$.
 Let $\G$ be a finite subgroup of $\PGL_2(F)$.  Then $\G$ acts by fractional linear
transformation on the field $F(x)$.  Using the notation of
Lemma~\ref{subgroups}, suppose
 that $\G\in\{\G_{D_{2n}}, \G_{A_4}, \G_{\PSL_2(\F_q)}\colon n>1\}$.  Then the subfield of $F(x)$ fixed
by the action of $\G$ is $F(t)$ where
\begin{enumerate}
  \item  \[t:=x^n+x^{-n}\] if $\G=\G_{D_{2n}}, n>1$.
  \item  \[t:=\frac{x^{12} - 33x^8 - 33x^4 + 1}{-x^{10} + 2x^6 - x^2}\]     if
  $\G=\G_{A_4}$.
  \item  \[t:=\frac{((x^q
  -x)^{q-1}+1)^{\frac{q+1}{2}}}{(x^q-x)^{\frac{q^2-q}{2}}}\]  if
  $\G=\G_{\PSL_2(\F_q)}$.
\end{enumerate}
\end{lem}
\begin{proof}  This is shown in the proof of Theorem~\ref{mainhyp}.
\end{proof}

\begin{lem}\label{hyp}
 Let $F$ be an algebraically closed field of characteristic not equal to
$2$.  Let $\G$ be as in Lemma~\ref{fixedfield}.  Then
\[y^2=f_{\alpha}(x)\] gives the equation of a hyperelliptic curve $X$ defined over $F$, with
$\Aut(X)/\la\iota\ra\cong \G$ where $\iota$ is the hyperelliptic
involution of $X$ and
\begin{enumerate}
  \item  \[f_{\alpha}(x):=x(x^4-1)(x^4-\alpha x^2+1),\; \alpha\in F^{\times}-\{\pm 2\}\] if $\G=\G_{D_{4}}$.
  \item
  \[f_{\alpha}(x):=x^{2n}-\alpha x^{n}+1,\; \alpha\in F^{\times}-\{\pm 2\}\] if $\G=\G_{D_{2n}}, n>2$.
  If $n=4$, assume also that $\alpha\ne 14$.
  If $n=3$, assume also that $\alpha\ne\pm\sqrt{-50}$.
  \item  \[f_{\alpha}(x):=x^{12} - 33x^8 - 33x^4 +1 +\alpha(x^{10} - 2x^6 + x^2),\; \alpha\in F^{\times}\]     if
  $\G=\G_{A_4}$.  If the characteristic of $F$ is not $5$, assume also that \[\alpha\ne\pm
  (22/5)(r^3+2r^2-5r-1)\] where $r$ is a zero of
  \[t^4+2t^3-6t^2-2t+1.\]
  \item  \[f_{\alpha}(x):=((x^q
  -x)^{q-1}+1)^{\frac{q+1}{2}}-\alpha (x^q-x)^{\frac{q^2-q}{2}}, \;\alpha\in F^{\times}\]   if
  $\G=\G_{\PSL_2(\F_q)}$.
\end{enumerate}
\end{lem}
\begin{proof}
Let $g_{\alpha}(X_0,X_1)$ be the homogenization of $f_{\alpha}(x)$
with even degree.  For example, if
$f_{\alpha}(x):=x(x^4-1)(x^4-\alpha x^2+1)$, then let
\[g_{\alpha}(X_0,X_1):=X_0X_1(X_0^4-X_1^4)(X_0^4-\alpha X_0^2X_1^2+X_1^4).\]  It can be
deduced from Lemmas~\ref{fixedfield} and~\ref{minform} that
$g_{\alpha}(X_0,X_1)$ is a squarefree $\G$-invariant form.  Then
$f_{\alpha}(x)$ is squarefree and the equation
\[y^2=f_{\alpha}(x)\] gives the equation of a hyperelliptic curve $X$.

Let $\iota$ be the hyperelliptic involution of $X$.  As discussed in
Chapter~3, $\Aut(X)/\la\iota\ra$ is given by a finite subgroup of
$\PGL_2(F)$. Let $\mf{H}=\Aut(X)/\la\iota\ra$.  Observe that
$\mf{H}=\Stab(g_{\alpha})$.
  {}From
our choice of $f_{\alpha}(x)$, it is clear that $\G\subseteq\mf{H}$.
We need to show $\G=\mf{H}$.  There are four cases.
\begin{enumerate}
  \item $\G=\G_{D_{4}}$.    Suppose that $\G\subsetneq \mf{H}$.
    Then by Lemma~\ref{propersubgroups}(a), $\mf{H}=\G_{A_4}$,
   $\mf{H}\cong\G_{S_4}$, $\mf{H}\cong
    \G_{A_5}$, $\mf{H}\cong D_{4n}$ for some $n>1$, $\mf{H}\cong\PSL_2(\F_q)$ where $q>3$,
    or $\mf{H}\cong\PGL_2(\F_q)$ for some finite field $\F_q$.

A computation shows that $g_{\alpha}$ is not
    $\G_{A_4}$-invariant, so $\mf{H}\ne \G_{A_4}$.

    By Lemma~\ref{minform},  if $\mf{H}\cong S_4$ then
    $\mf{H}$-minimal forms have degrees $6$, $8$, $12$, or
    $24$.  So a product of $\mf{H}$-minimal forms can never
    have degree equal to
    $10$.  Since the degree of $g_{\alpha}$
    is $10$, $\mf{H}\not\cong S_4$.

    If $\mf{H}\cong A_5$, then by Lemma~\ref{minform}, the
    degree of $g_{\alpha}$ must be greater than or equal to $12$.  So since the degree of $g_{\alpha}$
    is $10$, $\mf{H}\not\cong A_5$.

    By Lemma~\ref{minform}, $\mf{H}$-minimal forms have degrees $q+1$, $q(q-1)$, or
    $(q^3-q)/2$ if $\mf{H}\cong \PSL_2(\F_q)$ and have degrees $q+1$, $q(q-1)$, or
    $(q^3-q)$ if $\mf{H}\cong \PGL_2(\F_q)$.  Since the degree of $g_{\alpha}$
    is $10$, one can show that the cases $\mf{H}\cong\PSL_2(\F_q)$ where
    $q>3$,
    and $\mf{H}\cong\PGL_2(\F_q)$ for odd $q$ cannot occur.

   Assume that $\mf{H}\cong D_{4n}$ with $n>1$.  Then
   there is an element $A\in\mf{H}$ of order $2n$ with $A^n$
   equal to one of the elements of $\G$.  A computation shows that for any
   \[M:=\begin{pmatrix}
   a&b\\
   c&d\\
   \end{pmatrix}\in\GL_2(F),\] with $\ol{M}\in\G_{S_4}$ the form
   $(g_{\alpha}(X_0,X_1))M:=g_{\alpha}(aX_0+bX_1,cX_0+dX_1)$ is equal to
   $\lambda g_{\alpha'}(X_0,X_1)$ for some $\lambda\in F^{\times}$ and
   \[\alpha'\in\left\{\alpha, -\alpha, \frac{2\alpha+12}{\alpha-2},
   \frac{2\alpha-12}{-\alpha-2}, \frac{2\alpha-12}{\alpha+2},
    \frac{2\alpha+12}{-\alpha+2} \right\}.\]  Since all elements of
    order $2$ in $\G_{D_4}$ are conjugate by an element of
    $\G_{S_4}$, after replacing $g_{\alpha}$ with $(g_{\alpha})M:=g_{\alpha'}$
    where $M$ is an appropriate element of $\GL_2(F)$ such that $\ol
    M\in\G_{S_4}$,
    we may assume that
    \[A^n=\begin{bmatrix}
    -1&0\\
    0&1\\
    \end{bmatrix}.\]  Since $A$ is of order $2n$ and commutes with
    $A^n$, a computation shows that
    \[A=\begin{bmatrix}
    \zeta&0\\
    0&1\\
    \end{bmatrix},\] where $\zeta$ is a $2n^{th}$ root of unity.  It can be
    verified immediately that $g_{\alpha}$ is not $A$-invariant,
    so we get a contradiction.

Therefore $\mf{H}=\G_{D_4}$.
  \item $\G=\G_{D_{2n}}$, $n>2$.   Suppose that $\G\subsetneq \mf{H}$.
    Then $\mf{H}$ is given by one of the groups listed in
    Lemma~\ref{propersubgroups}(b)-(e).

    Since $\alpha\ne 0$, it can be shown with a computation that
    $\mf{H}$ is not of the form $\G_{D_{2n'}}$ for any
    $n'>n$.

    Suppose that $\mf{H}$ is isomorphic to $\PSL_2(\F_q)$ or
    $\PGL_2(\F_q)$ where $n|q-1$.
      By Lemma~\ref{minform}, $\mf{H}$-minimal forms have degrees $q+1$, $q(q-1)$, or
    $(q^3-q)/2$ if $\mf{H}\cong \PSL_2(\F_q)$ and have degrees $q+1$, $q(q-1)$, or
    $(q^3-q)$ if $\mf{H}\cong \PGL_2(\F_q)$.  Since $g_{\alpha}$ is
    squarefree,
     the degree of $g_{\alpha}$
    is $2n$, and since $2n$ is strictly less than $q(q-1)$, $(q^3-q)/2$, and
    $(q^3-q)$, by Lemma~\ref{minform} we must have $2n=q+1$.  Write
    $q-1=nm$, with $m\in\Z^+$.  Then $n(2-m)=2$.  This is a contradiction since we are assuming that
    $n>2$.

    By Lemma~\ref{propersubgroups}, the only possibilities left are: $n=3$ and $\mf{H}\cong S_4$ or $\mf{H}\cong
    A_5$, $n=4$ and $\mf{H}=\G_{S_4}$, or $n=5$ and $\mf{H}\cong
    A_5$.

    If $n=4$ and $\mathfrak{H}=\G_{S_4}$, a direct computation shows
    that $g_{\alpha}$ is $\G_{S_4}$-invariant if and only if
    $\alpha=14$.  Since we are assuming $\alpha\ne 14$ when $n=4$, $g_{\alpha}$
    is not $\G_{S_4}$-invariant.

    If $n=5$ and $\mf{H}\cong
    A_5$, then by Lemma~\ref{minform}, any $\mf{H}$-invariant
    form has degree greater than or equal to $12$.  So we get a
    contradiction.

    Suppose that $n=3$.  Then $g_{\alpha}$ has degree $6$.   By
    Lemma~\ref{minform}, it cannot be an $\mf{H}$-invariant
    form if $\mf{H}\cong A_5$.  Suppose that
    $\mf{H}\cong S_4$.  The group $\G_{D_6}$ acts on $F\cup\{\infty\}$ by fractional
    linear transformation.  The zeros of $f_{\alpha}$ consist of
    the $\G_{D_6}$-orbit $\mf{O}$ of a point $P$.  Write
    \[\mf{O}:=\{P,\omega P, \omega^2 P, 1/P, \omega/P, \omega^2/P\}\]
    where $\omega$ is a primitive cube root of unity.  By Lemma~\ref{minform}, there is
exactly one orbit  $\mf{O}':=\{0,\infty,\pm 1,\pm i\}$ of
$F\cup\{\infty\}$ under the action of $\G_{S_4}$ of size $6$.  One
can verify that for any element $g\in \G_{S_4}$ of order $3$, there
exists an element $h\in \G_{A_4}\subset\G_{S_4}$ of order $2$ and
$P',Q'\in \mf{O}'$ such that $\mf{O}'=\{P',g(P'),g^2(P'), Q', g(Q'),
g^2(Q')\}$, $h(P')=P'$, $h(Q')=Q'$,$h(g(P'))=g(Q')$, and
$h(g^2(P'))=g^2(Q')$.

Write
\[v:=\begin{bmatrix}
     \omega &0\\
     0&1\\
     \end{bmatrix}\in\mf{H}\] and note that $v$ has order $3$.
     Since $\mf{H}$ is conjugate to $\G_{S_4}$ there exists
$M\in\PGL_2(F)$ such that $M\mf{H}M^{-1}=\G_{S_4}$.  Let
$g=MvM^{-1}\in\G_{S_4}$.  We must have
\[\mf{O}'=\{M(P),M(\omega P), M(\omega^2 P), M(1/P), M(\omega/P), M(\omega^2/P)\}\]
\[=\{(M(P), g(M(P)),g^2(M(P)),M(1/P),g(M(1/P)),g^2(M(1/P))\}.\]  Let
$h\in\G_{S_4}$ be an
   element of order $2$ such that  $h(P')=P'$, $h(Q')=Q'$,$h(g(P'))=g(Q')$, and
$h(g^2(P'))=g^2(Q')$, where $P'=g^i(M(P))$ and $Q'=g^j(M(1/P))$ for
some $i, j$ with $0\le i,j\le 2$.  Let $u=M^{-1}hM\in\mf{H}$.  Then
$u(\omega^i P)=\omega^i P$, $u(\omega^j/P)=\omega^j/ P$,
$u(\omega^{i+1} P)=\omega^{j+1}/P$, and $u(\omega^{i+2}
P)=\omega^{j+2}/P$.  Replacing $P$ with $\omega^i P$, we may assume
that $u(P)=P$, $u(\omega^k/P)=\omega^k/P$, $u(\omega
P)=\omega^{k+1}/P$, and $u(\omega^2 P)=\omega^{k+2}/P$, for some $k$
with $0\le k\le 2$. Any element of order $2$ in $\PGL_2(F)$ is
conjugate to the image of the matrix
\[\begin{pmatrix}
-1 & 0\\
0 & 1\\
\end{pmatrix},\] and so fixes exactly $2$ points of $F\cup\{\infty\}$
and is the image of a matrix with trace $0$.  Since $u$ does not fix
$\infty$, $u$ is the image in $\PGL_2(F)$ of a matrix of the form
\[\begin{pmatrix}
a & b\\
1 & -a\\
\end{pmatrix}.\]  Write $Q:=\omega^k/P$.
Solving \[P=\frac{aP+b}{P-a}\] and \[Q=\frac{aQ+b}{Q-a}\] for $a$
and $b$ we see that $u$ is the image in $\PGL_2(\C)$ of
\[\begin{pmatrix}
     \frac{P+Q}{2}&-PQ\\
     1&-\frac{P+Q}{2}\\
     \end{pmatrix}.
     \]
Since
\[\zeta Q=\frac{\zeta P\left (\frac{P+Q}{2}\right )-PQ}{\zeta P-\frac{P+Q}{2}},\]
we have
\[Q^2+4PQ+P^2=0\]
Substituting $Q=\omega^k/P$, we obtain
\[(P/w^{2k})^4+4(P/\omega^{2k})^2+1=0.\]  A computation shows that the resultant of
\[f_{\alpha}(x)=x^6-\alpha x^3+1\] and
\[x^4+4x^2+1\] is
\[\alpha^4+100\alpha^2+2500.\]  Since $P/w^{2k}$ is a zero of $f_{\alpha}$, we must have
\[\alpha^4+100\alpha^2+2500=0.\]  It follows that $\alpha=\pm\sqrt{-50}$.  So we get a contradiction.

Therefore $\mf{H}=\G_{D_{2n}}$.

  \item $\G=\G_{A_{4}}$.    Suppose that $\G\subsetneq \mf{H}$.
    Then by Lemma~\ref{propersubgroups}, $\mf{H}=S_4$, $\G'\cong
    A_5$,
    $\G'\cong\PSL_2(\F_q)$,
    or $\G'\cong\PGL_2(\F_q)$ for some finite field $\F_q$.

    Since $\alpha\ne 0$, a computation shows that
    $\mf{H}\ne\G_{S_4}$.

    Suppose that $\mf{H}\cong A_5$.  Let
    \[N:=\begin{bmatrix}
    r&-1\\
    1&r\\
    \end{bmatrix},\]
    where $r$ is a zero of $t^4+2t^3-6t^2+1$.  A computation shows
    that $N\G_{A_4}N^{-1}\subset\G_{A_5}$.  By Lemma~\ref{subgroups}, there
    exists $M\in\PGL_2(F)$ such that $M\mf{H}M^{-1}=\G_{A_5}$.
     Since all subgroups of $A_5$ isomorphic to $A_4$ are conjugate
     in $A_5$,
     we may assume without loss of generality that
     $M\G_{A_4}M^{-1}=N\G_{A_4}N^{-1}$.  So
     $N^{-1}M$ is in $N(\G_{A_4})$.  By Lemma~\ref{normalizersub},
     $N(\G_{A_4})=\G_{S_4}$.  So $M=NA$, with $A\in\G_{S_4}$.
     By Lemma~\ref{minform}, any $\G_{A_5}$-invariant form of degree
     $12$ is a scalar multiple of
     \[h(X_0,X_1):=X_0X_1(X_0^{10}+11X_0^5X_1^5+X_1^{10}).\]
     Since $M^{-1}\G_{A_5}M=\mf{H}$, we must have
     \[([h])M=[g_{\alpha}],\] and so
     \[([h])M=([h])NA=[g_{\alpha}].\]  Then
     \[([h])N=([g_{\alpha}])A^{-1}.\]
     Since $A^{-1}$ is in $\G_{S_4}$, a computation shows that
     \[([g_{\alpha}])A^{-1}=[g_{\pm\alpha}].\]  Another computation shows
     that
     \[([h])N=[g_{\beta}],\] where
     \[\beta=22/5(r^3+2r^2-5r-1).\]  This implies that $\alpha=\pm 22/5(r^3+2r^2-5r-1)$, which is a contradiction
      by our choice of $\alpha$.

     Suppose that $\mf{H}$ is isomorphic to $\PSL_2(\F_q)$ or
    $\PGL_2(\F_q)$.   By Lemma~\ref{minform}, $\mf{H}$-minimal
     forms have degrees $q+1$, $q(q-1)$, or
    $(q^3-q)/2$ if $\mf{H}\cong \PSL_2(\F_q)$ and have degrees $q+1$, $q(q-1)$, or
    $(q^3-q)$ if $\mf{H}\cong \PGL_2(\F_q)$.  As in Lemma~\ref{subgroups}, the $\G_{A_4}$ case
    presumes that the characteristic
    of $F$ is greater than $3$.  Then  \[12\in\{q+1, q(q-1), (q^3-q)/2, q^3-q\}\]
    implies that
    $12=q+1$.  So
      $\mf{H}$ must be isomorphic to $\PSL_2(\F_{11})$ or
    $\PGL_2(\F_{11})$.

     Using Magma,
    one can show that any subgroup of $\PGL_2(\F_{11})$ isomorphic to
    $A_4$ is contained in a subgroup of $\PSL_2(\F_{11})$ isomorphic
    to $A_5$. It follows from the above argument showing
    $\mf{H}\not \cong A_5$, that we get a contradiction.

         Therefore, $\mf{H}=\G_{A_4}$.
    \item $\G=\G_{\PSL_2(\F_q)}$.    Suppose that $\G\subsetneq
    \mf{H}$.  Then by Lemma~\ref{subgroups},
    $\mf{H}\cong \PSL_2(\F_{q'})$ or $\PGL_2(\F_{q'})$ for some
    $q'$.  It is easily deduced that $q|q'$.  Write $q'=q^r$, $r>1$.  The degree of $g_{\alpha}$
    is $(q^3-q)/2$.   By Lemma~\ref{minform}, $\mf{H}$-minimal
     forms have degrees $q^r+1$, $q^r(q^r-1)$, or
    $(q^{3r}-q^r)/2$ if $\mf{H}\cong \PSL_2(\F_{q^r})$ and have degrees $q^r+1$, $q^r(q^r-1)$, or
    $(q^{3r}-q^r)$ if $\mf{H}\cong \PGL_2(\F_{q^r})$.  Since
    $l=(q^3-q)/2$ with $l\in\{q^r+1, q^r(q^r-1), (q^{3r}-q^r)/2, (q^{3r}-q^r)
    \}$ has no solution, we get a contradiction.

    Therefore,  $\mf{H}=\G_{\PSL_2(\F_q)}$.
\end{enumerate}

\end{proof}
\section{The main result for plane curves}
\begin{lem}\label{diagisom}  Let $K$ be a perfect field, let $F$ be an algebraic closure of $K$, and let
$\Gamma=\Gal(F/K)$. Let $X$ be a smooth plane curve of degree $d>3$
given by an equation \[f(X_0,X_1,X_2)=0\] and defined over $F$ with
$K$ as field of moduli relative to $F/K$.  Suppose there exists
$i\in\{0,1,2\}$ such that for all $\sigma\in\Gamma$ there exists an
isomorphism $X\to\el{\sigma}{X}$ given by
  \[X_i\mapsto \epsilon_{\sigma}X_i\]
  \[X_j\mapsto X_j\] for $j\in\{0,1,2\}-\{i\}$.
  Then $X$ can be defined over $K$.
\end{lem}
\begin{proof}  After a permutation of coordinates, we may assume that
for all $\sigma\in\Gamma$, there exists an isomorphism
$X\to\el{\sigma}{X}$ given by an element of the form
\[E_{\sigma}:=\begin{bmatrix}
1&0&0\\
0&1&0\\
0&0&\epsilon_{\sigma}\end{bmatrix}.\]  Write
\[f(X_0,X_1,X_2)=\sum_{i=0}^d a_if_i(X_0,X_1)X_2^{i},\]
where, for $0\le i\le d$, $a_i\in F^{\times}$ and $f_i(X_0,X_1)$ is
a form of degree $d-i$. For all $i$ with $f_i(X_0,X_1)\ne 0$, write
\[f_i(X_0,X_1):=\sum_{j=0}^{d-i} b_{i,j} X_0^jX_1^{d-i-j}.\]
Assume also that the $a_i$ were chosen in such a way that
$b_{i,j}=1$ for some $j$.
 Let $i_d$ be the greatest
$i\in\{0,\ldots,d\}$ such that $a_if_i(X_0,X_1)X_2^{i}\ne 0$.  Then
after multiplying $f$ by $1/a_{i_d}$ we may assume that
\[f(X_0,X_1,X_2)=c_{i_d}f_{i_d}(X_0,X_1)X_2^{i_d}+\sum_{i=0}^{i_d-1} c_if_i(X_0,X_1)X_2^{i}\]
 where $c_i=a_i/a_{i_d}$ for $0\le i\le i_d$.

For each $\sigma\in\Gamma$, we have
\[f(X_0,X_1,X_2)=\lambda_{\sigma} f^{\sigma}(X_0,X_1,\epsilon_{\sigma}X_2),\]
for some $\lambda_{\sigma}\in F^{\times}$.
 Since for each $\sigma\in\Gamma$
 \[f^{\sigma}(X_0,X_1,\epsilon_{\sigma}X_2)
 =\epsilon_{\sigma}^{i_d}(\sigma(c_{i_d})f_{i_d}^{\sigma}(X_0,X_1)X_2^{i_d}+\sum_{i=0}^{i_d-1}\epsilon_{\sigma}^{i-i_d} \sigma(c_i)f_i^{\sigma}(X_0,X_1)X_2^{i}),\]
 by our
choice of the $a_i$ we must have
$f_i^{\sigma}(X_0,X_1)=f_i(X_0,X_1)$ for all $i$.  So since $c_{i_d}=1$,  for each
$\sigma\in\Gamma$ we must have
$\lambda_{\sigma}=\epsilon_{\sigma}^{i_d}$. So for each
$\sigma\in\Gamma$, $\sigma(c_i)=\epsilon_{\sigma}^{i_d-i}c_i$ for
all $0\le i\le i_d$ with $f_i(X_0,X_1)\ne 0$.   Let $h$ be the
greatest common divisor of \[\{i_d-i\colon 0\le i\le i_d-1 \mbox{
and }f_i(X_0,X_1)\ne 0 \}.\] Then there exists $c$ in the
multiplicative group generated by
\[\{c_i\colon0\le i\le i_d-1 \mbox{ and }f_i(X_0,X_1)\ne 0
\}\]  such that $\sigma(c)=\epsilon_{\sigma}^h c$ for all
$\sigma\in\Gamma$.  Note that if $\zeta$ is an $h^{th}$ root of
unity then $f(X_0,X_1,\zeta X_2)=f(X_0,X_1, X_2)$.  Let $\zeta_h$ be
a primitive $h^{th}$ root of unity.  So the
element \[A:=\begin{bmatrix} 1&0&0\\
0&1&0\\
0&0&\zeta_h\\
\end{bmatrix}\]
is an automorphism of $X$.  Choose $\gamma\in F^{\times}$ so that
$\gamma^h=c$.  Then for all $\sigma\in\Gamma$ we have
$\sigma(\gamma)/\gamma=\epsilon_{\sigma}\zeta_h^k$ for some
$k\in\Z$.   For each $\sigma\in\Gamma$ define an isomorphism
$\varphi_{\sigma}\colon X\to\el{\sigma}{X}$ by \[
\begin{bmatrix}
 1&0&0\\
0&1&0\\
0&0&\sigma(\gamma)/\gamma\\
\end{bmatrix}\]  Then it is easily verified that the family
$\{\varphi_{\sigma}\}_{\sigma\in\Gamma}$ satisfies Weil's cocycle
condition of Theorem~\ref{weil}.  So by Theorem~\ref{weil}, $X$ is
definable over $K$.
\end{proof}

\begin{prop}\label{intransitive}  Let $K$ be a perfect field of characteristic
$p\ne 2$ and let $F$ be an algebraic closure of $K$.  Let $X$ be a
smooth plane curve of degree $d>3$ defined over $F$.  Suppose that
$\Aut(X):=\mf{P}_A\rtimes\mf{I}$ is a group of the type given in
Lemma~\ref{gps}(r) or an intransitive group whose image in
$\PGL_2(F)$ is not equal to one of the groups listed in
Lemma~\ref{subgroups}(a) or~(f) (so we allow $\mf{P}_A=\{Id\}$.)
Then $X$ can be defined over its field of moduli relative
  to the extension $F/K$.
\end{prop}
\begin{proof}  By~Proposition~\ref{closed}, we may assume that $K$ is the field of
moduli of $X$ relative to the extension $F/K$.  By
Corollary~\ref{finitefield}, we may assume that $K$ is infinite. Let
$\Gamma=\Gal(F/K)$.

Fix an equation $f(X_0,X_1,X_2)=0$ for $X$.   Write
\[f(X_0,X_1,X_2)=\sum_{i=0}^d f_i(X_0,X_1)X_2^{d-i},\]
where for $1\le i\le d$, $f_i(X_0,X_1)$ is a form of degree $i$.
Let $A$ be in $\mf{I}\subseteq \mf{P}_A\rtimes\mf{I}$.
  Write
\[A:=\begin{bmatrix}
a&b&0\\
c&d&0\\
0&0&1\\
\end{bmatrix}.\]
Then since $A$ gives an automorphism of $X$, we must have
\[\sum_{i=0}^d f_i(aX_0+bX_1,cX_0+dX_1)X_2^{d-i}=\lambda
f(X_0,X_1,X_2)\] for some $\lambda\in F^{\times}$.  So we must have
\[f_i(aX_0+bX_1,cX_0+dX_1)=\lambda f_i(X_0,X_1)\] for all $i$.  Since $A$ is an arbitrary
element of $\mf{I}$, it follows that for all $i$, $f_i(X_0,X_1)$ is
an $\ol {\mf{I}}$-invariant form.

The function field $F(X)$ of $X$ is given by $F(x,z)$ with the
relation $f(x,1,z)=0$. The inclusion $F(x)\hookrightarrow F(x,z)$
corresponds to a map $X\to\mathbb{P}^1_F$.  If $A\in\mf{I}\subseteq
\mf{P}_A\rtimes\mf{I}$ is an automorphism of $X$ then $A$ gives an
automorphism of $F(x,z)$ fixing $F$ and the restriction of $A$ to
$F(x)$ induces an automorphism of $F(x)$, fixing $F$, given by the
image $\ol A$ of $A$ in $\PGL_2(F)$.

By Lemma~\ref{conjcurve}(a) and~(r), for all $\sigma\in\Gamma$ there
exists an isomorphism $X\to\el{\sigma}{X}$  given by
$M_{\sigma}\in\PGL_3(F)$ where $M_{\sigma}$ is an element of an
intransitive group $\mf{I}'$ with $\ol{M_{\sigma}}\in
N(\ol{\mf{I}})$, the normalizer of $\ol{\mf{I}}$ in $\PGL_2(F)$.
{}From now on, unless otherwise stated,
 when an isomorphism $X\to\el{\sigma}{X}$ is mentioned, we will
 assume that it is given by a lift to $\PGL_3(F)$ of an element of
 $N(\ol{\mf{I}})$.

Suppose that $M_{\sigma}$ gives an isomorphism $X\to\el{\sigma}{X}$.
 Then $M_{\sigma}$ induces an automorphism of $\mathbb{P}^1_F$ given
by the image $\ol{M_{\sigma}}$ of $M_{\sigma}$ in $\PGL_2(F)$.  We
show that there exists $\mu\in\PGL_2(F)$ such that for all
$\sigma\in\Gamma$, $(\mu^{-1})^{\sigma}\mu\in\PGL_2(F)$ lifts to an
element of $\PGL_3(F)$ that gives an isomorphism
$X\to\el{\sigma}{X}$.

If $\ol{\mf{I}}=\G_{S_4}$, $\G_{A_5}$, or $\PGL_2(\F_q)$ for some
finite field $\F_q$ then by Lemma~\ref{normalizersub},
$N(\ol{\mf{I}})=\ol{\mf{I}}$. So for all $\sigma\in\Gamma$, the
identity lifts to an isomorphism $X\to\el{\sigma}{X}$.  So if
$\ol{\mf{I}}=\G_{S_4}$, $\G_{A_5}$, or $\PGL_2(\F_q)$ for some
finite field $\F_q$, then we may take $\mu=Id$. So assume $\G$ is in
$\{\G_{D_{2n}}, \G_{A_4}, \PSL_2(\F_q)\colon n>1\}$. Let $t$ be as
in Lemma~\ref{fixedfield}. So the subfield of $F(x)$ fixed by
$\ol{\mf{I}}$ acting by fractional linear transformation on $F(x)$
is $F(t)$.

 Then the inclusion
$F(t)\hookrightarrow F(x,z)$ corresponds to a map $\rho\colon X\to
C$, where $C$ is a curve of genus $0$.  An isomorphism
$\varphi_{\sigma}\colon X\to\el{\sigma}{X}$ given by $M_{\sigma}$
with $M_{\sigma}\in N(\ol{\mf{I}})$  induces an isomorphism
$\tilde{\varphi_{\sigma}}\colon C\to\el{\sigma}{C}$ that makes the
following diagram commute:
\[
\begin{CD}
{X} @> {\varphi_{\sigma}} >> {\el{\sigma}{X}} \\
@ V{\rho} VV @ VV {\rho^{\sigma}} V \\
{C} @>> {\tilde{\varphi_{\sigma}}} >
{\el{\sigma}{C}} \\
\end{CD}
\]
For any $\sigma\in\Gamma$, any two isomorphisms  from $X$ to
$\el{\sigma}{X}$,  given by $M_{\sigma}$, $M_{\sigma}'$, are equal
up to multiplication by an element of $\mf{I}$.  This implies that
for each $\sigma\in\Gamma$ there exists a unique isomorphism
$\tilde{\varphi_{\sigma}}\colon C\to\el{\sigma}{C}$ making the above
diagram commute. Then the family
$\{\tilde{\varphi_{\tau}}\}_{\tau\in\Gamma}$ satisfy Weil's cocycle
condition
$\tilde{\varphi_{\tau}}^{\sigma}\,\tilde{\varphi_{\sigma}}=\tilde{\varphi_{\sigma\tau}}$
given in Theorem~\ref{weil}.  Let $B$ be the $K$-model associated
with $\{\tilde{\varphi_{\tau}}\}_{\tau\in\Gamma}$ and let
$\psi\colon C\to B_F$ be an isomorphism such that
$(\psi^{-1})^{\sigma}\psi=\tilde{\varphi_{\sigma}}$.

 It can be verified by the proof of
Theorem~\ref{mainhyp} that $B$ has a $K$-rational point. Since $B$
has genus $0$ and since $K$ is infinite, $B$ has infinitely many
rational points.   As shown in Chapter~4, we have a homomorphism
$\Gamma\to\Aut(F(t)/K)$  given by $\sigma\mapsto\sigma^*$ where
\[\sigma^*(t)=\frac{at+b}{ct+d},\; {\sigma^*(\alpha)}=\sigma(\alpha),\; \alpha\in F,\]
 where $\tilde{\varphi_{\sigma}}$ is given by
\[t\mapsto \frac{at+b}{ct+d}.\]
The curve $B$ is the variety over $K$ corresponding to the fixed
field of $\Gamma^*=\{\sigma^*\}_{\sigma\in\Gamma}$.

Since $C=\mathbb{P}^1_F$, we can identify $C(F)$ with
$F\cup\{\infty\}$.  Note that $P\in B(K)$ if and only if for all
$\sigma\in\Gamma$,
\[\sigma(Q)=\tilde{\varphi_{\sigma}}(Q)=\frac{aQ+b}{cQ+d},\]
where $\psi(Q)=P$ and $\tilde{\varphi_{\sigma}}$ is given by
\[t\mapsto \frac{at+b}{ct+d}.\]

  Let $Q$ be a point of
$C(F)$ such that $\psi(Q)$ is a $K$-rational point of $B$, such that
$\tilde{\varphi_{\sigma}}(Q)=Q$ if and only if $\sigma^*(t)=t$, and
such that $y^2=f_Q(x)$, where $f_Q(x)$ is given in Lemma~\ref{hyp},
gives the equation of a hyperelliptic curve $Y$ with
$\Aut(Y)/\la\iota\ra\cong\ol{\mf{I}}$.   Then the function field of
$Y$, $F(Y)$ is a quadratic extension of $F(x)$ and the function
field of $Y/\Aut(Y)$ is equal to $F(t)$.  For all $\sigma\in\Gamma$,
the curve $\el{\sigma}{Y}$ is given by the equation
$y^2=f_{\sigma(Q)}$.  It is easily verified that if $\sigma$ is in
$\Gamma$ and if $M_{\sigma}$ is in $\PGL_3(F)$ gives an isomorphism
$X\to\el{\sigma}{X}$, then $\ol{M_{\sigma}}\in\PGL_2(F)$ lifts to an
isomorphism $Y\to\el{\sigma}{Y}$.  So $K$ is the field of moduli of
$Y$ relative to the extension $F/K$.  By Proposition~\ref{canform},
there exists $\mu\in\PGL_2(F)$ such that for all $\sigma\in\Gamma$,
$(\mu^{-1})^{\sigma}\mu$ lifts to an isomorphism
$Y\to\el{\sigma}{Y}$. By construction, $(\mu^{-1})^{\sigma}\mu$
lifts to an isomorphism $X\to\el{\sigma}{X}$.

Let
\[M:=\begin{bmatrix}
a&b&0\\
c&d&0\\
0&0&1\\
\end{bmatrix}\in\PGL_3(F)\] be any lift of $\mu$.  So for each $\sigma\in\Gamma$ then there exists
\[E_{\sigma}:=\begin{bmatrix}
1&0&0\\
0&1&0\\
0&0&\epsilon_{\sigma}\\
\end{bmatrix}\in\PGL_3(F)\] such that $E_{\sigma}(M^{-1})^{\sigma}M$ gives an isomorphism
$X\to\el{\sigma}{X}$. Let $X'$ be the plane curve given by
\[f(dX_0-bX_1,-cX_0+aX_1,(ad-bc)X_2)=0,\]
so $M$ gives an isomorphism $X\to X'$.  So for each
$\sigma\in\Gamma$,
\[M^{\sigma}(E_{\sigma}(M^{-1})^{\sigma}M) M^{-1}=E_{\sigma}\]
gives an isomorphism $X'\to\el{\sigma}{X'}$.

By Lemma~\ref{diagisom}, $X'$ can be defined over $K$.  Therefore,
$X$ can be defined over $K$.
\end{proof}

\begin{lem}\label{split}  Let $K$ be a perfect field of characteristic
$p\ne 2$, let $F$ be an algebraic closure of $K$, and let
$\Gamma=\Gal(F/K)$. Let $X$ be a smooth plane curve of degree $d>3$
defined over $F$.  Suppose that
  $\Aut(X)$ is given by a group $\G$ listed in Lemma~\ref{gps}.
  Suppose that for all $\sigma\in\Gamma$, $\G^{\sigma}=\G$, the
  normalizer $N(\G)$ of $\G$ in $\PGL_3(F)$ is of the form
  $\G\rtimes \mf{H}$ for some subgroup $\mf{H}\subset\PGL_3(F)$, and that
  $\mf{H}=\mf{H}^{\sigma}$ for all $\sigma\in\Gamma$.
  Then $X$ can be defined over its field of moduli relative to the
  extension $F/K$.
\end{lem}
\begin{proof}  By~Proposition~\ref{closed}, we may assume that $K$ is the field of
moduli of $X$ relative to the extension $F/K$.  By
Lemma~\ref{isnorm}, for all $\sigma\in\Gamma$ any isomorphism
$X\to\el{\sigma}{X}$ is given by an element of $N(\G)=\G\rtimes
\mf{H}$.  Since any two isomorphisms are equal up to composition
with an element of $\Aut(X)$, for each $\sigma\in\Gamma$ there is
exactly one element of $\mf{H}$ that gives an isomorphism
$X\to\el{\sigma}{X}$.

For each $\sigma\in\Gamma$, let $f_{\sigma}$ be the isomorphism
$X\to\el{\sigma}{X}$ given by an element of $\mf{H}$.  Since
$\mf{H}^{\sigma}=\mf{H}$ for all $\sigma\in\Gamma$,
 we must have
\[f_{\sigma}^{\tau}f_{\tau}=f_{\tau\sigma}\]
for all $\sigma$, $\tau\in\Gamma$.

So the family $\{f_{\sigma}\}_{\sigma\in\Gamma}$ satisfies Weil's
cocycle condition given in Theorem~\ref{weil}.  By
Theorem~\ref{weil}, $X$ can be defined over $K$.
\end{proof}

\begin{prop}\label{typec}  Let $K$ be a perfect field of characteristic
$p\ne 2$ and let $F$ be an algebraic closure of $K$.  Let $X$ be a
smooth plane curve of degree $d>3$ defined over $F$.  Suppose that
$\Aut(X)$ is a group of the type $\mf{C}$ given in
Lemma~\ref{gps}(b).
      Then $X$ can be defined over its field of moduli relative
  to the extension $F/K$.
\end{prop}
\begin{proof}By~Proposition~\ref{closed}, we
may assume that $K$ is the field of moduli of $X$ relative to the
extension $F/K$. Let $\Gamma=\Gal(F/K)$.  Let $\G:=\Aut(X)$. We
follow the notation used in Lemma~\ref{normalizergps}.  There are
three cases.
\begin{enumerate}

  \item[i.] Suppose $\G=\G_9$.  By
Lemma~\ref{conjcurve}, for all $\sigma\in\Gamma$, any isomorphism
$X\to\el{\sigma}{X}$ is given by an element of $N(\G_9)$.  By
Lemma~\ref{normalizergps}, $N(\G_9)=\G_9\rtimes\la US^2, V\ra$.  One
can verify that for all $\sigma\in\Gamma$, $\la US^2, V\ra=\la US^2,
V\ra^{\sigma}$.  So by Lemma~\ref{split}, $X$ can be defined over
$K$.

 \item[ii.]  Suppose $S_{m,1}\in\mf{H}$ and $\mf{H}\ne\la S_{3,2}\ra$.  By
  Lemma~\ref{conjcurve} and Lemma~\ref{normalizergps}, for all
  $\sigma\in\Gamma$, any isomorphism
$X\to\el{\sigma}{X}$ is given by an element of $N(\G)\subseteq\la\G,
R, S_{3m,2}\ra$.    Suppose that $[1\colon 0\colon 0]$ is in $X(F)$.
 Since $\G$ is normal in $N(\G)$, every element of $N(\G)$ can be written as
   an element of $\G$ times an element of $\la R,S_{3m,2}\ra$.
Any element of $\la R,S_{3m,2}\ra$ fixes $[1\colon 0\colon 0]$.
 For all $\sigma\in\Gamma$ we have $\Aut(\el{\sigma}{X}) = \G^{\sigma} = \G$.  Thus for all $\sigma\in\Gamma$, any isomorphism $X\to\el{\sigma}{X}$
maps the point $[1\colon 0\colon 0]$  into the
$\Aut(\el{\sigma}{X})$-orbit of  $[1\colon 0\colon 0]$.  So if the
point $[1\colon 0\colon 0]$ is in $X(F)$, it will map to a rational
point of the canonical $K$-model of $X/\Aut(X)$ given in
Theorem~\ref{bk}.  So in this case, by Corollary~\ref{pqpt}, $X$ can
be defined over $K$.

{}From now on we assume that $[1\colon 0\colon 0]$ is not in $X(F)$.
 Recall that
\[S_{3m,2}=\begin{bmatrix}
\zeta_{3m}&0&0\\
0&\zeta_{3m}^{-1}&0\\
0&0&1\\
\end{bmatrix},\; T=\begin{bmatrix}
0&1&0\\
0&0&1\\
1&0&0\\
\end{bmatrix},\;\mbox{and }R=\begin{bmatrix}
1&0&0\\
0&0&1\\
0&1&0\\
\end{bmatrix},\] where $T\in\G$ and $\zeta_{3m}$ is a primitive $(3m)^{th}$ root of unity where $m=3^l$ with
$l\ge 0$.  Define
\[S:=TRS_{3m,2}(TR)^{-1}=\begin{bmatrix}
1&0&0\\
0&\zeta_{3m}^{-1}&0\\
0&0&\zeta_{3m}\\
\end{bmatrix}.\]  So $N(\G)\subseteq\la\G,
R, S_{3m,2}\ra=\la\G, R, S\ra$.  We have $R^2=Id$, $S^3\in\G$, and
$RSR^{-1} = S^{-1}$,
 so we have a surjection $S_3 \to\la\G,R,S\ra/\G$,
 where $S_3$ is the symmetric group on $3$ letters.
 This surjection is an isomorphism,
 since every nontrivial normal subgroup of $S_3$ contains $A_3$,
 but $S\not\in\G$.  It follows that for any $\sigma\in\Gamma$, an isomorphism
$X\to\el{\sigma}{X}$ is given by $S^sR^r$, for some $s$ and $t$ with
$0\le s\le 2$ and $0\le r\le 1$.  Furthermore, it follows that for
$0\le s_1, s_2\le 2$ and $0\le r_1, r_2\le 1$ the element
$(S^{s_1}R^{r_1})^{-1}S^{s_2}R^{r_2}$ is an automorphism of $X$ if
and only if $s_1=s_2$ and $r_1=r_2$. Thus $s$ and $r$ are uniquely
determined by $\sigma$.

Fix an equation
\[f(X_0,X_1,X_2)=0\] for $X$.  Since $[1\colon 0\colon 0]$ is not in $X(F)$,
after multiplying $f(X_0,X_1,X_2)$ by an element of $F^{\times}$ we
may assume that \[f(X_0,X_1,X_2)=X_0^d+\sum_{i=0}^{d-1}
f_i(X_1,X_2)X_0^{i},\] where $f_i(X_1,X_2)$ is a form of degree
$d-i$.  For all $i$, write
\[f_i(X_1,X_2):=\sum_{j=0}^{d-i} b_{i,j} X_1^jX_2^{d-i-j}.\]
Since the elements
\[A_1:=\begin{bmatrix}
1&0&0\\
0&\zeta_m&0\\
0&0&1\\
\end{bmatrix}\;\mbox{and }A_2:=\begin{bmatrix}
1&0&0\\
0&1&0\\
0&0&\zeta_m\\
\end{bmatrix}\] are automorphisms of $X$, for all $i,j$ with
$b_{i,j}\ne 0$, we must have $j\equiv 0\pmod{m}$ and
$d-i-j\equiv0\pmod{m}$.  Since $S$ is not an isomorphism of $X$,
there exists $i_1$ and $j_1$ such that $b_{i_1,j_1}\ne 0$ and such
that $2j_1+i_1-d\not\equiv 0\pmod{3m}$.  Let $a=b_{i_1,j_1}$ and let
$b=b_{i_1,(d-i_1-j_1)}$. Note that since $2j_1+i_1-d\equiv
0\pmod{m}$, we have $(\zeta_{3m}^{2j_1+i_1-d})^3=1$ and since
$2j_1+i_1-d\not\equiv 0\pmod{3m}$ we have
$\zeta_{3m}^{2j_1+i_1-d}\ne 1$. So $\zeta_{3m}^{(2j_1-d+i_1)}$ is a
primitive cube root of unity. Note that
\[\zeta_{3m}^{(2(d-i_1-j_1)-d+i_1)}=\zeta_{3m}^{-(2j_1-d+i_1)}.\]

 Suppose that $a^3=b^3$, say $b=\omega a$ where $\omega^3=1$.  It
is shown in the proof of Lemma~\ref{normalizergps}(b) that
$S_{3m,2}\in N(\G)$.  So $S\in N(\G)$.  If necessary, after
replacing $f(X_0,X_1,X_2)$ with
$f(X_0,\zeta_{3m}^nX_1,\zeta_{3m}^{-n}X_2)$ with an appropriate
$n\in\{1,2\}$,  we may assume that if $a=b$.

Fix $\sigma\in\Gamma$ and suppose an
 isomorphism $X\to\el{\sigma}{X}$ is
 given by  $S^sR^r$, where
$0\le s\le 2$ and $0\le r\le 1$.

 We have
\[f^{\sigma}(X_0,X_1,X_2)=\lambda_{\sigma}f( X_0,(\zeta_{3m}^s)^{(-1)^r}X_{r+1},(\zeta_{3m}^{-s})^{(-1)^r}X_{2-r})\] for some
$\lambda_{\sigma}\in F^{\times}$. Since
\[f( X_0,(\zeta_{3m}^s)^{(-1)^r}X_{r+1},(\zeta_{3m}^{-s})^{(-1)^r}X_{2-r})\]
\[=X_0^d+\sum_{i=0}^{d-1}
f_i((\zeta_{3m}^s)^{(-1)^r}X_{r+1},(\zeta_{3m}^{-s})^{(-1)^r}X_{2-r})X_0^{i},\]
we must have $\lambda_{\sigma}=1$.  So
\[f^{\sigma}(X_0,X_1,X_2)=X_0^d+\sum_{i=0}^{d-1}
f_i((\zeta_{3m}^s)^{(-1)^r}X_{r+1},(\zeta_{3m}^{-s})^{(-1)^r}X_{2-r})X_0^{i}.\]
For all $i$, we have
\[f_i^{\sigma}(X_1,X_2)=\sum_{j=0}^{d-i}(\zeta_{3m}^{s(2j-d+i)})^{(-1)^r} b_{i,j} X_{r+1}^jX_{2-r}^{d-i-j}.\]
  If $r=0$, we have
\[\sigma(a)=\zeta_{3m}^{s(2j_1-d+i_1)} a\] and
\[\sigma(b)=\zeta_{3m}^{-s(2j_1-d+i_1)}b\]. If $r=1$ then
\[\sigma(a)=\zeta_{3m}^{s(2j_1-d+i_1)} b\] and
\[\sigma(b)=\zeta_{3m}^{-s(2j_1-d+i_1)}a.\] There are three cases to
consider.

Suppose that $a^3=b^3$.  As shown earlier, we have $a=b$. Since
$\zeta_{3m}^{(2j_1-d+i_1)}$ is a primitive cube root of unity, we
must have $s=0$. Since $\sigma$ is arbitrary, this implies that any
isomorphism $X\to\el{\tau}{X}$ with $\tau\in\Gamma$ is given by an
element of $\la R\ra$.  For $\tau\in\Gamma$ define $\varphi_{\tau}$
by $R$ if $\el{\tau}{X}\ne X$ and by $R^2=Id$ if $\el{\tau}{X}= X$.
The family $\{\varphi_{\tau}\}_{\tau\in\Gamma}$ satisfies Weil's
cocycle condition of Theorem~\ref{weil}.  So by Theorem~\ref{weil},
$X$ is definable over $K$.

  Now suppose that $a^3\ne b^3$ and that $b\ne 0$.  Write $2j_1-d+i_1=3^{l}u$ where $u\in\Z$
   and $(u,3)=1$.  Choose
$\alpha,\beta\in F$ such that $\alpha^{3^{l}}=a^u$ and
$\beta^{3^{l}}=b^u$.  If $r=0$, then
\[\sigma(\alpha^{3^{l}})=\zeta_{3m}^{s(2j_1-d+i_1)u}\alpha^{3^{l}}=(\zeta_{3m}^{s3^{l}})^{u^2}\alpha^{3^{l}}=\zeta_{3m}^{s3^{l}}\alpha^{3^{l}}\] and
\[\sigma(\beta^{3^{l}})=\zeta_{3m}^{-s(2j_1-d+i_1)u}\beta^{3^{l}}=(\zeta_{3m}^{-s3^{l}})^{u^2}\beta^{3^{l}}=\zeta_{3m}^{-s3^{l}}\beta^{3^{l}}.\]
It follows that
$\frac{\sigma(\alpha)}{\alpha}=(\zeta_{3m})^{3q_2}\zeta_{3m}^s$ and
$\frac{\sigma(\beta)}{\beta}=(\zeta_{3m})^{3q_1}\zeta_{3m}^{-s}$ for
some $q_1,q_2\in\Z$.

If $r=1$, then
\[\sigma(\alpha^{3^{l}})=\zeta_{3m}^{s(2j_1-d+i_1)u}\beta^{3^{l}}=(\zeta_{3m}^{s3^{l}})^{u^2}\beta^{3^{l}}=\zeta_{3m}^{s3^{l}}\beta^{3^{l}}\] and
\[\sigma(\beta^{3^{l}})=\zeta_{3m}^{-s(2j_1-d+i_1)u}\alpha^{3^{l}}=(\zeta_{3m}^{-s3^{l}})^{u^2}\alpha^{3^{l}}=\zeta_{3m}^{-s3^{l}}\alpha^{3^{l}}.\]
It follows that
$\frac{\sigma(\beta)}{\alpha}=(\zeta_{3m})^{3p_1}\zeta_{3m}^{-s}$
and $\frac{\sigma(\alpha)}{\beta}=(\zeta_{3m})^{3p_2}\zeta_{3m}^{s}$
for some $p_1,p_2\in\Z$.  Recall that $A_1, A_2\in\G$.  For some
$a_1,a_2\in\Z$, we have
\[S^sR^rA_1^{a_1}A_2^{a_2}=\begin{bmatrix}
1&0&0\\
0& \frac{\sigma(a^3)-b^3}{a^3-b^3}\frac{\sigma(\beta)}{\beta}&\frac{\sigma(a^3)-a^3}{b^3-a^3}\frac{\sigma(\beta)}{\alpha}\\
0&\frac{\sigma(b^3)-b^3}{a^3-b^3}\frac{\sigma(\alpha)}{\beta}&\frac{\sigma(b^3)-a^3}{b^3-a^3}\frac{\sigma(\alpha)}{\alpha}\\
\end{bmatrix}.\]
 For each $\tau\in\Gamma$, define
$\varphi_{\tau}$ by
\[\begin{bmatrix}
1&0&0\\
0& \frac{\tau(a^3)-b^3}{a^3-b^3}\frac{\tau(\beta)}{\beta}&\frac{\tau(a^3)-a^3}{b^3-a^3}\frac{\tau(\beta)}{\alpha}\\
0&\frac{\tau(b^3)-b^3}{a^3-b^3}\frac{\tau(\alpha)}{\beta}&\frac{\tau(b^3)-a^3}{b^3-a^3}\frac{\tau(\alpha)}{\alpha}\\
\end{bmatrix}.\]
Then for each $\tau\in\Gamma$, $\varphi_{\tau}$ is an isomorphism
$X\to\el{\tau}{X}$.  It can be verified that family
$\{\varphi_{\tau}\}_{\tau\in\Gamma}$ satisfies Weil's cocycle
condition of Theorem~\ref{weil}.  So by Theorem~\ref{weil}, $X$ is
definable over $K$.

Now assume that $b=0$.  Since $a\ne 0$, we must have $r=0$.  Since
$\sigma$ is arbitrary, this implies that any isomorphism
$X\to\el{\tau}{X}$ with $\tau\in\Gamma$ is given by $S^{s'}$ where
$0\le s'\le 2$. Let $\alpha$ be defined as in the case where $a^3\ne
b^3$ and $b\ne 0$. For some $a_1,a_2\in\Z$ we have
\[S^sA_1^{a_1}A_2^{a_2}=\begin{bmatrix}
1&0&0\\
0&\frac{\alpha}{\sigma(\alpha)}&0\\
0&0&\frac{\sigma(\alpha)}{\alpha}\\
\end{bmatrix}.\] For each
$\tau\in\Gamma$, define $\varphi_{\tau}$ by
\[\begin{bmatrix}
1&0&0\\
0&\frac{\alpha}{\tau(\alpha)}&0\\
0&0&\frac{\tau(\alpha)}{\alpha}\\
\end{bmatrix}.\]  Then $\varphi_{\tau}$ gives an isomorphism
$X\to\el{\tau}{X}$ and it can be verified that family
$\{\varphi_{\tau}\}_{\tau\in\Gamma}$ satisfies Weil's cocycle
condition of Theorem~\ref{weil}.  So by Theorem~\ref{weil}, $X$ is
definable over $K$.

 \item[iii.] Suppose $S_{m,1}\not\in\mf{H}$ and $\mf{H}\ne\la S_{3,2}\ra$.  By
  Lemma~\ref{conjcurve} and Lemma~\ref{normalizergps}, for all
  $\sigma\in\Gamma$, any isomorphism
$X\to\el{\sigma}{X}$ is given by an element of $N(\G)\subseteq\la\G,
R, S_{m,1}\ra$.  Suppose that $[1\colon 0\colon 0]$ is in $X(F)$.
 Then for all $\sigma\in\Gamma$,
any isomorphism $X\to\el{\sigma}{X}$ maps the point $[1\colon
0\colon 0]$  into the $\Aut(\el{\sigma}{X})$-orbit of  $[1\colon
0\colon 0]$.  So if the point $[1\colon 0\colon 0]$ is in $X(F)$, it
will correspond to a rational point of the canonical $K$-model of
$X/\Aut(X)$ given in Theorem~\ref{bk}.  So in this case, by
Corollary~\ref{pqpt}, $X$ can be defined over $K$.

{}From now on we assume that $[1\colon 0\colon 0]$ is not in $X(F)$.
Fix an equation
\[f(X_0,X_1,X_2)=0\]
for $X$.  Let $Y$ be the smooth plane curve defined over $F$ given
by
\[f(X_0,-X_1+X_2,X_1+X_2)=0.\] So an isomorphism $X\to Y$ is given by
\[M:=\begin{bmatrix}
1&0&0\\
0&-1/2&1/2\\
0&1/2&1/2\\
\end{bmatrix}\] and $\Aut(Y)$ is given by $M\G M^{-1}$.   It is enough to show that $Y$ is definable over $K$.
 We see that $[1\colon 0\colon 0]$ is not in $Y(F)$.

For all $\sigma\in\Gamma$, if $E_{\sigma}\in\PGL_3(F)$ gives an
isomorphism $X\to\el{\sigma}{X}$ then
\[M^{\sigma}E_{\sigma}M^{-1}=ME_{\sigma}M^{-1}\] gives an isomorphism
$Y\to\el{\sigma}{Y}$. For any $\sigma\in\Gamma$ there exists an
isomorphism $X\to\el{\sigma}{X}$ given by $S_{m,1}^sR^r$ with $0\le
s\le 2$ and $0\le r\le 1$.  A computation shows that
$MS_{m,1}M^{-1}=S_{m,1}$ and another computation shows that
$MRM^{-1}=D$ where
\[D:=\begin{bmatrix}
1&0&0 \\
0&-1&0\\
0&0&1\\
\end{bmatrix}.\]  So for any $\sigma\in\Gamma$ there exists an
isomorphism $Y\to\el{\sigma}{Y}$ given by $S_{m,1}^sD^r$ with $0\le
s\le 2$ and $0\le r\le 1$.   Furthermore, it is easily verified that
for $0\le s_1, s_2\le 2$ and $0\le r_1, r_2\le 1$ the element
$(S_{m,1}^{s_1}D^{r_1})^{-1}S_{m,1}^{s_2}D^{r_2}$ gives and
automorphism of $Y$ if and only if $s_1=s_2$ and $r_1=r_2$. Thus $s$
and $r$ are uniquely determined by $\sigma$.

Let $g(X_0,X_1,X_2):=f(X_0,-X_1+X_2,X_1+X_2)$.  After multiplying
$g(X_0,X_1,X_2)$ by an element of $F^{\times}$, we may assume that
\[g(X_0,X_1,X_2)=X_0^d+\sum_{i=0}^{d-1} a_ig_i(X_1,X_2)X_0^{i},\]
where, for $1\le i\le d$, $a_i$ is in $F^{\times}$ and
$g_i(X_1,X_2)$ is a form of degree $d-i$. For all $i$, write
\[g_i(X_1,X_2):=\sum_{j=0}^{d-i} b_{i,j} X_1^jX_2^{d-i-j}.\]
Assume also that the $a_i$ are chosen in such a way that $b_{i,j}=1$
for some $j$ with $j$ odd if possible.   Note that this is possible
when $g_i(X_1,X_2)\ne g_i(-X_1,X_2)$.    Since $R$ is not an
automorphism of $X$, $D=MRM^{-1}$ is not an automorphism of $X_1$,
so this is possible for at least one $g_i$.  Also note that
$a_ig_i(X_0,X_1)=0$ if and only if $a_i=0$.

Recall that
\[S_{m,1}=\begin{bmatrix}
\zeta_m&0&0\\
0&1&0\\
0&0&1\\
\end{bmatrix}\] where $\zeta_m$ is a primitive $(3^l)^{th}$ root of unity with
$l>0$ and that $S_{m,1}^3$ is in $\G$.   Since
$MS_{m,1}^3M^{-1}=S_{m,1}^3$, the element $S_{m,1}^3$ gives an
automorphism of $Y$.  Also, since $S_{m,1}$ does not give an
automorphism of $X$, $MS_{m,1}M^{-1}=S_{m,1}$ does not give an
automorphism of $Y$.

 Fix $\sigma\in\Gamma$ and suppose an
 isomorphism $Y\to\el{\sigma}{Y}$ is
 given by  $S_{m,1}^sD^r$, with
$0\le s\le 2$ and $0\le r\le 1$.   We have
\[g^{\sigma}(X_0,X_1,X_2)=\lambda_{\sigma}g(\zeta_m^{-s} X_0,(-1)^rX_1,X_2)\] for some
$\lambda_{\sigma}\in F^{\times}$.  Since
\[g(\zeta_m^{-s} X_0,(-1)^rX_1,X_2)
=\zeta_m^{-sd}(X_0^d+\sum_{i=0}^{d-1}
a_i\zeta_m^{s(d-i)}g_i((-1)^rX_1,X_2)X_0^{i}),\] we must have
$\lambda_{\sigma}=\zeta_m^{-sd}$.  So
\[g^{\sigma}(X_0,X_1,X_2)=X_0^d+\sum_{i=0}^{d-1}
(-1)^{i_t}\zeta_m^{s(d-i)}a_ig_i^{\sigma}(X_1,X_2)X_0^{i}\] where
\[i_t=\left\{
  \begin{array}{ll}
    0, & \hbox{if $g_i((-1)^rX_1,X_2)=g_i(X_1,X_2)$} \\
    1, & \hbox{if $g_i((-1)^rX_1,X_2)\ne g_i(X_1,X_2)$.}
  \end{array}\right.\]
  Then we must have $\sigma(a_i)=(-1)^{i_t}\zeta_m^{s(d-i)}a_i$.  Note
  that since $S_{m,1}^3$ gives an automorphism of $X$,
  $(\zeta_m^3)^{(d-i)}=1$ for all $i$ with $a_i\ne 0$.
  Since $S_{m,1}$ does not give an automorphism of $Y$, we can choose $i_1$ so that  $a_{i_1}\ne 0$  and so that
  $\zeta_m^{d-i_1}$ is a primitive cube root of unity.  So we can write $d-i_1=3^{l-1}u$ where $u\in\Z_{>0}$
  and $(u,3)=1$.  Choose $\alpha\in F$ so that
  $\alpha^{3^{l-1}}=a_{i_1}^{4u}$.  Note that our choice of $\alpha$
  is independent of $\sigma$.
   Then \[\sigma(\alpha^{3^{l-1}})=\sigma(a_{i_1}^{4u})=(\zeta_m^{s(d-i_1)})^{4u}\alpha^{3^{l-1}}=(\zeta_m^{s3^{l-1}})^{4u^2}\alpha^{3^{l-1}}=\zeta_m^{s3^{l-1}}\alpha^{3^{l-1}}.\]
   Since
   $\left(\frac{\sigma(\alpha)}{\alpha}\right)^{3^{l-1}}=(\zeta_m^s)^{3^{l-1}}$,
   we must have $\frac{\sigma(\alpha)}{\alpha}=\zeta^s_m(\zeta_m^3)^q$ for some $q\in\Z$.
      Choose $i_2$
   so that  $a_{i_2}\ne 0$ and  so that $g_{i_2}(-X_1,X_2)\ne g_{i_2}(X_1,X_2)$.  Let $\beta=a_{i_2}^3$.  Note that our
   choice of $\beta$ is independent of $\sigma$.  Then
   $\frac{\sigma(\beta)}{\beta}=(-1)^r$.  Thus there exists $q\in\Z$ so
   that
   \[S_{m,1}^sD^r(S_{m,1}^3)^q=\begin{bmatrix}
  \frac{\sigma(\alpha)}{\alpha}&0&0\\
  0&\frac{\sigma(\beta)}{(\beta)}&0\\
  0&0&1\\
  \end{bmatrix}.\]  Note that $(S_{m,1}^3)^q\in\Aut(Y)$.

  For all $\tau\in\Gamma$ define an isomorphism $X\to\el{\tau}{X}$
  by
  \[\begin{bmatrix}
  \frac{\tau(\alpha)}{\alpha}&0&0\\
  0&\frac{\tau(\beta)}{(\beta)}&0\\
  0&0&1\\
  \end{bmatrix}.\]
 The family
$\{\varphi_{\tau}\}_{\tau\in\Gamma}$ satisfies Weil's cocycle
condition of Theorem~\ref{weil}.  So by Theorem~\ref{weil}, $Y$ is
definable over $K$.

\end{enumerate}
\end{proof}

\begin{prop}\label{typed}  Let $K$ be a perfect field of characteristic
$p\ne 2$ and let $F$ be an algebraic closure of $K$.  Let $X$ be a
smooth plane curve of degree $d>3$ defined over $F$.  Suppose that
$\Aut(X)$ is a group of the type $\mf{D}$ given in
Lemma~\ref{gps}(c) and that $\Aut(X)$ is not given by $\G_{18}$.
      Then $X$ can be defined over its field of moduli relative
  to the extension $F/K$.
\end{prop}

\begin{proof}By~Proposition~\ref{closed}, we
may assume that $K$ is the field of moduli of $X$ relative to the
extension $F/K$. Let $\Gamma=\Gal(F/K)$.  Let $\G:=\Aut(X)$.  We
follow the notation used in Lemma~\ref{normalizergps}.  There are
two cases.
\begin{enumerate}
    \item[ii] Suppose $S_{m,1}\in\mf{H}$ and $\mf{H}\ne\la S_{3,2}\ra$.
    Since our argument is so similar the argument shown
in the proof of Proposition~\ref{typec}~ii, we omit some details.
Suppose that $[1\colon 0\colon 0]$ is in $X(F)$.  By
  Lemma~\ref{conjcurve} and Lemma~\ref{normalizergps}, for all
  $\sigma\in\Gamma$, any isomorphism
$X\to\el{\sigma}{X}$ is given by an element of $N(\G)\subseteq\la\G,
S_{3m,2}\ra$.
 Since $\G$ is normal in $N(\G)$, every element of $N(\G)$ can be written as
   an element of $\G$ times an element of $\la S_{3m,2}\ra$.
Any element of $\la S_{3m,2}\ra$ fixes $[1\colon 0\colon 0]$.
 For all $\sigma\in\Gamma$ we have $\Aut(\el{\sigma}{X}) = \G^{\sigma} = \G$.  Thus for all $\sigma\in\Gamma$, any isomorphism $X\to\el{\sigma}{X}$
maps the point $[1\colon 0\colon 0]$  into the
$\Aut(\el{\sigma}{X})$-orbit of  $[1\colon 0\colon 0]$.  So if the
point $[1\colon 0\colon 0]$ is in $X(F)$, it will map to a rational
point of the canonical $K$-model of $X/\Aut(X)$ given in
Theorem~\ref{bk}.  So in this case, by Corollary~\ref{pqpt}, $X$ can
be defined over $K$.

{}From now on assume that $[1\colon 0\colon 0]$ is not in $X(F)$.
Let $S:=TRS_{3m,2}(TR)^{-1}$.  So
$N(\G)\subseteq\la\G,S_{3m,2}\ra=\la\G,S\ra$.  We have
$\Z/3\Z\cong\la\G, S\ra/\la S\ra$, so for all $\sigma\in\Gamma$ an
isomorphism $X\to\el{\sigma}{X}$ is given by an element of the form
$S^s$ with $0\le s\le 2$.  It is easy to see that an element
$\sigma\in\Gamma$ uniquely determines $s\in\{0,1,2\}$.  Fix an
equation
\[f(X_0,X_1,X_2)=0\] for $X$.  Since $[1\colon 0\colon 0]$ is not in $X(F)$,
after multiplying $f(X_0,X_1,X_2)$ by an element of $F^{\times}$ we
may assume that \[f(X_0,X_1,X_2)=X_0^d+\sum_{i=0}^{d-1}
f_i(X_1,X_2)X_0^{i},\] where $f_i(X_1,X_2)$ is a form of degree
$d-i$.  For all $i$, write
\[f_i(X_1,X_2):=\sum_{j=0}^{d-i} b_{i,j} X_1^jX_2^{d-i-j}.\]

Since the elements $A_1$ and $A_2$ in the proof of
Lemma~\ref{typec}~ii are automorphisms of $X$, for all $i,j$ with
$b_{i,j}\ne 0$, we must have $j\equiv 0\pmod{m}$ and
$d-i-j\equiv0\pmod{m}$.  Since $S$ is not an automorphism of $X$,
there exists $i_1$ and $j_1$ such that $b_{i_1,j_1}\ne 0$ and such
that $2j_1+i_1-d\not\equiv 0\pmod{3m}$. Let $a=b_{i_1,j_1}$. Write
$2j_1+i_1-d=mu$ where $u\in\Z$ and $(u,3)=1$.  Choose $\alpha\in F$
such that $\alpha^m=a^u$.
 For each $\sigma\in\Gamma$, define
$\varphi_{\sigma}$ by
\[\begin{bmatrix}
1&0&0\\
0&\frac{\alpha}{\sigma(\alpha)}&0\\
0&0&\frac{\sigma(\alpha)}{\alpha}\\
\end{bmatrix}.\]  Then $\varphi_{\sigma}$ gives an isomorphism
$X\to\el{\sigma}{X}$ and it can be verified that family
$\{\varphi_{\sigma}\}_{\sigma\in\Gamma}$ satisfies Weil's cocycle
condition of Theorem~\ref{weil}.  So by Theorem~\ref{weil}, $X$ is
definable over $K$.

\item[iii] Suppose $S_{m,1}\not\in\mf{H}$ and $\mf{H}\ne\la S_{3,2}\ra$.
 Then by Lemma~\ref{isnorm}, any isomorphism
$X\to\el{\sigma}{X}$ is given by an element of $N(\G)$. By
Lemma~\ref{normalizergps}, $N(\G)\subseteq\la\G, S_{m,1}\ra$.  By
Lemma~\ref{diagisom}, $X$ is definable over $K$.
\end{enumerate}
\end{proof}

\begin{prop}\label{gpsl}  Let $K$ be a perfect field of characteristic
$p> 2$ and let $F$ be an algebraic closure of $K$.  Let $X$ be a
smooth plane curve of degree $d>3$ defined over $F$.  Suppose that
$\Aut(X)$ is given by the group $\G_{\PSL_2(q)}$ of
Lemma~\ref{gps}(n).  Then $X$ can be defined over its field of
moduli relative to the extension $F/K$.
\end{prop}
\begin{proof} By~Proposition~\ref{closed}, we may assume that $K$ is the field of
moduli of $X$ relative to the extension $F/K$.  Let
$\Gamma:=\Gal(F/K)$.  By Lemmas~\ref{isnorm}
and~\ref{normalizergps}, for all $\sigma\in\Gamma$ any isomorphism
$X\to\el{\sigma}{X}$ is given by an element of
\[N(\G)=\G_{\PGL_2(q)}=\la \G_{\PSL_2(q)}, M\ra\]
  where
  \[M:=\begin{bmatrix}
  \alpha^2&0&0\\
  0&\alpha&0\\
  0&0&1\\
  \end{bmatrix}\] for any $\alpha\in\F_q$ that is not a square.   If $-1$ is not a square then we may
  take $\alpha:=-1$.  In this case $X$ can be defined over $K$ by Lemma~\ref{diagisom}, or by
  Lemma~\ref{split}.  So assume that $-1$ is a square in $\F_q$.  Since
  $[\G_{\PGL_2(q)}\colon\G_{\PSL_2(q)}]=2$, and since if $\alpha$ is not a
  square in $\F_q$ and if $\sigma$ is  in $\Gamma$ then
  $\sigma(\alpha)$ is also not a square, the subgroup
  \[\Lambda:=\{\sigma\in\Gamma\colon X=\el{\sigma}{X}\}\] is a
  normal subgroup of $\Gamma$ and $[\Gamma\colon\Lambda]\le 2$.

  If $[\Gamma\colon\Lambda]=1$ there is nothing to prove so assume
  that $[\Gamma\colon\Lambda]=2$.  Let $L:=F^{\Lambda}$ be the
  subfield of $F$ fixed by $\Lambda$.  Then there exists $c\in L-K$
  with $c^2\in K$ and we have $L=K(c)$.  For any $\sigma\in\Gamma$
  we have
  \[\sigma(c)=\left\{
            \begin{array}{ll}
              c, & \hbox{if $\sigma|_L=Id$} \\
              -c, & \hbox{if $\sigma|_L\ne Id$.}
            \end{array}
          \right.\]
          Let $\alpha$ be in $\F_q$ and suppose that $\alpha$ is not a
          square in $\F_q$.  Let $n$ be the smallest integer such
          that $\alpha^n=1$.  Replacing $\alpha$ with $-\alpha$ if necessary, we
may assume that $n$ is even. Choose $\gamma\in F$ with
$\gamma^{n/2}=c$.  Let $\sigma$ be in $\Gamma$.  Then
$\sigma(\gamma)=\alpha^k\gamma$ for some $k$.  If
$X=\el{\sigma}{X}$, then $\sigma(c)=c$ and so
$(\alpha^k)^{\frac{n}{2}}=1$.  Since $n$ is the smallest integer
such that $\alpha^n=1$, we must have $\frac{kn}{2}\equiv 0\pmod{n}$.
It follows that $2\mid k$ and so $\alpha^k$ is a square in $\F_q$
and so\[ \begin{bmatrix}
\left(\frac{\sigma(\gamma)}{\gamma}\right)^2 &0&0\\
0&\frac{\sigma(\gamma)}{\gamma} &0\\
0&0&1\\
\end{bmatrix}\]  gives an automorphism of $X$.  Similarly, if
$X\ne\el{\sigma}{X}$, then $\sigma(c)=-c$ and so
$(\alpha^k)^{\frac{n}{2}}=-1$.  Since $n$ is the smallest integer
such that $\alpha^n=1$, we must have $2\nmid k$ and so $\alpha^k$ is
not a square in $\F_q$ and \[ \begin{bmatrix}
\left(\frac{\sigma(\gamma)}{\gamma}\right)^2 &0&0\\
0&\frac{\sigma(\gamma)}{\gamma} &0\\
0&0&1\\
\end{bmatrix}\]  gives an isomorphism $X\to\el{\sigma}{X}$.

Then for each $\sigma\in\Gamma$ the map $\varphi_{\sigma}$
 given by
\[ \begin{bmatrix}
\left(\frac{\sigma(\gamma)}{\gamma}\right)^2 &0&0\\
0&\frac{\sigma(\gamma)}{\gamma} &0\\
0&0&1\\
\end{bmatrix}\] gives an isomorphism $X\to\el{\sigma}{X}$.  The family
$\{\varphi_{\sigma}\}_{\sigma\in\Gamma}$ satisfies Weil's cocycle
condition of Theorem~\ref{weil}.  So by Theorem~\ref{weil}, $X$ is
definable over $K$.
\end{proof}

\begin{prop}\label{ga6}  Let $K$ be a perfect field of characteristic
$5$ and let $F$ be an algebraic closure of $K$.  Let $X$ be a smooth
plane curve of degree $d>3$ defined over $F$.  Suppose that
$\Aut(X)$ is given by the group $\G_{A_6}$ of Lemma~\ref{gps}(p).
Then $X$ can be defined over its field of moduli relative to the
extension $F/K$.
\end{prop}
\begin{proof}
By~Proposition~\ref{closed}, we may assume that $K$ is the field of
moduli of $X$ relative to the extension $F/K$.  Let
$\Gamma:=\Gal(F/K)$.  By Lemmas~\ref{isnorm} and~\ref{gps}(p), for
all $\sigma\in\Gamma$ any isomorphism $X\to\el{\sigma}{X}$ is given
by an element of $N(\G_{A_6})=\la\G_{A_6}, U\ra$
  where
  \[U:=\begin{bmatrix}
  1&0&0\\
  0&\alpha&\alpha\\
  0&\alpha&-\alpha\\
  \end{bmatrix}\] and where $\alpha^2=2$.  By Lemma~\ref{gps},
  $[N(\G_{A_6})\colon\G_{A_6}]=2$.

   Define
  \[\Lambda:=\{\sigma\in\Gamma\colon X=\el{\sigma}{X}\}.\]  Since
  $U^2$ gives an automorphism of $X$ we must have
  $[\Gamma\colon\Lambda]\le 2$.  So $\Lambda$ must be a normal
  subgroup of $\Gamma$.  If $[\Gamma\colon\Lambda]=1$, then clearly
  $X$ is definable over its field of moduli.  So assume
  $[\Gamma\colon\Lambda]=2$.  Let $L:=F^{\Lambda}$ be the subfield
  of $F$ fixed by $\Lambda$.  Write $L=K(\sqrt{c})$ with $c\in K^{\times}-(K^{\times})^2$.
    There are three cases.
  \begin{enumerate}
    \item [Case] I: $\alpha\in K$.  Choose
    $\gamma\in F^{\times}$ with $\gamma^4=c$.  Then $K(\gamma)$ is the splitting field of the
    polynomial $x^4-c=(x-\gamma)(x-2\gamma)(x-4\gamma)(x-3\gamma)$ over $K$.  Suppose $\sigma\in\Gamma$ and suppose
    that $\sigma|_{K(\gamma)}:=\ol{\sigma}$ generates $\Gal(K(\gamma)/K$), say
    \[\ol{\sigma}(\gamma)=2\gamma.\]
    For each $\tau\in\Gamma$ the map $\varphi_{\tau}$
 given by $U^i$ if $\tau|_{K(\gamma)}=\ol{\sigma}^i$, gives an
 isomorphism $X\to\el{\tau}{X}$.  Since $U^i$ is defined over $K$ for all
 $i$, and since $U$ has order $4$,
 it is easily verified that the family
$\{\varphi_{\sigma}\}_{\sigma\in\Gamma}$ satisfies Weil's cocycle
condition of Theorem~\ref{weil}.  So by Theorem~\ref{weil}, $X$ is
definable over $K$.
    \item [Case] II: $L=K(\alpha)$.  For each $\sigma\in\Gamma$ the map
$\varphi_{\sigma}$
 given by
\[\left\{
  \begin{array}{ll}
    U, & \hbox{if $\sigma|_{L}\ne Id$} \\
    Id, & \hbox{if $\sigma|_{L}= Id$}
  \end{array}
\right.\]

  gives an
 isomorphism $X\to\el{\sigma}{X}$.  Note that
\[U^{\sigma}=\left\{
  \begin{array}{ll}
    U^{-1}, & \hbox{if $\sigma|_{L}\ne Id$} \\
    U, & \hbox{if $\sigma|_{L}= Id$}
  \end{array}
\right.\]  {}From this, it is clear that the family
$\{\varphi_{\sigma}\}_{\sigma\in\Gamma}$ satisfies Weil's cocycle
condition of Theorem~\ref{weil}.  So by Theorem~\ref{weil}, $X$ is
definable over $K$.
    \item [Case] III: $\alpha\not\in L$.  Choose $\gamma\in F^{\times}$
    with $\gamma^4=c$.  Then $K(\gamma,\alpha)$ is the splitting
    field of the polynomials
    \[x^4-c\mbox{ and } x^2-2\]
    over $K$.  Kummer theory shows that
    $\Gal(K(\gamma,\alpha)/K)=\la\ol{\sigma},\ol{\tau}\ra$ where
\[\ol{\sigma}(\gamma)=2\gamma,\; \ol{\sigma}(\alpha)=\alpha,\]
      and
    \[\ol{\tau}(\gamma)=\gamma,\; \ol{\tau}(\alpha)=-\alpha.\]  The
    elements $\ol{\tau}$ and $\ol{\sigma}$ satisfy the relations
    \[\ol{\sigma}^4=\ol{\tau}^2=Id\] and
    \[\ol{\tau}\;\ol{\sigma}=\ol{\sigma}\;\ol{\tau}.\]

     Following the notation of Lemma~\ref{gps}(p), let \[Q:=VTV=\begin{bmatrix}
  1&0&0\\
  0&-1&0\\
  0&0&1\\
  \end{bmatrix}\in\G_{A_6}.\]
  Note that $U^{-1}=U^3=QUQ^{-1}$, and that $Q^2=Id$.  Note also
  that if $\omega$ is in $\Gamma$ and $\omega(\alpha)=-\alpha$ then
  $U^{\omega}=U^{-1}$.

  For all $\omega\in\Gamma$ define $\varphi_{\omega}\colon
  X\to\el{\omega}{X}$, by
  \[(U^i)^{\ol{\tau}^j}Q^j\]
  if \[\omega|_{K(\gamma,\alpha)}=\ol{\sigma}^i\;\ol{\tau}^j,\]
  where
  $0\le i\le 3$ and $0\le j\le 1$.

  Let $\omega_1$ and $\omega_2$ be in $\Gamma$ and suppose that
  \[\omega_1|_{K(\gamma,\alpha)}=\ol{\sigma}^i\;\ol{\tau}^j,\] and
  \[\omega_2|_{K(\gamma,\alpha)}=\ol{\sigma}^k\;\ol{\tau}^l\] where
$0\le i, k\le 3$ and $0\le j,l\le 1$.  Then
\[\varphi_{\omega_1}^{\omega_2}\varphi_{\omega_2}\] is given by
\[((U^i)^{\ol{\tau}^j}Q^j )^{\omega_2}(U^k)^{\ol{\tau}^l}Q^l \]
\[= [(U^i)^{\ol{\tau}^j}Q^j U^kQ^l]^{\ol{\tau}^l} \]
\[=\left\{
  \begin{array}{ll}
    \left[(U^i)^{\ol{\tau}^j}U^kQ^{j+l}\right]^{\ol{\tau}^l}, & \hbox{if $j=0$} \\
    \left[(U^i)^{\ol{\tau}^j}U^{-k}Q^{j+l}\right]^{\ol{\tau}^l}, & \hbox{if $j=1$}
  \end{array}
\right.\]
\[= [(U^{i+k})^{\ol{\tau}^j}Q^{j+l}]^{\ol{\tau}^l} \]
\[= (U^{i+k})^{\ol{\tau}^{j+l}}Q^{j+l} \] which equals the
isomorphism $\varphi_{\omega_2\omega_1}$.

 So the family
$\{\varphi_{\omega}\}_{\omega\in\Gamma}$ satisfies Weil's cocycle
condition of Theorem~\ref{weil}.  So by Theorem~\ref{weil}, $X$ is
definable over $K$.
\end{enumerate}
 \end{proof}

\begin{thm}\label{mainplane}  Let $K$ be a perfect field of characteristic
 $p$ where $p=0$ or
$p>2$ and let $F$ be an algebraic closure of $K$.  Let $X$ be a
smooth plane curve of degree $d>3$ over $F$.  Suppose that $\Aut(X)$
is not $\PGL_3(F)$-conjugate to a diagonal subgroup of $\PGL_3(F)$
or to one of the groups $\G_{18}$, $\G_{36}$, or a group
 given by Lemma~\ref{gps}(s).
Then $X$ can be defined over its field of moduli relative to the
extension $F/K$.
\end{thm}
\begin{proof}  By~Proposition~\ref{closed}, we may assume that $K$ is the field of
moduli of $X$ relative to the extension $F/K$.  Let $\G:=\Aut(X)$.
 By Lemma~\ref{gps}, we may assume that $\G$ is one of the groups in
 Lemma~\ref{gps} that is not a diagonal group, $\G_{18}$,$\G_{36}$,
   or a group
 given by Lemma~\ref{gps}(s).  Let $\Gamma:=\Gal(F/K)$.  We follow the notation of
 Lemma~\ref{gps}.  There are seventeen cases.
 \begin{enumerate}
    \item[(a)] $\G$ is an intransitive group whose image $\ol{\G}$
    in $\PGL_2(F)$ is equal to one of the groups in
    Lemma~\ref{subgroups} Case~I~(b)-(e).  Then $X$ is definable
    over $K$ by Proposition~\ref{intransitive}.
    \item[(b)] $\G$ is a group of type~$\mf{C}$.  Then $X$ is definable
    over $K$ by Proposition~\ref{typec}.
    \item[(c)] $\G\not=\G_{18}$ and $\G$ is a group of type~$\mf{D}$.  Then $X$ is definable
    over $K$ by Proposition~\ref{typed}.
    \item[(e)] $\G=\G_{72}$.  It is easily verified that
    $\G_{72}^{\sigma}=\G_{72}$ for all $\sigma\in\Gamma$.  By
    Lemma~\ref{normalizergps}, $N(\G_{72})=\G_{72}\rtimes\la U\ra$.
    Since $U^{\sigma}=U$ for all $\sigma\in\Gamma$, by
    Proposition~\ref{split}, $X$ is definable over $K$.
    \item[(f)] $\G=\G_{216}$.  It is easily verified that
    $\G_{216}^{\sigma}=\G_{216}$ for all $\sigma\in\Gamma$.  By
    Lemma~\ref{normalizergps}, $N(\G_{216})=\G_{216}$.
    By
    Proposition~\ref{split}, $X$ is definable over $K$.
    \item[(g)] $\G=\G_{60}$.  It is easily verified that
    $\G_{60}^{\sigma}=\G_{60}$ for all $\sigma\in\Gamma$.   By
    Lemma~\ref{normalizergps}, $N(\G_{60})=\G_{60}$.
    By
    Proposition~\ref{split}, $X$ is definable over $K$.
     \item[(h)] $\G=\G_{360}$.  By Lemma~\ref{conjcurve}, any
     isomorphism $X\to\el{\sigma}{X}$, with $\sigma\in\Gamma$ is
     given by an element of the form
     \[\begin{bmatrix}
     \lambda&0&0\\
     0&1&0\\
     0&0&1\\
     \end{bmatrix},\] where $\lambda\in\{1,\lambda_2\}$.  So by
     Lemma~\ref{diagisom}, $X$ is definable over $K$.
    \item[(i)] $\G=\G_{168}$.  It is easily verified that
    $\G_{168}^{\sigma}=\G_{168}$ for all $\sigma\in\Gamma$.   By
    Lemma~\ref{normalizergps}, $N(\G_{168})=\G_{168}$.
    By
    Proposition~\ref{split}, $X$ is definable over $K$.
    \item[(j)] $\G=\PSL_3(\F_q)$.  It is easily verified that
    $(\PSL_3(\F_q))^{\sigma}=\PSL_3(\F_q)$ for all $\sigma\in\Gamma$.   By
    Lemma~\ref{normalizergps}, $N(\PSL_3(\F_q))=\la\PSL_3(\F_q),
    M\ra$, where
    \[M:=\begin{bmatrix}
    \alpha&0&0\\
    0&1&0\\
    0&0&1\\
    \end{bmatrix}\] with $\alpha=1$ if $q\not\equiv 1\pmod 3$ and with
    $\alpha\in(\F_q)^3-\F_q$ if $q\equiv 1 \pmod 3$.  By Lemma~\ref{isnorm},
    if $\sigma\in\Gamma$, any isomorphism
    $X\to\el{\sigma}{X}$ is given by an element of $N(\PSL_3(\F_q))$.  So by
     Lemma~\ref{diagisom}, $X$ is definable over $K$.
     \item[(k)]$\G=\PGL_3(\F_q)$.  It is easily verified that
    $(\PGL_3(\F_q))^{\sigma}=\PGL_3(\F_q)$ for all $\sigma\in\Gamma$.   By
    Lemma~\ref{normalizergps}, $N(\PGL_3(\F_q))=\PGL_3(\F_q)$.   So by
     Lemma~\ref{split}, $X$ is definable over $K$.
      \item[(l)] $\G=\PSU_3(\F_q)$.  It is easily verified that
    $(\PSU_3(\F_q))^{\sigma}=\PSU_3(\F_q)$ for all $\sigma\in\Gamma$.   By
    Lemma~\ref{normalizergps}, $N(\PSU_3(\F_q))=\la\PSU_3(\F_q),
    M\ra$, where
    \[M:=\begin{bmatrix}
    \alpha&0&0\\
    0&1&0\\
    0&0&1\\
    \end{bmatrix}\] with $\alpha=1$ if $q\equiv 0\pmod 3$, with
    $\alpha\in(\F_q)^3-\F_q$ if $q\equiv 1 \pmod 3$, and with
    $\alpha\in(\F_{q^2})^{q-1}-\F_{q^2}^{3(q-1)}$ if $q\equiv 2 \pmod 3$.  By Lemma~\ref{isnorm},
    if $\sigma\in\Gamma$, any isomorphism
    $X\to\el{\sigma}{X}$ is given by an element of $N(\PSU_3(\F_q))$.  So by
     Lemma~\ref{diagisom}, $X$ is definable over $K$.
     \item[(m)]$\G=\PGU_3(\F_q)$.  It is easily verified that
    $(\PGU_3(\F_q))^{\sigma}=\PGU_3(\F_q)$ for all $\sigma\in\Gamma$.   By
    Lemma~\ref{normalizergps}, $N(\PGU_3(\F_q))=\PGU_3(\F_q)$.   So by
     Lemma~\ref{split}, $X$ is definable over $K$.
     \item[(n)] $\G=\G_{\PSL_2(q)}$.  By
 Proposition~\ref{gpsl}, $X$ can be defined over $K$.
  \item[(o)]$\G=\G_{\PGL_2(q)}$.  It is easily verified that
    $(\G_{\PGL_2(q)})^{\sigma}=\G_{\PGL_2(q)}$ for all $\sigma\in\Gamma$.   By
    Lemma~\ref{normalizergps}, $N(\G_{\PGL_2(q)})=\G_{\PGL_2(q)}$.   So by
     Lemma~\ref{split}, $X$ is definable over $K$.
     \item[(p)] $\G=\G_{A_6}$, $\G=N(\G_{A_6})$, or $\G=\G_{A_7}$.
     In each case, it is easily verified that $\G^{\sigma}=\G$ for
     all $\sigma\in\Gamma$.  If $\G=\G_{A_6}$, then by
 Proposition~\ref{ga6}, $X$ can be defined over $K$.

 If $\G\ne\G_{A_6}$, then by  Lemma~\ref{normalizergps}, $N(\G)=\G$.
   So by
     Lemma~\ref{split}, $X$ is definable over $K$.
     \item[(q)] $\G=\G_{60}^3$, $\G=\G_{168}^3$, $\G$ is a group of
     type $\mf{C}^3$ or $\G$ is a group of type $\mf{D}^3$.  In each case, it is easily verified that $\G^{\sigma}=\G$ for
     all $\sigma\in\Gamma$.  By Lemma~\ref{normalizergps}, if $\G$ is a group of type $\mf{C}^3$
     then $N(\G)=\G\rtimes\mf{K}$ where $\mf{K}\le \la R\ra$, and  in
     all other cases $N(\G)=\G$.  Note that $\mf{K}^{\sigma}=\mf{K}$
     for all $\sigma\in\Gamma$.   So by
     Lemma~\ref{split}, $X$ is definable over $K$.
     \item[(r)] $\G$ is a group given in Lemma~\ref{gps}(r).  By
     proposition~\ref{intransitive}, $X$ is definable over $K$.

 \end{enumerate}
\end{proof}
\begin{thm}\label{mainplane2}  Let $K$ be a field of characteristic
 not equal to $2$ and let $F$ be an algebraic closure of $K$.  Let $X$ be a
smooth plane curve of degree $d>3$ over $K$.  Suppose that $\Aut(X)$
is not $\PGL_3(F)$-conjugate to a diagonal subgroup of $\PGL_3(F)$
or to one of the groups $\G_{36}$, $\G_{18}$, or a group
 given by Lemma~\ref{gps}(s).
Then $X$ can be defined over its field of moduli.
\end{thm}
\begin{proof}This follows from Theorem~\ref{mainplane} and
Theorem~\ref{relmod}.
\end{proof}

%% file: ch7.tex
Let $K$ be a field of characteristic not equal to $2$, let $F$ be an
algebraic closure of $K$, and let $X$ be a smooth plane curve over
$K$. By Theorem~\ref{mainplane2}, $X$ can be defined over its field
of moduli if $\Aut(X)$ is not $\PGL_3(F)$-conjugate to an
intransitive group whose image in $\PGL_2(F)$ is cyclic, not
$\PGL_3(F)$-conjugate to $\G_{18}$, not $\PGL_3(F)$-conjugate to
$\G_{36}$, and not $\PGL_3(F)$-conjugate to a group of the type
listed in Lemma~\ref{gps}(s).  We now construct smooth plane curves
not definable over their field of moduli with automorphism groups
conjugate to all but the last type of group.

\section{Plane curves with diagonal automorphism groups}
We now construct smooth planes curves with diagonal automorphism
groups not definable over their fields of moduli.  These examples
are easily obtained from adjusting the examples of hyperelliptic
curves not definable over their fields of moduli.  We will construct
examples using the polynomial $f(x)$ given in Chapter~5.  Note that
we could just as easily work with the polynomial $g(x)$ given in
Chapter~5.

Suppose $n$, $r\in\Z_{>0}$.  Assume that $2nr>5$.   Assume also that
if $n$ is odd, then $r$ is also odd.  Let $z^c$ be the complex
conjugate of $z$ for any $z\in\C$.  Suppose $a_1,\ldots,a_r \in \C$.
Consider the form $f(X_0, X_1)\in\C[X_0, X_1]$ given by
\[f(X_0,X_1):=\prod_{i=1}^r(X_0^n-a_iX_1^n)(X_0^n+a_i^cX_1^n).\]
 Assume that $f(X_0,X_1)$ has no repeated zeros and
that it is not a form in $\R[X_0,X_1]$. Assume that the map
$[\alpha\colon\beta]\mapsto [\beta\colon\alpha]$ does not map the
zero set of $f$ into itself. For any
  root of unity $\zeta\ne 1$, assume
that  \[\{a_i,-1/a_i^c\}_{i=1}^r\ne\{\zeta
a_i,-\zeta/a_i^c\}_{i=1}^r.\] Lastly, if $n=3$ assume that the map
\[[\alpha\colon\beta]\mapsto[-\alpha+(\sqrt{3}+1)\beta\colon \alpha(\sqrt{3}+1)+\beta]\] does not map the
zero set of $f$ into itself.
 Define
 \[h(X_0,X_1,X_2):=X_2^{2nr}-f(X_0,X_1).\]

\begin{lem}\label{diagaut}Following the above notation, $h(X_0,X_1,X_2)=0$  gives the equation of a smooth
plane curve $X$ with $\Aut(X)$ given by a diagonal group generated
by $E$, $F$, and $H$ where
\[E:=\begin{bmatrix}
\zeta_{n}&0&0\\
0&1&0\\
0&0&1\\
\end{bmatrix},\;F:=\begin{bmatrix}
1&0&0\\
0&\zeta_{n}&0\\
0&0&1\\
\end{bmatrix},\;H:=\begin{bmatrix}
1&0&0\\
0&1&0\\
0&0&\zeta_{2nr}\\
\end{bmatrix},\] and where $\zeta_n$ and $\zeta_{2nr}$ are primitive
$n^{th}$ and $2nr^{th}$ roots of unity, respectively.
\end{lem}
\begin{proof}  We first show that $h(X_0,X_1,X_2)=0$ gives the equation of
a smooth plane curve. Since $h(X_0,0,X_2)=X_2^{2nr}-X_0^{2nr}$ has
$2nr$ distinct zeros, it is clear that no zero of $h(X_0,0,X_2)$ can
be a zero of $h_{X_0}(X_0,0,X_2)$.  We have
$h_{X_0}(X_0,1,X_2)=-f'(X_0,1)$ and
$h_{X_2}(X_0,1,X_2)=2nrX_2^{2nr-1}$. So $h_{X_2}(X_0,1,X_2)=0$ if
and only if $X_2=0$. Since $f(X_0,1)$ is square free, $f(X_0,1)$ and
$f'(X_0,1)$ have no common zeros.  It follows that
$h(X_0,X_1,X_2)=0$ gives the equation of a smooth plane curve $X$.

Since the degree of $h$ is greater than $3$, the genus of $X$ is
larger than $1$, so $\Aut(X)$ is finite.  Since $X$ is a smooth
plane curve of degree larger than $3$, by Theorem~\ref{planeisom}
$\Aut(X)$ is $\PGL_3(\C)$-conjugate to one of the group listed in
Lemma~\ref{gps} Case~I.  We now show that $\Aut(X)=\la E, F, H\ra$.
It is clear that $\la E, F, H\ra\subseteq\Aut(X)$.

  The element
$H$ has order equal to $2nr>5$. The only groups listed in
Lemma~\ref{gps} Case~I that have elements of even order larger than
$5$ are the groups of type~(a)-(c) or $\G_{216}$.  The orders of the
element of $\G_{216}$ are $1$, $2$, $3$, $4$, or $6$.  Using Magma,
one can verify that $\G_{216}$ has two conjugacy classes of elements
of order $6$. Representatives for each class are given by
\[M_1:=\begin{bmatrix}
\omega&0&0\\
0&0&1\\
0&1&0\\
\end{bmatrix}\quad\mbox{and}\quad M_2:=\begin{bmatrix}
\omega^2&0&0\\
0&0&1\\
0&1&0\\
\end{bmatrix}\] where $\omega$ is a primitive cube root of unity.  Any matrix in $\GL_3(\C)$
mapping to either $M_1$ or $M_2$ in $\PGL_3(\C)$ has three distinct
eigenvalues, so $M_1$ and $M_2$ are not $\PGL_3(\C)$-conjugate to
$H$. It follows that $\Aut(X)$ is not $\PGL_3(\C)$-conjugate to
$\G_{216}$.   So $\Aut(X)$ is $\PGL_3(\C)$-conjugate to a group of
type~(a)-(c) of Lemma~\ref{gps}.

   Recall that an intransitive element of
$\PGL_3(\C)$ is an element of the form
\[\begin{bmatrix}
a&b&0\\
c&d&0\\
0&0&1\\
\end{bmatrix}.\]  A computation shows that the normalizer of $\la H\ra$ in
 $\PGL_3(\C)$ is equal the subgroup of intransitive elements of $\PGL_3(\C)$.  Note
that $X/\la H^2\ra$ is one of the hyperelliptic curves given in
Lemma~\ref{hypce}, and that an intransitive element of $\PGL_3(\C)$
that gives an automorphism of $X$ induces an automorphism of the
hyperelliptic curve $X/\la H^2\ra$.  By Lemma~\ref{hypce}, we have
$\Aut(X/\la H^2\ra)\cong\Z/2\Z\times\Z/n\Z$.  Observe that $\la
E,F,H\ra/\la H^2\ra\cong \Z/2\Z\times\Z/n\Z$.  It follows that the
normalizer of $\la H\ra$ in $\Aut(X)$ is equal to $\la E, F, H\ra$.

  Suppose $M\in\PGL_3(\C)$ and that $\G':=M\Aut(X) M^{-1}$ is a group of
type~(a), (b), or (c) of Lemma~\ref{gps}. Write $H':=MHM^{-1}$.

First suppose that $\G'$ is a group of type~(b) or (c).  Then,
following the notation of Lemma~\ref{gps}, any element of $\G'$ can
be written as $DR^iT^j$, with $0\le i\le 2$ and $0\le j\le 2$ and
where $D$ is the image in $\PGL_3(\C)$ of a diagonal matrix. A
computation shows that for $1\le j\le 2$, an element of the form
$DT^j$ has order $3<2nr$. Another computation shows that for $0\le
j\le 2$, an element of the form $DRT^j$ is the image in $\PGL_3(\C)$
of a matrix with only $2$ distinct eigenvalues only if it has order
$2<2nr$. It follows that $H'$ must be the image in $\PGL_3(\C)$ of a
diagonal matrix.  Note that $TH'T^{-1}\in\G'$ is an element of order
$2nr$ that commutes with $H'$ and that $\la TH'T^{-1}\ra\cap\la
H'\ra=Id$. It follows that there exist an element $H_2$ of order
$2nr$ in the normalizer of $\la H\ra$ in $\Aut(X)$ such that $\la
H_2\ra\cap\la H\ra=Id$.  This is a contradiction since if $A\in\la
E, F, H\ra$ has order $2nr$, then $A^{nr}=H^{nr}\ne Id$.

Lastly, suppose that $\G'$ is a group of type~(a).  Then there
exists an intransitive element $M'\in\PGL_3(\C)$ such that
$H'':=M'H'(M')^{-1}$ is the image in $\PGL_3(\C)$ of a diagonal
matrix with only two distinct eigenvalues.  Since
$\G'':=M'\G'(M')^{-1}$ is an intransitive group, there exists a
natural map $\psi\colon \G''\to\PGL_2(\C)$.  It is easy to see that
either $H''$ is in the kernel of $\psi$ or $\psi(H'')$ has order
$2nr$. If $\psi(H'')$ has order $2nr$, then since $2nr>5$, by
Lemma~\ref{subgroups}, $\psi(\G'')$ is either a cyclic or dihedral
group.  In any case, $\la\psi(H'')\ra$ is a normal subgroup of
$\psi(\G'')$.  Since the kernel of $\psi$ is contained in the center
of $\G''$, it must be true that $\la H''\ra$ is a normal subgroup of
$\G''$.  So $\la H\ra$ is a normal subgroup of $\Aut(X)$.  It
follows from that $\Aut(X)=\la E, F, H\ra$.
\end{proof}

\begin{prop}\label{diag} Following the above notation,
let $X$ by the smooth plane curve over $\C$ given by
$h(X_0,X_1,X_2)=0$. Then the field of moduli of $X$ relative to the
extension $\C/\R$ is $\R$ and is not a field of definition for $X$.
\end{prop}
\begin{proof} We follow the notation of Lemma~\ref{diagaut}.   There exists
 an isomorphism $\mu\colon X\to\el{c}{X}$ given by
\[\begin{bmatrix}
0&1&0\\
\zeta_{2n}&0&0\\
0&0&\gamma\\
\end{bmatrix}\] where
\[\gamma^{2nr}=\prod_{i=1}^r-a_i^c/a_i\] and $\zeta_{2n}$ is a
primitive $2n^{th}$ root of unity.  In particular, the field of
moduli of $X$ relative to $\C/\R$ is $\R$. Suppose that $X$ is
definable over $\R$. Then by Theorem~\ref{weil}, there exists an
isomorphism $\mu'\colon X\to\el{c}{X}$ that satisfies Weil's cocyle
condition. Since $\zeta_n'\zeta_{2n}$ is a $2n^{th}$ root of unity
not equal to $1$ for any $n^{th}$ root of unity $\zeta_n'$, since
$(\zeta_{2nr}'\gamma)^{2nr}=\prod_{i=1}^r-a_i^c/a_i$ for any
$2nr^{th}$ root of unity $\zeta_{2nr}'$, and since $\Aut(X)=\la E,
F, H\ra$, we may assume that
\[\mu':=\begin{bmatrix}
0&1&0\\
\zeta'&0&0\\
0&0&\gamma'\\
\end{bmatrix}\]where $\zeta'$ is a $2n^{th}$ root of unity not equal to $1$ and where
  \[(\gamma')^{2nr}=\prod_{i=1}^r-a_i^c/a_i.\]    Since
$(\mu')^c\mu'=Id$, we have
\[Id=\begin{bmatrix}
0&1&0\\
(\zeta')^{-1}&0&0\\
0&0&(\gamma')^c\\
\end{bmatrix}\begin{bmatrix}
0&1&0\\
\zeta'&0&0\\
0&0&\gamma'\\
\end{bmatrix}=\begin{bmatrix}
1&0&0\\
0&(\zeta')^{-2}&0\\
0&0&(\zeta')^{-1}(\gamma')^c\gamma'\\
\end{bmatrix}.
\]  So $(\zeta')^{-1}(\gamma')^c\gamma'=1$.  Note that $|\gamma'|=1$,
 so $(\gamma')^c\gamma'=1$.  So we must have $(\zeta')^{-1}=1$.
This is a contradiction. Therefore, $X$ is not definable over $\R$.
\end{proof}
\begin{rem}An alternate and slightly longer proof of Proposition~\ref{diag} shows
that the definability of $X$ over $\R$ implies the definability of
the hyperelliptic curve $X/\la H^2\ra$ over $\R$, contradicting
Proposition~\ref{hypnotdef}.
\end{rem}
\section{Plane curves with automorphism groups given by $\G_{18}$}

We now construct plane curves with automorphism groups given by
$\G_{18}$ that are not definable over their fields of moduli.  All
other examples of curves not definable over their fields of moduli
in this thesis are of curves over $\C$ with field of moduli $\R$
relative to $\C/\R$ but not definable over $\R$. It is not too hard
to show that a plane curve $X$ over $\C$ with $\Aut(X)=\G_{18}$ is
always definable over its field of moduli relative to $\C/\R$.  So
our examples in this section are somewhat different from the others.

\begin{lem}\label{g18orbits}  Let $F$ be an algebraically closed
field of characteristic $0$ or of characteristic greater than $3$.
Suppose  $P\in\mathbb{P}^2(F)$.  Then the orbit $\mf{O}(P)$ of $P$
under the action of $\G_{18}$ has $18$ points unless $P$ is in the
orbit of
\[P_3^1:=[1\colon 0\colon 0],\] \[P^2_3:=[1\colon 1\colon 1],\]
\[P^3_3:=[1\colon 1\colon \omega],\]
\[P^4_3:=[1\colon 1\colon \omega^2],\]
 or
\[P_{9, \delta}:=[\delta\colon 1\colon 1],\]
where $\omega$ is a primitive cube root of unity and where
$\delta\in F-\{1,\omega,\omega^2\}$. We have $|\mf{O}(P_3^i)|=3$ for
$1\le i\le 4$, and $|\mf{O}(P_{9, \delta})|=9$.
\end{lem}
\begin{proof}  We follow the notation of Lemma~\ref{gps}.
Let $P$ be in $\mathbb{P}^2(F)$ and suppose that the orbit of $P$
under the action of $\G_{18}$ has less than $18$ points. It follows
that for some $g\in \G_{18}-\{Id\}$, $g(P)=P$.  If $h$ is in
$\G_{18}$, then $(h^{-1}gh)(h^{-1}(P))=h^{-1}(P)$.  So any element
in $\G_{18}$ which is conjugate to $g$ fixes a point in the orbit of
$P$ under the action of $\G_{18}$.  It is easily shown that any non
identity element of $\G_{18}$ is $\G_{18}$-conjugate to one of $S$,
$T$, $ST$, $ST^2$, or $R$.

A computation shows the following.
 If $S(P)=P$ then \[P\in\{[1\colon 0\colon 0], [0\colon 1\colon
0], [0\colon 0\colon 1]\}.\] If  $T(P)=P$ then
\[P\in\{[1\colon 1\colon 1], [1\colon \omega\colon \omega^2], [1\colon \omega^2\colon \omega]\}.\]
If  $ST(P)=P$ then
\[P\in\{[1\colon 1\colon \omega^2], [1\colon \omega\colon \omega], [1\colon \omega^2\colon 1]\}.\]
If  $ST^2(P)=P$ then
\[P\in\{[1\colon \omega\colon 1], [1\colon 1\colon \omega], [1\colon \omega\colon \omega^2]\}.\]
 If $R(P)=P$ then
$P=[\delta \colon \gamma\colon \gamma]$ where $\delta$ and $\gamma$
are not both zero.   Note that if $\gamma=0$ then $P=P^1_3$,  and if
$\gamma\ne 0$ and $\frac{\delta}{\gamma} =\omega^j$, with $0\le j\le
2$, then $P\in \mf{O}(P_3^i)$, where  $2\le i\le 4$.

A simple computation shows that each of these points is in the orbit
of one of the points listed in the lemma.  Another computation using
the stabilizers of each point proves the last statement of the
lemma.
\end{proof}

\begin{df}Let $K$ be a field.  A \emph{quaternion extension} of $K$
is a Galois extension $F$ of $K$ such that $\Gal(F/K)$ is isomorphic
to the quaternion group of order $8$.
\end{df}

The following lemma will be helpful in constructing examples of
smooth plane curves with automorphism groups isomorphic to $\G_{18}$
and not definable over their fields of moduli.
\begin{lem}\label{quaternion} Let $K$ be a field of level $2$: the element $-1$ is not a square in $K$ but
it is a sum of two squares in $K$.  Let $u, v\in K^{\times}-
(K^{\times})^2$ be such that $uv \not\in(K^{\times})^2$. Then
$K(\sqrt{u},\sqrt{v})$ is embeddable into a quaternion extension of
$K$ if and only if $-u$ is a norm from $K(\sqrt{-v})$ to $K$, that
is if $-u=x^2+vy^2$ for some $x,y\in K$.
\end{lem}
\begin{proof}This follows from Theorem~I.3.3 on page~166  and
 Proposition~I.1.6 on page~160 of~\cite{Jensen}.
\end{proof}

Throughout the rest of this section   let $K=\Q(\omega)$ where
$\omega$ is a primitive cube root of unity.
\begin{rem}Note that $K$ is of level $2$
since $(\omega^2)^2+\omega^2=-1$ and $\sqrt{-1}\not\in K$. So we can
use Lemma~\ref{quaternion} to see if $\Z/2\Z\times\Z/2\Z$ extensions
of $K$ are embeddable into quaternion extensions of $K$. For
example, let $L:=K(\sqrt{-13})$. It is easily shown that both $2$
and $-2$ are not norms from $L$ to $K$.  So by
Lemma~\ref{quaternion}, neither $K(\sqrt{-2},\sqrt{13})$ nor
$K(\sqrt{2},\sqrt{13})$ are not embeddable into quaternion
extensions of $K$.
\end{rem} We define the following:
 \begin{align*}
  \phi&:= X_0X_1X_2,\\
 \psi&:=X_0^3+X_1^3+X_2^3,\\
 \chi&:=X_0^3X_1^3 + X_1^3 X_2^3 + X_2^3 X_0^3.\\
\end{align*}\vspace*{-.76in}

Suppose $\alpha_1,\alpha_2,\alpha_3, u, v\in \mathbb{Q}^{\times}$
and suppose that $L:=K(\sqrt{u},\sqrt{v})$ is a
 $\Z/2\Z\times\Z/2\Z$ extension of $K$ that cannot be embedded
into a quaternion extension of $K$.  Let
\begin{align*}
  c_{\phi^2}&:=\alpha_1\omega\sqrt{u}  +\alpha_2\sqrt{v} + \alpha_3\omega^2\sqrt{uv},\\
 c_{\phi\psi}&:=\alpha_1\omega^2\sqrt{u} +\alpha_2\sqrt{v} + \alpha_3\omega\sqrt{uv}),\\
 c_{\psi^2}&:=-1/12+\alpha_1\sqrt{u} +\alpha_2\sqrt{v}  + \alpha_3\sqrt{uv},\\
\end{align*}\vspace*{-.76in}
\newline and let
 \begin{align*}
 f_{\sqrt{u},\sqrt{v}}(X_0,X_1,X_2) &:=c_{\psi^2}\psi^2-6c_{\phi\psi}\phi\psi-18c_{\phi^2}\phi^2+\chi\\
                                    &\,\,=c_{\psi^2}(X_0^6+X_1^6+X_2^6)-6c_{\phi\psi}(X_0^3+X_1^3+X_2^3)X_0X_1X_2\\
&\,\,- 18c_{\phi^2}(X_0X_1X_2)^2+ (2c_{\psi^2}+1)(X_0^3 X_1^3 +
X_1^3
X_2^3 + X_2^3 X_0^3).\\
 \end{align*}\vspace*{-.76in}

 Assume also that $f_{\sqrt{u},\sqrt{v}}(X_0,1,1)$ is squarefree.
(This is possible: for example if $\alpha_i=1$ for $1\le i\le 3$,
$u=2$, and $v=13$, a computation using gp shows that the resultant
of $f_{\sqrt{2},\sqrt{13}}(X_0,1,1)$ and $\frac{\partial
f_{\sqrt{2},\sqrt{13}}}{\partial X_0}(X_0,1,1)$ is nonzero.)

For the rest of this section, fix $f:=f_{\sqrt{u},\sqrt{v}}$ and $L$
as above and let $\ol{\Q}$ be an algebraic closure of $\Q$
containing $L$.

\begin{lem}\label{g18eq}   Following the above notation,
$f=0$ gives the equation of a smooth plane curve $X$ over $\ol{\Q}$
with $\Aut(X)$ given by $\G_{18}$. Furthermore, the field of moduli
of $X$  is equal to $K$ and any isomorphism $X\to\el{\sigma}{X}$
with $\sigma\in\Gal(\ol{\Q}/K)$ is given by an element of $\G_{72}$.
\end{lem}
\begin{proof}
We now show that $f=0$ gives the equation of a smooth plane curve.
We first show that if $g\in\ol{\Q}[X_0,X_1,X_2]$ is a form of
positive degree and if $g\mid f$ then $g^2\nmid f$. Suppose there
exists a form $g$ of positive degree such that $g^2\mid f$. Then
$g\mid f_{X_0}$ and $g\mid f_{X_1}$.  By B\'{e}zout's Theorem,
$g(X_0,X_1,X_2)=0$ and $X_2=0$ have a nonempty intersection.  We
have
\[f_{X_0}(X_0,X_1,0)\]
 \[=6X_0^2((\alpha_1\sqrt{u} + \alpha_2\sqrt{v} + \alpha_3\sqrt{uv} - 1/12)X_0^3
  + (\alpha_1 \sqrt{u}+ \alpha_2\sqrt{v} + \alpha_3\sqrt{uv} +5/12)X_1^3)\] and
\[f_{X_1}(X_0,X_1,0)\]
\[=6X_1^2((\alpha_1\sqrt{u} + \alpha_2\sqrt{v} + \alpha_3\sqrt{uv} -
1/12)X_1^3 + (\alpha_1 \sqrt{u}+ \alpha_2\sqrt{v} +
\alpha_3\sqrt{uv} + 5/12)X_0^3).\]

The zeros of $f_{X_0}(X_0,X_1,0)$ are
\[\left\{[0\colon 1\colon 0],[1\colon \lambda\colon 0]\mid \lambda^3
=-\frac{\alpha_1 \sqrt{u}+ \alpha_2\sqrt{v} + \alpha_3\sqrt{uv}
-1/12}{\alpha_1\sqrt{u} + \alpha_2\sqrt{v} + \alpha_3\sqrt{uv} +
5/12}\right\}\] and the zeros of $f_{X_1}(X_0,X_1,0)$ are
\[\left\{[1\colon 0\colon 0],[\lambda\colon 1 \colon 0]\mid \lambda^3
=-\frac{\alpha_1 \sqrt{u}+ \alpha_2\sqrt{v} + \alpha_3\sqrt{uv}
-1/12}{\alpha_1\sqrt{u} + \alpha_2\sqrt{v} + \alpha_3\sqrt{uv} +
5/12}\right\}.\]

If $f_{X_0}(X_0,X_1,0)$ and $f_{X_1}(X_0,X_1,0)$ share a common
zero,
 it must be $[1\colon\lambda\colon 0]$ where $1/\lambda^3=\lambda^3$.
   So $\lambda^3 = \pm 1$.
If we let $s =\alpha_1 \sqrt{u} + \alpha_2 \sqrt{v} + \alpha_3
\sqrt{uv}$,
 then $\pm 1 = - \frac{s-1/12}{s+5/12}$.  So $s \in \Q$, contradicting the
definition of $s$.

We will say that $f$ is of type $(d_1,\ldots, d_n)$ if $f$ can be
factored into irreducible forms of positive degrees $f_1,\ldots,
f_n$ of degrees $d_1,\ldots, d_n$ respectively, where $d_i\le d_j$
if $i\le j$. Suppose $f$ is of type $(d_1,\ldots, d_n)$.  Note that
since the degree of $f$ is $6$, $n\in\{1,2,3,6\}$. Write
$f=\prod_{i=1}^n f_i$ where $f_i$ has degree $d_i$. Let
$g_i:=\frac{(d_i-1)(d_i-2)}{2}$ be the arithmetic genus of the curve
given by $f_i=0$.  We will make use of the fact that a geometrically
integral curve of arithmetic genus $g$ has at most $g$
singularities. Let $P$ be in $\mathbb{P}^2(\ol{\Q})$. Then
$f(P)=f_{X_0}(P)=f_{X_1}(P)=f_{X_2}(P)=0$ if and only if $P$ is a
singular point of a curve given by $f_i=0$ for some $i$ or if $P$ is
an intersection point of two or more of the $f_i$. Using
B\'{e}zout's Theorem we obtain a bound $S(f)$ for the number of
points $P\in\mathbb{P}^2(\ol{\Q})$ such that
$f(P)=f_{X_0}(P)=f_{X_1}(P)=f_{X_2}(P)=0$, namely
\[S(f):=\sum_{i=1}^n g_i +\sum_{i<j}d_id_j. \]
Write
\[G(f):=\sum_{i=1}^n g_i\]
and
\[I(f):=\sum_{i<j}d_id_j. \] Checking each partition of deg$(f)=6$ shows
 that $G(f)\in\{0, 1, 2, 3, 6, 10\}$. Another computation shows that
if $G(f)=0$ then $I(f)\le 15$, if $G(f)=1$ then $I(f)\le 12$, if
$G(f)=2$ then $I(f)\le 9$, if $G(f)=3$ then $I(f)\le 9$, if $G(f)=6$
then $I(f)\le 5$, and if $G(f)=10$ then $I(f)=0$.  It follows that
$S(f)\le 15$.

Let $P$ be in $\mathbb{P}^2(\ol{\Q})$ and suppose that
\[f(P)=f_{X_0}(P)=f_{X_1}(P)=f_{X_2}(P)=0.\] It is easy to see that $f$ is $\G_{18}$-invariant.
So for any $Q\in \mf{O}_{18}(P)$ we must have
\[f(Q)=f_{X_0}(Q)=f_{X_1}(Q)=f_{X_2}(Q)=0,\] where $\mf{O}_{18}(P)$
is the orbit of $P$ under the action of $\G_{18}$.
    By Lemma~\ref{g18orbits}, $\mf{O}_{18}(P)$
must have size $3$, or $9$.  Since $f(X_0,1,1)$ is squarefree, by
Lemma~\ref{g18orbits}, the orbit of $P$ under the action of
$\G_{18}$ must have size $3$.  One can verify directly that $f(Q)\ne
0$ if $Q\in\{P_3^i\}_{1\le i\le 4}$. So no such $P$ exists.
 It follows that $f=0$ gives the equation of a smooth plane
curve $X$.

Since the degree of $f$ is $6$, the genus of $X$ is $10$, so
$\Aut(X)$ is finite.  Since $X$ is a smooth plane curve of degree
larger than $3$, by Theorem~\ref{planeisom}, $\Aut(X)$ is
$\PGL_3(\ol{\Q})$-conjugate to one of the groups listed in
Lemma~\ref{gps}.  We now show that $\Aut(X)=\G_{18}$.  It can be
deduced from Lemma~\ref{gps} that $\Aut(X)$ is
$\PGL_3(\ol{\Q})$-conjugate to either a group of type $\mf{D}$, a
subgroup of $\G_{216}$, both a group of type $\mf{D}$ and a subgroup
of $\G_{216}$, or a subgroup of $\G_{360}$.

  For $1\le i\le 3$, define
\[f^{(i)}_{\sqrt{u},\sqrt{v}}:=c_{\psi^2}\psi^2-6\omega^{i-1} c_{\phi\psi}\phi\psi-18\omega^{2(i-1)} c_{\phi^2}\phi^2+\chi.\]
Then, $f=f^{(1)}_{\sqrt{u},\sqrt{v}}$ and we have
$[f^{(i)}_{\sqrt{u},\sqrt{v}}]=([f^{(1)}_{\sqrt{u},\sqrt{v}}])U^i$,
where $U\in\G_{216}$ is given in Lemma~\ref{gps}.  Recall that
$\G_{216}=\la\G_{72},U\ra$.  A computation shows that
\[\mf{O}_{72}([f^{(i)}_{\sqrt{u},\sqrt{v}}])=\{[f^{(i)}_{\pm\sqrt{u},\pm\sqrt{v}}]\},\]
for $1\le i\le 3$, where
$\mf{O}_{72}([f^{(i)}_{\sqrt{u},\sqrt{v}}])$ is the orbit of
$[f^{(i)}_{\sqrt{u},\sqrt{v}}]$ under the action of $\G_{72}$, and
that $\la U\ra$ cyclically permutes these three orbits. Note that
\[\mf{O}_{216}([f])=\bigsqcup_{i=1}^3\mf{O}_{72}([f^{(i)}_{\sqrt{u},\sqrt{v}}]),\]
where $\mf{O}_{216}([f])$ is the orbit of $[f]$ under the action of
$\G_{216}$.  By looking at the coefficients of the
$f^{(i)}_{\pm\sqrt{u},\pm\sqrt{v}}$, we see that for all
$\sigma\in\Gal(\ol{\Q}/\Q)$ and all $[g]\in\mf{O}_{216}([f])$, we
have $[g^{\sigma}]\in\mf{O}_{216}([f])$ if and only if
$\sigma(\omega)=\omega$.  We also see by looking at the coefficients
of the $f^{(i)}_{\pm\sqrt{u},\pm\sqrt{v}}$ that for all
 $[g]\in\mf{O}_{216}([f])$, we
have $[h]\in\mf{O}_{72}([g])$ if and only if $[h]=[g^{\sigma}]$ for
some $\sigma\in\Gal(\ol{\Q}/K)$.

Suppose $M\in\PGL_3(\ol{\Q})$ and suppose that $M\Aut(X) M^{-1}$ is
a subgroup of $\G_{216}$.  Then since $\G_{18} \subseteq\Aut(X)$, we
have $M\G_{18}M^{-1} \subseteq M \Aut(X) M^{-1} \subseteq \G_{216}$.
  According to Magma, $\G_{216}$ has two
conjugacy classes of subgroups of order $18$, one consisting of the
group $\G_{18}$ and the other containing $12$ conjugate subgroups of
order $18$.  It is easy to see that the subgroup $\la R,
U\ra\subset\G_{216}$ has order $18$.  Since $\la R, U\ra$ is
$\PGL_3(\ol{\Q})$-conjugate to an intransitive subgroup of
$\PGL_3(\ol{\Q})$ and $\G_{18}$ is not, $\la R, U\ra$ and $\G_{18}$
are not $\PGL_3(\ol{\Q})$-conjugate.  It follows that $M\in
N(\G_{18})$.  By Lemma~\ref{normalizergps}, $N(\G_{18})=\G_{216}$.
So $\Aut(X)$ must be a subgroup of $\G_{216}$.  By the previous
paragraph, the subgroup of $\G_{216}$ consisting of elements that
stabilize $[f]$ equals $\G_{18}$.
    Thus if $\Aut(X)$ is
$\PGL_3(\ol{\Q})$-conjugate to a subgroup of $\G_{216}$ then
$\Aut(X)=\G_{18}$.

  Suppose that $\Aut(X)$ is $\PGL_3(\ol{\Q})$-conjugate to a group of type
  $\mf{D}$.  We will call the union of $3$ lines in
$\mathbb{P}^2$ a triangle. By \S113 of~\cite{Miller:1}, a group of
type $\mf{D}$ containing $\G_{18}$ but not equal to $\G_{18}$ fixes
exactly one triangle in $\mathbb{P}^2$: the triangle given by
$X_0X_1X_2=0$.  By \S113 of~\cite{Miller:1}, the group $\G_{18}$
fixes four triangles in $\mathbb{P}^2$:
\[X_0X_1X_2=0,\]
\[(X_0+X_1+\theta X_2)(X_0+\omega X_1+\omega^2\theta
X_2)(X_0+\omega^2X_1+\omega\theta X_2)=0,\] where
$\theta\in\{1,\omega, \omega^2\}$.  Since $\Aut(X)$ is
$\PGL_3(\ol{\Q})$-conjugate to a group of type $\mf{D}$ and since
$\G_{18}\subseteq\Aut(X)$, $\Aut(X)$ must fix at least one of the
four $\G_{18}$-invariant triangles.  By \S115 of~\cite{Miller:1},
the four invariant triangles of $\G_{18}$ are permuted transitively
by $\G_{216}$.  So for some $A\in\G_{216}$, the group $A^{-1}\Aut(X)
A$ must fix the triangle $X_0X_1X_2=0$. Following the notation of
Lemma~\ref{gps}, it can easily be deduced that a finite subgroup
$\G\subset\PGL_3(\ol{\Q})$ that fixes the triangle $X_0X_1X_2=0$ is
a group of type $\mf{D}$ if and only if $T, R\in\G$.  By
Lemma~\ref{normalizergps}, $N(\G_{18})=\G_{216}$. So
$\G_{18}\subseteq A^{-1}\Aut(X) A$.  In particular we must have $T,
R\in A^{-1}\Aut(X) A$.  So $A^{-1}\Aut(X) A$ is a group of type
$\mf{D}$.  Choose $f_A\in\ol{\Q}[X_0,X_1,X_2]$ so that
$[f_A]=([f])A\in\mf{O}_{216}([f])$ and let $Y$ be the curve given by
$f_A=0$.  We have just seen that $[g]\in\mf{O}_{72}([f_A])$ if and
only if and if $[g]=[f_A^{\sigma}]$ for some
$\sigma\in\Gal(\ol{\Q}/K)$.  So every element of $\G_{72}$ gives an
isomorphism $Y\to\el{\sigma}{Y}$ for some
$\sigma\in\Gal(\ol{\Q}/K)$.  So by Lemma~\ref{conjcurve}, we must
have $\G_{72}\subseteq N(A^{-1}\Aut(X)A)$.  By
Lemma~\ref{normalizergps}, a subgroup $\G$ of type $\mf{D}$ has
$\G_{72}\subseteq N(\G)$ if and only if $\G=\G_{18}$. Hence
$A^{-1}\Aut(X) A=\G_{18}$.  Since $A\in N(\G_{18})$, we must have
$\Aut(X)=\G_{18}$.  Thus if $\Aut(X)$ is $\PGL_3(\ol{\Q})$-conjugate
to a group of type $\mf{D}$, then $\Aut(X)=\G_{18}$.

Lastly, let $M\in\PGL_3(\ol{\Q})$ and suppose that $M\Aut(X) M^{-1}$
is a subgroup of $\G_{360}$.   According to Magma, the maximal
subgroups of $\G_{360}$ are isomorphic to $S_4$, $A_5$, or
$\G_{36}$.  Of these three groups, only $\G_{36}$ has a subgroup
isomorphic to $\G_{18}$, and such a subgroup of $\G_{36}$ must be
normal in $\G_{36}$ since it is of index $2$.  So if
$\Aut(X)\ne\G_{18}$, we must have
$M\G_{18}M^{-1}\subsetneq(M\Aut(X)M^{-1}\cap N(M\G_{18}M^{-1})$. It
follows that $\G_{18}\subsetneq(\Aut(X) \cap N(\G_{18}))$.   But we
have just seen that $([f])A\ne [f]$ for any $A\in
N(\G_{18})-\G_{18}$. Therefore we must have $\Aut(X)=\G_{18}$.

  By Lemma~\ref{conjcurve}, any
isomorphism $X\to\el{\sigma}{X}$ with $\sigma\in\Gal(\ol{\Q}/\Q)$ is
given by an element of $N(\G_{18})=\G_{216}$.  We have seen that for
$\sigma\in\Gal(\ol{\Q}/\Q)$, we have already
 $[f^{\sigma}]=([f])A$ for some $A\in\G_{216}$ if and only if $\sigma|_K=Id$ and
 $A\in\G_{72}$.
\end{proof}

\begin{prop} Let $f$ be as above and let $X$ be
the smooth plane curve over $\ol{\Q}$ given by $f=0$.  Then $X$ is
not definable over its field of moduli.
\end{prop}
\begin{proof}  By Lemma~\ref{g18eq}, the field of moduli of $X$
 is $K$.  Also by Lemma~\ref{g18eq}, any
isomorphism $X\to\el{\sigma}{X}$, with $\sigma\in\Gal(\ol{\Q}/K)$,
is given by an element of $\G_{72}$.

Suppose that $X$ is definable over $K$.  Recall that
$L=K(\sqrt{u},\sqrt{v})$. Then there exists a finite Galois
extension field $F$ of $K$ and a collection of isomorphisms
$\mathscr{C}:=\{\varphi_{\sigma}\}_{\sigma\in\Gal(F/K)}$ that
satisfy Weil's cocycle condition.  Without loss of generality we may
assume that $L\subseteq F$.   Note that if $\sigma, \tau
\in\Gal(F/K)$ and $\varphi_{\sigma}, \varphi_{\tau}\in\mathscr{C}$,
then we have
\[\varphi_{\sigma^{-1}}=\varphi_{\sigma}^{-1}\] and \[\varphi_{\tau}\varphi_{\sigma}=\varphi_{\sigma\tau},\]
since $\Gal(F/K)$ is finite and since $\G_{72}\subset\PGL_3(K)$.  If
$\sigma\in\Gal(F/K)$, let $M_{\sigma}\in\G_{72}$ be the element of
$\PGL_3(K)$ that gives the corresponding isomorphism
$\varphi_{\sigma}\in\mathscr{C}$. Then we have a homomorphism
$\Psi\colon\Gal(F/K)\to\G_{72}$ defined by
\[\sigma\mapsto M_{\sigma}^{-1}.\]

 Lemma~\ref{g18eq} shows that $\Psi^{-1}(\G_{18}) = \Gal(F/L)$.
Thus $\Psi$ induces an injection $\Gal(L/K)\to \G_{72}/\G_{18}$,
which because of orders must be an isomorphism.  According to Magma,
$\G_{72}$ has a unique subgroup of order $9$ which, by
Lemma~\ref{gps}, must be $\G_9$.  Also according to Magma, $\G_{72}$
has three conjugacy classes of elements of order $4$, each
consisting of $18$ elements. Following the notation of
Lemma~\ref{gps}, since $\G_{18}$ has no elements of order $4$ and is
normal in $\G_{72}$, since $U\not\in\G_{72}$, and since $V$ has
order $4$ in $\G_{72}$, these conjugacy classes are given by
$V\G_{18}$, $U^{-1}VU\G_{18}$, and $UVU^{-1}\G_{18}$.  So any
subgroup of $\G_{72}$ which surjects onto $\G_{72}/\G_{18}$ must
contain an element from each conjugacy class of order $4$. According
to Magma, the  proper subgroups of $\G_{72}$ which contain an
element from each conjugacy class of order $4$ are maximal and of
order $8$. Thus, any subgroup of $\G_{72}$ which surjects onto
$\G_{72}/\G_{18}$ must also surject onto $\G_{72}/\G_9$.  In
particular, the image of $\Psi$ must surject onto $\G_{72}/\G_{9}$.
Let
 $L'$ be the subfield of $F$ fixed by the normal subgroup $\Psi^{-1}(\G_9)$
 of $\Gamma$.  Then $\Gal(L'/K)$ is isomorphic to $\G_{72}/\G_9$,
   which is a quaternion group.  Since  $\G_9 \subset \G_{18}$ we have $\Psi^{-1}(\G_9) \subset
\Psi^{-1}(\G_{18})$,
   so $L'$ contains $L$.  This is a contradiction, since by assumption $L$ cannot be
embedded into a quaternion extension of $K$. Therefore, no such $F$
exists and $X$ is not definable over $K$.
\end{proof}

\section{Plane curves with automorphism groups given by $\G_{36}$}

\begin{lem}\label{g36orbits}  Let $F$ be an algebraically closed
field of characteristic $0$ or of characteristic greater than $3$.
Let $P$ be in $\mathbb{P}^2(F)$.  Then the orbit $\mf{O}(P)$ of $P$
under the action of $\G_{36}$ has $36$ points unless $P$ is in the
orbit of \[P_6:=[1\colon 0\colon 0],\] \[P_9:=[0\colon 1\colon
-1],\]
\[Q_9:=[1-\sqrt{3}\colon 1\colon 1],\]
\[Q'_9:=[1+\sqrt{3}\colon 1\colon 1],\] or
\[P_{18, \delta}:=[\delta\colon 1\colon 1],\]
with $\delta\not\in\{1, 1\pm\sqrt{3}\}$.  We have $|\mf{O}(P_6)|=6$,
$|\mf{O}(P_9)|=9$, $|\mf{O}(Q_9)|=9$, $|\mf{O}(Q'_9)|=9$, and
$|\mf{O}(P_{18, \delta})|=18$.
\end{lem}
\begin{proof}  We follow the notation of Lemma~\ref{gps}.
Let $P$ be in $\mathbb{P}^2(F)$ and suppose that the orbit of $P$
under the action of $\G_{36}$ has less than $36$ points. It follows
that for some $g\in \G_{36}-\{Id\}$, $g(P)=P$.  If $h$ is in
$\G_{36}$, then $(h^{-1}gh)(h^{-1}(P))=h^{-1}(P)$.  So any element
in $\G_{36}$ which is conjugate to $g$ fixes a point in the orbit of
$P$ under the action of $\G_{36}$.  It is easily shown that any non
identity element of $\G_{36}$ is $\G_{36}$-conjugate to one of $S$,
$S^2$, $V$, $V^2$, or $V^3$.

A computation shows the following.
 If $S(P)=P$ or $S^2(P)=P$ then \[P\in\{[1\colon 0\colon 0], [0\colon 1\colon
0], [0\colon 0\colon 1]\}.\] If  $V(P)=P$ or $V^3(P)=P$ then
\[P\in\{[0\colon 1\colon -1], [1-\sqrt{3}\colon 1\colon 1], [1+\sqrt{3}\colon 1\colon 1]\}.\] If $V^2(P)=P$ then
$P=[\delta \colon \gamma\colon \gamma]$ where $\delta$ and $\gamma$
are not both zero.   Note that if $\gamma=0$ then $P=P_6$,  if
$\gamma\ne 0$ and $\delta=\gamma$ then $P\in\mf{O}(P_6)$, and if
$\gamma\ne 0$ and $\frac{\delta}{\gamma} =1\pm\sqrt{3}$ then
$P\in\{Q_9, Q'_9\}$.

A simple computation shows that each of these points is in the orbit
of one of the points listed in the lemma.  Another computation using
the stabilizers of each point proves the last statement of the
lemma.
\end{proof}
 We define the following:
 \begin{align*}
  \phi&:= X_0X_1X_2,\\
 \psi&:=X_0^3+X_1^3+X_2^3,\\
 \chi&:=X_0^3X_1^3 + X_1^3 X_2^3 + X_2^3 X_0^3.\\
\end{align*}\vspace*{-.76in}
\newline Let $a$ be in $\C$ and let
\[f_a(X_0,X_1,X_2):=\chi - 18a\phi^2+ (a - 1/12)\psi^2 - 6a\psi\phi.\]

\begin{lem}\label{g36eq} Following the above notation, write $f=f_a$,
 where $a:=i\alpha\in\C$ with $\alpha\in \mathbb R^{\times}$.
Then $f=0$ gives the equation of a smooth  plane curve $X$ defined
over $\C$ with $\Aut(X)$ given by $\G_{36}$.
\end{lem}
\begin{proof} We show that $f=0$ gives the equation of a smooth curve $X$.
 We first now show that if $g\in\C[X_0,X_1,X_2]$ is a form of
positive degree and if $g\mid f$ then $g^2\nmid f$.  Suppose there
exists a form $g$ of positive degree such that $g^2\mid f$.  Then
$g\mid f_{X_0}$ and $g\mid f_{X_1}$.  By B\'{e}zout's Theorem,
$g(X_0,X_1,X_2)=0$ and $X_2=0$ have a nonempty intersection.  We
have
\[f_{X_0}(X_0,X_1,0)=X_0^2((6a-1/2)X_0^3+(6a + 5/2)X_1^3) \]
and
\[f_{X_1}(X_0,X_1,0)=X_1^2((6a + 5/2)X_0^3+(6a-1/2)X_1^3) \]
The zeros of $f_{X_1}(X_0,X_1,0)$ are
\[\left\{[1\colon 0\colon 0],[\beta\colon 1\colon 0]\mid \beta^3=\frac{6a+5/2}{-6a+1/2}\right\}\]
and the zeros of $f_{X_0}(X_0,X_1,0)$ are
\[\left\{[0\colon 1\colon 0],[1\colon \beta \colon 0]\mid \beta^3=\frac{6a+5/2}{-6a+1/2}
\right\}.\]  Since $f_{X_0}(X_0,X_1,0)$ and $f_{X_1}(X_0,X_1,0)$
share a common zero we must have $1/\beta^3=\beta^3$.  So
$\beta^3=\pm 1$.  Then $\pm 1=\frac{6a+5/2}{-6a+1/2}$.  This implies
$a\in\Q$ which contradicts our choice of $a$.

Since $f$ is of degree $6$, by an argument identical to one given in
the proof of Lemma~\ref{g18eq}, a bound for the number of points
$P\in\mathbb{P}^2(\C)$ such that
\[f(P)=f_{X_0}(P)=f_{X_1}(P)=f_{X_2}(P)=0\] is $15$.
Let $P$ be in $\mathbb{P}^2(\C)$ and suppose that
\[f(P)=f_{X_0}(P)=f_{X_1}(P)=f_{X_2}(P)=0.\]  A simple
computation shows that $f$ is $\G_{36}$-invariant.  So for any $Q\in
\mf{O}(P)$ we must have
\[f(Q)=f_{X_0}(Q)=f_{X_1}(Q)=f_{X_2}(Q)=0.\]
    By Lemma~\ref{g36orbits}, the orbit of $P$ under the action of $\G_{36}$
must have size $9$, or $6$. One can verify directly  that $f(Q)\ne
0$ if $Q\in\{P_6, P_9\}$ and that $f_{X_0}(Q)\ne 0$ if $Q\in\{Q_9,
Q'_9\}$.
   It follows that $f=0$ gives the equation of a smooth
plane curve $X$.

Since the degree of $f$ is $6$, the genus of $X$ is $10$, so
$\Aut(X)$ is finite.  Since $X$ is a smooth plane curve of degree
larger than $3$, by Theorem~\ref{planeisom}, $\Aut(X)$ is
$\PGL_3(\ol{\C})$-conjugate to one of the groups listed in
Lemma~\ref{gps}.   We now show that $\Aut(X)=\G_{36}$.  We use the
notation of Lemma~\ref{gps}. According to Magma, any subgroup of
$\G_{360}$ of order $36$ is maximal.  Since
$\G_{36}\subseteq\Aut(X)$, Lemma~2.3.7 implies that $\Aut(X)$ equals
$\G_{36}$, $M\G_{72}M^{-1}$,
 $M\G_{216}M^{-1}$, or $M\G_{360}M^{-1}$
 for some $M\in\PGL_3(\C)$.

 Suppose $\Aut(X) = M\G_{216}M^{-1}$.  Since $\G_{18}$ is contained in $\Aut(X)$,
 $M^{-1}\G_{18}M$ is contained in $M^{-1}\Aut(X)M =\G_{216}$.
Magma shows that the only subgroup of $\G_{216}$ isomorphic to
$\G_{18}$ is $\G_{18}$ itself,
 so $M^{-1}\G_{18}M =\G_{18}$.  Thus $M\in N(\G_{18}) =\G_{216}$,
 so $\Aut(X) =M\G_{216}M^{-1}=\G_{216}$.   Similarly, if $\Aut(X)=M\G_{72}M^{-1}$,
  the same argument shows that $M\in\G_{216} = N(\G_{72})$,
 so $\Aut(X) =\G_{72}$.  Following the notation of Lemma~\ref{gps},
 in both cases we must have $UVU^{-1}\in\Aut(X)$.  A computation
 shows that $([f])UVU^{-1}\ne[f]$.  So we cannot have
$\Aut(X)=\G_{72}$, nor can we have $\Aut(X)=\G_{216}$.

Lastly, suppose that $\Aut(X)$ is $\PGL_3(\C)$-conjugate to
$\G_{360}$.  Recall that $\G_{36}=\la S, T, V\ra$ and that $V^2=R$.
Define $\G_{360}':=\la R, T, E_1', S(E_2')S^{-1}\ra$, where
\[E_1':=\begin{bmatrix}
1&0&0\\
0&-1&0\\
0&0&-1\\
\end{bmatrix} \mbox{ and }
E_2':=\begin{bmatrix}
-1&\mu_2&\mu_1\\
\mu_2&\mu_1&-1\\
\mu_1&-1&\mu_2\\
\end{bmatrix}\] with $\mu_1=1/2(-1+\sqrt{5})$ and
$\mu_2=1/2(-1-\sqrt{5})$.  By \S123 of~\cite{Miller:1},
$S^{-1}\G_{360}'S$ is $\PGL_3(\C)$-conjugate to $\G_{360}$.
According to Magma, $\G_{360}'$ has two conjugacy classes of
subgroups of order $6$.  Note that the subgroups $\mf{H}'_6:=\la S,
R\ra\subset\G_{36}$ and $\mf{H}_6:=\la T, R\ra\subset\G_{36}$ both
have order $6$ and are not $\PGL_3(\C)$-conjugate, since $\mf{H}'_6$
is $\PGL_3(\C)$-conjugate to an intransitive subgroup and $\mf{H}_6$
is not.  Choose $M'\in\PGL_3(\C)$ so that
$M'\Aut(X)(M')^{-1}=\G_{360}'$.   The subgroups
$M'\mf{H}_6(M')^{-1}$ and $M'\mf{H}_6'(M')^{-1}$ of $\G_{360}'$ must
be in distinct conjugacy classes under $\G_{360}'$.  Since
$\mf{H}_6\subset\G_{360}'$ and since $\mf{H}_6$ is not
$\PGL_3(\C)$-conjugate to $M'\mf{H}_6'(M')^{-1}$, $\mf{H}_6$ is in
the same conjugacy class as $M'\mf{H}_6(M')^{-1}$ under $\G_{360}'$.
So without loss of generality we may assume that $M'\in
N(\mf{H}_6)$. By Lemma~\ref{normalizergps}(c),
$N(\mf{H}_6)=\mf{H}_6$.  Thus $M'\in\G_{360}'$.  So
$\Aut(X)=\G_{360}'$.  A computation shows that
 $([f])E_1'\ne [f]$, so we get a contradiction.  Therefore, we must
 have $\Aut(X)=\G_{36}$.
\end{proof}

\begin{prop} Let $f_a$ be as above where $a=\beta i\in\C$ with $\beta\in\R^{\times}$
 and let $X$ be the smooth plane
curve defined over $\C$ given by $f_a=0$.  Then the field of moduli
of $X$ relative to the extension $\C/\R$ is $\R$ and is not a field
of definition for $X$.
\end{prop}
\begin{proof}  We follow the notation of Lemma~\ref{gps}.
Let $\Gal(\C/\R):=\la\sigma\ra$.  By Lemmas~\ref{isnorm}
and~\ref{normalizergps}, any isomorphism $X\to\el{\sigma}{X}$ is
given by and element of \[N(\G_{36})=\G_{72}=\la\G_{36},
UVU^{-1}\ra.\] A computation shows that
\[([f_a])UVU^{-1}=[f_{-a}].\]  So $UVU^{-1}$ gives an isomorphism
$X\to\el{\sigma}{X}$.  A computation using gp shows that for all $M$
in the nontrivial coset of $\G_{72}/\G_{36}$, $M^{\sigma}M\ne Id$.
Therefore, Weil's cocycle condition of Theorem~\ref{weil} does not
hold.  It follows that $X$ cannot be defined over $\R$.
\end{proof}